\documentclass [reqno]{amsart}
\usepackage {mathptmx}
\usepackage {amsmath,amssymb}
\usepackage {color,graphicx}
\usepackage [utf8]{inputenc}
\usepackage [dvips,centering,includehead,width=14cm,height=22.5cm]{geometry}

\newcommand {\ignore}[1]{}

\newcommand {\df}[1]{\textsl {#1}}

\newcommand {\gen}[1]{\langle #1 \rangle}
\newcommand {\M}[1]{M_{\text {\rm (#1)}}}
\newcommand {\n}[1]{n_{\text {\rm (#1)}}}
\newcommand {\Mrs}{M_{(r,s)}}
\newcommand {\Nrs}{N_{(r,s)}}
\newcommand {\ori}{{\text {\rm or}}}
\newcommand {\even}{{\text {\rm even}}}
\newcommand {\SM}{{\mathbb S}}

\DeclareMathOperator {\val}{val}
\DeclareMathOperator {\sign}{sign}

\DeclareMathOperator {\ev}{ev}
\DeclareMathOperator {\ft}{ft}
\DeclareMathOperator {\mult}{mult}

\newcommand {\isom}{\cong}

\sloppypar

\renewenvironment {enumerate}%
  {\rule{1mm}{0mm}\begin {oldenumerate}%
    \parskip1ex plus0.5ex \itemsep 0mm \parindent 0mm}%
  {\end {oldenumerate}}

\renewenvironment {itemize}%
  {\rule{1mm}{0mm}\begin {olditemize}%
    \parskip1ex plus0.5ex \itemsep 0mm \parindent 0mm}%
  {\end {olditemize}}

\newlength {\sidetextwidth}
\newlength {\sidepicwidth}
\newsavebox {\sidepicbox}
\newenvironment {sidepic}[1]{%
  \par%
  \def\message##1{}%
  \hbadness=10000%
  \setlength{\sidetextwidth}{\textwidth}%
  \savebox{\sidepicbox}{\input{pics/#1}}%
  \settowidth{\sidepicwidth}{\usebox{\sidepicbox}}%
  \addtolength{\sidetextwidth}{-6mm}%
  \addtolength{\sidetextwidth}{-\sidepicwidth}%
  \begin {minipage}{\sidetextwidth}}%
  {\end {minipage}%
  \hfuzz=10000pt%
  \hspace {4mm}%
  \begin {minipage}{\sidepicwidth}\usebox{\sidepicbox}\end {minipage}%
  \hfuzz=0pt%
  \hbadness=1000%
  }

\theoremstyle {plain}
\newtheorem {theorem}{Theorem}[section]
\newtheorem {proposition}[theorem]{Proposition}
\newtheorem {lemma}[theorem]{Lemma}

\newtheorem {corollary}[theorem]{Corollary}
\theoremstyle {definition}
\newtheorem {definition}[theorem]{Definition}
\theoremstyle {remark}
\newtheorem {remark}[theorem]{Remark}
\newtheorem {example}[theorem]{Example}

\newtheorem {convention}[theorem]{Convention}
\newtheorem {notation}[theorem]{Notation}
\newtheorem {algo}[theorem]{Algorithm}

\newcommand{\N}{\mathbb{N}}
\newcommand{\NN}{\mathbb{N}}
\newcommand{\ZZ}{\mathbb{Z}}

\newcommand{\RR}{\mathbb{R}}

\newcommand{\CC}{\mathbb{C}}
\newcommand{\PP}{\mathbb{P}}
\newcommand{\calM}{\mathcal{M}}
\newcommand{\calP}{\mathcal{P}}

\title {Broccoli curves and the tropical invariance of Welschinger numbers}
\author {Andreas Gathmann \and Hannah Markwig \and Franziska Schroeter}
\address {Andreas Gathmann, Fachbereich Mathematik, Technische Universit\"at
  Kaiserslautern, Postfach 3049, 67653 Kaiserslautern, Germany}
  \email{andreas@mathematik.uni-kl.de}
\address {Hannah Markwig, Fachrichtung Mathematik, Universit\"at des
  Saarlandes, Postfach 151150, 66041 Saarbr\"ucken, Germany}
  \email{hannah@math.uni-sb.de}
\address {Franziska Schroeter, Cluster of Excellence $M^2CI$, Fachrichtung Mathematik,
  Universit\"at des Saarlandes, Postfach 151150, 66041 Saarbr\"ucken, Germany}
  \email{schroeter@math.uni-sb.de}

\thanks {\emph {2010 Mathematics Subject Classification:} 14T05, 14N10}
\keywords {Tropical geometry, enumerative geometry, Welschinger numbers}

\begin {document}

  \begin {abstract}
    In this paper we introduce broccoli curves, certain plane tropical curves
    of genus zero related to real algebraic curves. The numbers of these
    broccoli curves through given points are independent of the chosen points
    --- for arbitrary choices of the directions of the ends of the curves,
    possibly with higher weights, and also if some of the ends are fixed. In
    the toric Del Pezzo case we show that these broccoli invariants are equal
    to the Welschinger invariants (with real and complex conjugate point
    conditions), thus providing a proof of the independence of Welschinger
    invariants of the point conditions within tropical geometry. The general
    case gives rise to a tropical Caporaso-Harris formula for broccoli curves
    which suffices to compute all Welschinger invariants of the plane.
  \end {abstract}

  \maketitle

  \section {Introduction} \label {sec-intro}

\subsection {Background on tropical Welschinger numbers}

Welschinger invariants of real toric unnodal Del Pezzo surfaces count real
rational curves, weighted with $\pm 1$ depending on the nodes of the curve,
belonging to an ample linear system $D$ and passing through a generic
conjugation invariant set $\calP$ of $-K_{\Sigma}\cdot D-1$ points. It was
shown in \cite{Wel03} and \cite{Wel05} that these numbers depend only on the
number of real points in $ \calP $, i.e.\ are invariant under movements of the
points in $ \calP $. They can be thought of as real analogues of the numbers of
complex rational curves belonging to a fixed linear system and satisfying point
conditions, which in the case of $\PP^2$ are the genus-$0$ Gromov-Witten
invariants.

By Mikhalkin's Correspondence Theorem \cite{Mik03}, Gromov-Witten invariants
of the plane (resp.\, the complex enumerative numbers of other toric surfaces)
can be determined via tropical geometry, by counting tropical curves of a fixed
degree and satisfying point conditions. Each tropical curve has to be counted
with a ``complex multiplicity'' which reflects how many complex curves map to
it under tropicalization.

Welschinger invariants can be computed via tropical geometry in a similar way:
one can define a certain count of tropical curves and prove a Correspondence
Theorem stating that this tropical count equals the Welschinger invariant. For
the case when $\calP$ consists of only real points, such a Correspondence
Theorem is proved in \cite{Mik03}, the general case is proved in \cite{Shu06}.

If $\calP$ consists of only real points, the tropical curves we have to count
to get Welschinger invariants are exactly the same as the ones we need to count
to determine complex enumerative numbers --- we just have to count them with a
different, ``real'' multiplicity. The lattice path algorithm of \cite{Mik03}
enumerates the tropical curves we have to count. If $\calP$ also contains pairs
of complex conjugate points, we have to count tropical curves satisfying some
more special conditions. The lattice path algorithm is generalized in
\cite{Shu06} to an algorithm that computes the corresponding Welschinger
invariants.

It follows from the Correspondence Theorem and the fact that Welschinger
invariants are independent of the point conditions that the corresponding
tropical count is also invariant, i.e.\ does not depend on the position of the
points that we require the tropical curves to pass through. 

Still, it is interesting to find an argument within tropical geometry that
proves the invariance of the tropical numbers. For the case when $\calP$
consists of only real points, such a statement follows easily since the
corresponding tropical count can be shown to be locally invariant, i.e.\
invariant around a codimension-1 cone of the corresponding moduli space of
curves. In addition, such a codimension-1 cone is specified by a $4$-valent
vertex of a tropical curve, and it is sufficient to consider the curves locally
around this $4$-valent vertex. This tropical invariance statement was proved in
\cite{IKS06}, and generalized to a relative situation where we count tropical
curves with ends of higher weights with their real multiplicity. In
\cite{GM05c}, tropical curves with ends of higher weights counted with their
complex multiplicity are shown to determine relative Gromov-Witten invariants
of the plane, i.e.\ numbers of complex plane curves satisfying point conditions
and tangency conditions to a given line $L$. Thus one could imagine that the
tropical relative real count corresponds to numbers of real curves satisfying
point and tangency conditions. This is true only for real curves near the
tropical limit however \cite{Mik03}. The tropical proof of the invariance in
this situation thus led to the construction of new tropical invariant numbers
whose real counterparts are yet to be better understood.

Also, because of the invariance of the tropical relative real count one can
establish a Caporaso-Harris formula for Welschinger invariants for which
$\calP$ consists of only real points. Originally, Caporaso and Harris developed
their algorithm to determine the numbers of complex curves satisfying point
conditions \cite{CH98}. They defined the above mentioned relative Gromov-Witten
invariants and specialized one point after the other to lie on the line $L$.
Since a curve of degree $d$ intersects $L$ in $d$ points, after some steps the
curves become reducible and the line $L$ splits off as a component. One then
collects the contributions from all the components and thus produces recursive
relations among the relative Gromov-Witten invariants that finally suffice to
compute the numbers of complex curves satisfying point conditions. A tropical
counterpart of this algorithm has been established in \cite{GM05c}. There, one
moves one point after the other to the far left part of the plane (but still in
general position). The tropical curves then do not become reducible, but in a
sense decompose into two parts, leading to recursive relations. The left part,
passing through the moved point, is called a floor \cite{BM08}. In \cite{IKS06}
the authors use the same idea to specialize points and consider tropical curves
decomposing into a floor and another part, only now they have to deal with the
real multiplicity for these tropical curves. The formula one thus obtains
computes tropical Welschinger numbers which are equal to their classical
counterparts by the Correspondence Theorem. Since this formula is recursive it
is much more efficient for the computation of Welschinger invariants than the
lattice path algorithm mentioned above. There is also work in progress to
compute Welschinger invariants without tropical methods \cite{Sol07}.

Now let us discuss the situation when $\calP$ does not only contain real
points, but also pairs of complex conjugate points. As already mentioned, also
here a Correspondence Theorem exists to relate these Welschinger invariants to
a certain count of tropical curves, and one can count the tropical curves with
a generalized lattice path algorithm \cite{Shu06}. In addition, it follows of
course again from the Correspondence Theorem together with the Welschinger
Theorem that the tropical count is invariant. However, the tropical count is no
longer locally invariant in the moduli space, and thus there was no known
tropical proof for the (global) invariance of the tropical count. Even worse,
if we try to generalize the tropical count to relative numbers, i.e.\ to curves
with ends of higher weight, then these numbers are no longer invariant.
However, one can still pick a special configuration of points, namely the
result after applying the Caporaso-Harris algorithm as many times as possible.
Then each point is followed by a point which is far more left, and the curves
totally decompose into floors. They can then be counted by means of floor
diagrams. Although the tropical relative count is not invariant, the floor
diagram count leads to a Caporaso-Harris type formula which is sufficient to
compute all Welschinger invariants of the plane \cite{ABM08}. 

\subsection{The content of this paper} \label{sec-content}

The aim of this paper is to give a tropical proof of the invariance of tropical
Welschinger numbers for real and complex conjugate points. As an additional
result this will allow us to construct corresponding tropical invariants in the
relative setting (or more generally for any choice of directions for the ends
of the curve). Using this result, we can then establish a Caporaso-Harris
formula for rational curves in a much simpler way than in \cite{ABM08}. 

The key idea to achieve this is to modify (and in fact also simplify) the class
of tropical curves that we count in order to obtain the invariants. This
modification is small enough so that the (weighted) number of these curves
through given points remains the same in the toric Del Pezzo case, but big
enough so that their count becomes locally invariant in the moduli space.

Let us explain this modification in more detail. For this it is important to
distinguish between odd and even edges of a tropical curve, i.e.\ edges whose
weight is odd resp.\ even. In our pictures we will always draw odd edges as
thin lines and even edges as thick lines. Moreover, we will draw real points as
thin dots and complex points (i.e.\ those corresponding to a pair of complex
conjugate points in the algebraic case) as thick dots. All our curves will be
of genus zero.

The tropical curves that are usually counted to obtain the Welschinger
invariants --- we will call them Welschinger curves --- then have the property
that each connected component of even edges is connected to the rest of the
curve at exactly one point (we can think of such a component as an end tree).
Moreover, real points cannot lie on end trees, and each complex point is either
on an end tree or at a 4-valent vertex \cite {Shu06}. Below on the left we have
drawn a typical (schematic) picture of such a Welschinger curve, with the end
trees marked blue. Note that the marking lying on a point is itself an edge, so
that the 4-valent complex markings away from the end trees look like 3-valent
vertices in the picture.

\begin {center} \input {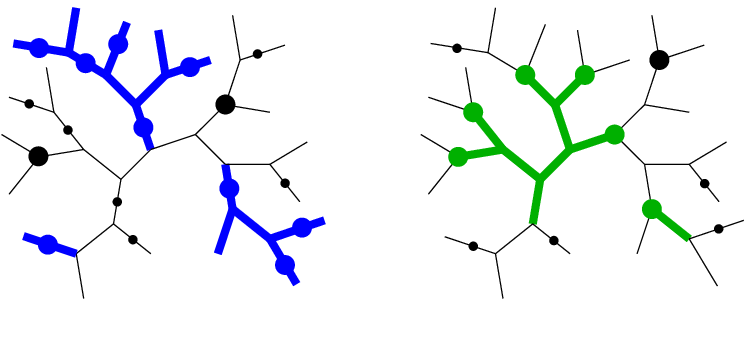} \end {center}

We now change this condition slightly to obtain a different class of curves
that we call broccoli curves: each connected component of even edges can now be
connected to the rest of the curve at several points, of which exactly one is a
3-valent vertex without marking as before (the ``broccoli stem''), and the
remaining ones are complex points (the ``broccoli florets''). The even part of
the curve (the ``broccoli part'') may not contain any marked points in its
interior, whereas away from this part we can have real points at 3-valent and
complex points at 4-valent vertices as before. The picture above on the right
shows a typical schematic example of a broccoli curve, with the broccoli part
drawn in green. Note that, in contrast to Welschinger curves, complex points
are always at 4-valent vertices in broccoli curves.

Broccoli curves have the advantage that their count (with suitably defined
multiplicities) is locally invariant in the moduli space, similarly to the
situation mentioned above when we count complex curves or Welschinger curves
through only real points. Hence counting these curves we obtain well-defined
broccoli invariants --- even for curves with directions of the ends for which
the corresponding Welschinger count would not be invariant of the position of
the points.

In addition, we show that in the toric Del Pezzo case broccoli invariants equal
Welschinger numbers, thereby giving a new and entirely tropical proof of the
invariance of Welschinger numbers. We prove this by constructing bridges
between broccoli curves and Welschinger curves which show that their numbers
must be equal. To illustrate this concept of bridges in an easy example we have
drawn in the picture below a Welschinger curve (which is not a broccoli curve)
and a broccoli curve (which is not a Welschinger curve) of degree $3$ through
the same two real and three complex points. They can be connected by the bridge
drawn below those curves: starting from the Welschinger curve we first split
the vertical end of weight $2$ into two edges of weight $1$ until the rightmost
complex point becomes 4-valent (in the picture at the bottom), and then split
the other end of weight $2$ in a similar way until we arrive at the broccoli
curve.

\begin {center} \input {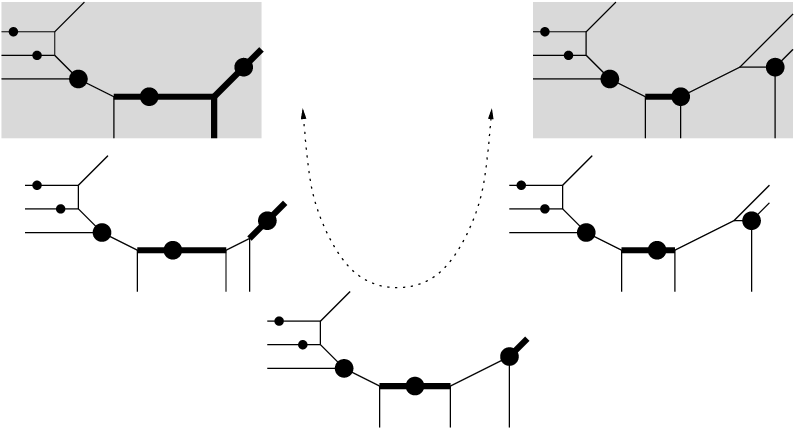} \end {center}

It should be noted that this example is a particularly simple bridge as it
connects a Welschinger curve to a unique corresponding broccoli curve. In
general, traversing bridges will involve creating and resolving higher-valent
vertices of curves along 1-dimensional families --- and as there are usually
several possibilities for such resolutions this means that bridges may ramify
on their way from the Welschinger to the broccoli side. Bridge curves will be
assigned a multiplicity (in a similar way as for Welschinger and broccoli
curves), and at each point of the bridge it is just the weighted number of
incoming Welschinger and outgoing broccoli curves that is the same --- not
necessarily the absolute number of them. In particular, bridges do in general
not provide a bijection between Welschinger and broccoli curves, in fact not
even a well-defined map in either direction.

Another technical thing to note is that we have twice split an even end of
weight 2 into two odd ends of weight $1$ on the bridge above. This might look
like a discontinuous change in the underlying graph of the tropical curve. In
order to avoid this inconvenience we will usually parametrize even ends of
Welschinger curves as two ends of half the weight (which we call double ends).
This way no further end splitting takes place on bridges.

It would certainly be very interesting to see if one could prove a
Correspondence Theorem for broccoli curves that relates these tropical curves
directly to certain real algebraic ones. So far there is no such statement
known; in particular there is no algebraic counterpart to broccoli invariants
for directions of the ends of the curves when the corresponding Welschinger
number is not an invariant.

This paper is organized as follows. In section \ref {sec-oriented} we review
basic notions of tropical curves and their moduli spaces. In particular, we
introduce the notion of oriented curves (i.e.\ tropical curves with the edges
oriented in a certain way), a tool which simplifies proofs in the rest of the
paper. The next three sections are dedicated to the different kinds of tropical
curves mentioned above: section \ref {sec-broccoli} deals with broccoli curves;
the main result here is theorem \ref {thm-broccoliinvariant} which states that
the counts of broccoli curves do not depend on the position of the points. In a
very analogous way, section \ref {sec-welschinger} considers Welschinger curves
and shows that their counts yield the Welschinger invariants. We then introduce
bridge curves in section \ref {sec-bridge} and use them in corollary \ref
{cor-bridge} to prove that Welschinger and broccoli invariants agree in the
toric Del Pezzo case, and thus that the Welschinger invariants then do not
depend on the choice of point conditions (corollary \ref{cor-welinv}). Finally,
the existence of well-defined broccoli invariants also in the relative case
enables us to prove a Caporaso-Harris formula for Welschinger invariants of the
plane in section \ref {sec-ch}.

\subsection {Acknowledgments}

We would like to thank Eugenii Shustin and Inge Sandstad Skrondal for helpful
discussions. Part of this work was accomplished at the Mathematical Sciences
Research Institute (MSRI) in Berkeley, CA, USA, during the one-semester program
on tropical geometry in fall 2009, and part at the Mittag-Leffler Institute in
Stockholm, during the semester program in spring 2011 on ``Algebraic Geometry
with a View towards Applications''. The authors would like to thank both
institutes for hospitality and support. In particular, Andreas Gathmann was
supported by the Simons Professorship of the MSRI.

  \section {Oriented marked curves} \label {sec-oriented}

Let us start by introducing the tropical curves that we will deal with in this
paper. As all our curves will be tropical we usually drop this attribute in the
notation. All curves will be in $ \RR^2 $ (parametrized and labeled in the
sense of \cite {GKM07} section 4), connected, and of genus $0$. Let us quickly
recall the definition of these tropical curves, already making the distinction
between real and complex markings resp.\ odd and even edges that we will later
need to consider real enumerative invariants.

\begin {definition}[Marked curves] \label {def-curve}
  Let $ r,s \in \NN $. An \df {$(r,s)$-marked (plane tropical) curve} is a
  tuple $ C = (\Gamma,x_1,\dots,x_{r+s},y_1,\dots,y_n,h) $ for some $ n \in \NN
  $ such that:
  \begin {enumerate}
  \item \label {def-curve-a}
    $ \Gamma $ is a connected rational metric graph, with unbounded edges
    allowed, and such that each vertex has valence at least 3. The unbounded
    edges of $ \Gamma $ will be called the \df {ends} of $C$.
  \item \label {def-curve-b}
    $ h: \Gamma \to \RR^2 $ is a continuous map that is integer affine
    linear on each edge of $ \Gamma $, i.e.\ on each edge $E$ it is of the form
    $ h(t) = a + t \, v $ for some $ a \in \RR^2 $ and $ v \in \ZZ^2 $. If we
    parametrize $E$ starting at the vertex $ V \in \partial E $ the vector $v$
    in this equation will be denoted $ v(E,V) $ and called the \df {direction
    (vector)} of $E$ starting at $V$. For an end $E$ we will also write $
    v(E) $ instead of $ v(E,V) $, where $V$ is the unique vertex of $E$. We say
    that an edge is \df {contracted} if its direction is $0$.
  \item \label {def-curve-c}
    At each vertex $V$ of $ \Gamma $ the \df {balancing condition}
      \[ \qquad \quad
         \sum_{E: \, V \in \partial E} v(E,V) = 0 \]
    holds.
  \item \label {def-curve-d}
    $ x_1,\dots,x_{r+s} $ is a labeling of the contracted ends, $ y_1,\dots,
    y_n $ a labeling of the non-contracted ends of $C$. We call $
    x_1,\dots,x_{r+s} $ the \df {markings} or \df {marked ends}; more
    specifically the $r$ ends $ x_1,\dots,x_r $ are called the \df {real
    markings}, the $s$ ends $ x_{r+1},\dots,x_{r+s} $ the \df {complex
    markings} of $C$. The other ends $ y_1,\dots,y_n $ are called the \df
    {unmarked ends}; the collection $ (v(y_1),\dots,v(y_n)) $ of their
    directions will be called the \df {degree} $ \Delta = \Delta(C) $ of $C$.
    We denote the number $n$ of vectors in $ \Delta $ by $ |\Delta| $.
  \end {enumerate}
  The set of all $ (r,s) $-marked curves of degree $ \Delta $ will be denoted
  $ \Mrs(\Delta) $.
\end {definition}

\begin {definition}[Even and odd edges, weights] \label {def-weight}
  Let $C$ be a marked curve.
  \begin {enumerate}
  \item \label {def-weight-a}
    A vector in $ \ZZ^2 $ will be called \df {even} if both its coordinates are
    even, and \df {odd} otherwise. We say that an edge of $C$ is even resp.\
    odd if its direction vector is even resp.\ odd.
  \item \label {def-weight-b}
    If we write the direction vector of an edge $E$ of $C$ as a non-negative
    multiple $ \omega(E) $ of a primitive integral vector we call this number $
    \omega(E) $ the \df {weight} of $E$. Note that $E$ is even resp.\ odd if
    and only if its weight is even resp.\ odd.
  \end {enumerate}
\end {definition}

\begin {convention} \label {conv-curve}
  When drawing a marked curve $ C=(\Gamma,x_1,\dots,x_{r+s},y_1,\dots,y_n,h) $
  we will usually only show the image $ h(\Gamma) \subset \RR^2 $, together
  with the image points $ h(x_1),\dots,h(x_{r+s}) $ of the markings. These
  image points will be drawn as small dots for real markings and as big dots
  for complex markings. The other edges will always be displayed as thin lines
  for odd edges and as thick lines for even edges. Unmarked contracted edges
  would not be visible in these pictures, but (although allowed) they will not
  play a special role in this paper.
\end {convention}

\begin {sidepic}{curve} \begin {example} \label {ex-curve}
  Using convention \ref {conv-curve}, the picture on the right shows a
  $(1,1)$-marked plane curve of degree $ ((-2,1),(0,-1),(1,-1),(1,1)) $. It has
  two 3-valent vertices and one 4-valent vertex. The thick edge has direction $
  (-2,0) $ starting at the complex marking. For clarity we have labeled all the
  ends in the picture, but in the future we will usually omit this as the
  actual labeling will not be relevant for most of our arguments.
\end {example} \end {sidepic}

\begin {remark} \label {rem-curve}
  Note that our set $ \Mrs(\Delta) $ is precisely the moduli space $
  \calM_{0,r+s,\text{trop}}^{\text{lab}} (\RR^2,\Delta) $ of $ (r+s) $-marked
  plane labeled tropical curves of \cite {GKM07} definition 4.1. As such it is
  a polyhedral complex, and in fact even a tropical variety (see \cite {GKM07}
  proposition 4.7). In this paper we will not need its structure as a tropical
  variety however, but only consider $ \Mrs(\Delta) $ as an abstract
  polyhedral complex with polyhedral structure induced by the combinatorial
  types of the curves. Let us quickly establish this notation.
\end {remark}

\begin {definition}[Combinatorial types] \label {def-combtype}
  Let $ C = (\Gamma,x_1,\dots,x_{r+s},y_1,\dots,y_n,h) \in \Mrs(\Delta) $ be a
  marked curve. The \df {combinatorial type} of $C$ is the data of the
  non-metric graph $ \Gamma $, together with the labeling $ x_1,\dots,x_{r+s},
  y_1,\dots,y_n $ of the ends and the directions of all edges. For such a
  combinatorial type $ \alpha $ we denote by $ \Mrs^\alpha(\Delta) $ the
  subspace of $ \Mrs(\Delta) $ of all marked curves of type $ \alpha $.
\end {definition}

\begin {remark}[$ \Mrs(\Delta) $ as a polyhedral complex]
    \label {rem-combtype}
  In the same way as in \cite {GM05b} example 2.13 the moduli spaces $ \Mrs
  (\Delta) $ are abstract polyhedral complexes in the sense of \cite {GM05b}
  definition 2.12, i.e.\ they can be obtained by glueing finitely many real
  polyhedra along their faces. The open cells of these complexes are exactly
  the subspaces $ \Mrs^\alpha(\Delta) $, where $ \alpha $ runs over all
  combinatorial types of curves in $ \Mrs(\Delta) $. The curves in such a cell
  (i.e.\ for a fixed combinatorial type) are parametrized by the position in $
  \RR^2 $ of a chosen root vertex and the lengths of all bounded edges (which
  need to be positive). Hence $ \Mrs^\alpha(\Delta) $ can be thought of as an
  open polyhedron whose dimension is equal to 2 plus the number of bounded
  edges in the combinatorial type $ \alpha $. We will call this dimension the
  \df {dimension} $ \dim \alpha $ of the type $ \alpha $.
\end {remark}

Let us now consider enumerative questions for our curves. In addition to the
usual incidence conditions we want to be able to require that some of the
unmarked ends are fixed, i.e.\ map to a given line in $ \RR^2 $. To count such
curves we will now introduce the corresponding evaluation maps. Moreover, to
be able to compensate for the overcounting due to the labeling of the non-fixed
unmarked ends we will define the group of permutations of these ends that keep
the degree fixed.

\begin {definition}[Evaluation maps and $ G(\Delta,F) $] \label {def-ev}
  Let $ r,s \ge 0 $, let $ \Delta = (v_1,\dots,v_n) $ be a collection of
  vectors in $ \ZZ^2 \backslash \{0\} $, and let $ F \subset \{1,\dots,n\} $.
  \begin {enumerate}
  \item \label {def-ev-a}
    The \df {evaluation map} $ \ev_F $ (with \df {set of fixed ends} $F$) on $
    \Mrs(\Delta) $ is defined to be
      \[ \qquad \quad
         \begin {array}{rrcl}
           \ev_F: & \Mrs(\Delta)
                  & \longrightarrow
                  & \displaystyle
                    (\RR^2)^{r+s} \times
                    \prod_{i \in F} \, \big(\RR^2 / \gen {v_i} \big)
                    \qquad \isom \RR^{2(r+s)+|F|} \\[2.5ex]
                  & (\Gamma,x_1,\dots,x_{r+s},y_1,\dots,y_n,h)
                  & \longmapsto
                  & \big( (h(x_1),\dots,h(x_{r+s})), (h(y_i): i \in F) \big).
         \end {array} \]
    In our pictures we will indicate ends that we would like to be considered
    fixed with a small orthogonal bar at the infinite side.
  \item \label {def-ev-b}
    We denote by $ G(\Delta,F) $ the subgroup of the symmetric group $ \SM_n $
    of all permutations such that $ \sigma(i)=i $ for all $ i \in F $ and
    $ v_{\sigma(i)} = v_i $ for all $ i=1,\dots,n $.
  \end {enumerate}
  For the case $ F=\emptyset $ of no fixed ends we denote $ \ev_F $ simply by $
  \ev $ and $ G(\Delta,F) $ by $ G(\Delta) $.
\end {definition}

\begin {remark} \label {rem-ev}
  As in \cite {GM05b} example 3.3 these evaluation maps are morphisms of
  polyhedral complexes in the sense that they are continuous maps that are
  linear on each cell $ \Mrs^\alpha(\Delta) $ of $ \Mrs(\Delta) $. Note that
  $ G(\Delta,F) $ acts on $ \Mrs(\Delta) $ by permuting the unmarked ends, and
  that $ \ev_F $ is invariant under this operation. By definition, if
    \[ \calP = \big( (P_1,\dots,P_{r+s}), (Q_i: i \in F) \big) 
       \;\in\; (\RR^2)^{r+s} \times
       \prod_{i \in F} \big( \RR^2/\gen {v_i} \big) \]
  then the inverse image $ \ev_F^{-1}(\calP) $ consists of all $ (r+s)
  $-marked curves $ (\Gamma,x_1,\dots,x_{r+s},y_1,\dots,y_n,h) $ of degree $
  \Delta $ that pass through $ P_i \in \RR^2 $ at the marked point $ x_i $ for
  all $ i=1,\dots,r+s $ and map the $i$-th unmarked end $ y_i $ to the line $
  Q_i \in \RR^2 / \gen {v_i} $ for all $ i \in F $. We call $ \calP $ a \df
  {collection of conditions} for $ \ev_F $.
\end {remark}

Of course, when counting curves we must assume that the conditions we impose
are in general position so that the dimension of the space of curves satisfying
them is as expected. Let us define this notion rigorously.

\begin {definition}[General and special position of points]
    \label {def-general}
  Let $ N \in \NN $, and let $ f: M \to \RR^N $ be a morphism of polyhedral
  complexes (as e.g.\ the evaluation map $ \ev_F $ of definition \ref
  {def-ev} \ref {def-ev-a}). Then the union $ \bigcup_\alpha f(M^\alpha)
  \subset \RR^N $, taken over all cells $ M^\alpha $ of $M$ such that the
  polyhedron $ f(M^\alpha) $ has dimension at most $ N-1 $, is called the locus
  of points \df {in special position} for $f$. Its complement is denoted the
  locus of points \df {in general position} for $f$.
\end {definition}

\begin {remark} \label {rem-general-dense}
  Note that the locus of points in general position for a morphism $ f: M \to
  \RR^N $ is by definition the complement of finitely many polyhedra of
  positive codimension in $ \RR^N $. In particular, it is a dense open subset
  of $ \RR^N $.
\end {remark}

\begin {example} \label {ex-general}
  Let $ M \subset \Mrs(\Delta) $ be a polyhedral subcomplex, and let $ F
  \subset \{ 1,\dots,|\Delta| \} $. Then a collection of conditions $ \calP \in
  \RR^{2(r+s)+|F|} $ as in remark \ref {rem-ev} is in general position for $
  \ev_F : M \to \RR^{2(r+s)+|F|} $ if and only if for each curve in $M$
  satisfying the conditions $ \calP $ and every small perturbation of these
  conditions we can still find a curve of the same combinatorial type
  satisfying them.
\end {example}

Collections of conditions in general position for the evaluation map have a
special property that will be crucial for the rest of the paper: in \cite
{GM05b} remark 3.7 it was shown that every 3-valent curve $ C=(\Gamma,
x_1,\dots,x_{r+s},y_1,\dots,y_n,h) \in \Mrs(\Delta) $ through a collection of $
r+s=|\Delta|-1 $ points in general position for the evaluation map $ \ev:
\Mrs(\Delta) \to \RR^{2(r+s)} $ without fixed ends has the property that each
connected component of $ \Gamma \backslash (x_1 \cup \cdots \cup x_{r+s}) $
contains exactly one unmarked end. For the purposes of this paper we need the
following generalization of this statement to curves that are not
necessarily 3-valent and evaluation maps that may have fixed ends.

\begin {lemma} \label {lem-general}
  Let $ M \subset \Mrs(\Delta) $ be a polyhedral subcomplex, and let $ \calP $
  be a collection of conditions in general position for the evaluation map $
  \ev_F: M \to \RR^{2(r+s)+|F|} $. Assume that there is a curve $
  C=(\Gamma,x_1,\dots,x_{r+s},y_1,\dots,y_n,h) \in \ev_F^{-1}(\calP) $
  satisfying these conditions. Then:
  \begin {enumerate}
  \item \label {lem-general-a}
    Each connected component of $ \, \Gamma \backslash (x_1 \cup \cdots \cup
    x_{r+s}) $ has at least one unmarked end $ y_i $ with $ i \notin F $.
  \item \label {lem-general-b}
    If the combinatorial type of $C$ has dimension $ 2(r+s)+|F| $ and every
    vertex of $C$ that is not adjacent to a marking is 3-valent then every
    connected component of $ \, \Gamma \backslash (x_1 \cup \cdots \cup
    x_{r+s}) $ as in \ref {lem-general-a} has exactly one unmarked end $ y_i $
    with $ i \notin F $.
  \end {enumerate}
\end {lemma}

\begin {proof}
  Consider a connected component of $ \Gamma \backslash (x_1 \cup \cdots
  \cup x_{r+s}) $ and denote by $ \Gamma' $ its closure in $ \Gamma $. We can
  consider $ \Gamma' $ as a graph, having a certain number $a$ of unbounded
  fixed ends, $b$ unbounded non-fixed ends, and $c$ bounded ends (i.e.\
  1-valent vertices) at markings of $C$. The statement of part \ref
  {lem-general-a} of the lemma is that $ b \ge 1 $, with equality holding in
  case \ref {lem-general-b}. For an example, in the picture below on the right
  $ \Gamma' $ consists of the solidly drawn lines; the curve continues in some
  way behind the dashed lines. Recall that fixed ends are indicated by small
  bars at the infinite sides. Hence in our example we have $ a=1 $, $ b=1 $,
  and $ c=2 $.

  \begin {sidepic}{general}
    By the same argument as in remark \ref {rem-combtype}, the graph $ \Gamma'
    $ as well as the map $ h|_{\Gamma'} $ is fixed by the position of a root
    vertex in $ \Gamma' $ and the lengths of all bounded edges of $ \Gamma' $.
    But an easy combinatorial argument shows that the number of bounded edges
    of $ \Gamma' $ is equal to $ a+b+2c-3-\sum_V (\val V-3) $, with the sum
    taken over all vertices $V$ that are not adjacent to a marking. Hence $
    \Gamma' $ and its image $ h|_{\Gamma'} $ can vary with $ a+b+2c-1-\sum_V
    (\val V-3) $ real parameters in $M$.
  \end {sidepic}

  On the other hand, $ \Gamma' $ together with $ h|_{\Gamma'} $ fixes $ a+2c $
  coordinates in the image of the evaluation map, namely the positions of the
  $a$ fixed ends and the $c$ markings in $ \Gamma' $.

  Hence $ b=0 $ is impossible: then these $ a+2c $ coordinates of the
  evaluation map would vary with fewer than $ a+2c $ coordinates of $M$,
  meaning that the image of $ \ev_F $ on the cell of $C$ cannot be
  full-dimensional and thus $ \calP $ cannot have been in general position.
  This proves \ref {lem-general-a}. But in case \ref {lem-general-b} $ b>1 $ is
  impossible as well: then by assumption we have $ \val V = 3 $ for all $V$ as
  above, and thus one could fix a position for the fixed ends and markings at $
  \Gamma' $ in $ \RR^2 $ and still obtain a $ (b-1) $-dimensional family for $
  \Gamma' $ and $ h|_{\Gamma'} $. As a movement in this family does not change
  anything away from $ \Gamma' $ this means that $ \ev_F $ is not injective on
  the cell of $M$ corresponding to $C$. But $ \ev_F $ is surjective on this
  cell as $ \calP $ is in general position. This is a contradiction since by
  assumption the source and the target of the restriction of $ \ev_F $ to the
  cell corresponding to $C$ have the same dimension.
\end {proof}

\begin {sidepic}{oriented}
  \begin {remark} \label {rem-general}
    The important consequence of lemma \ref {lem-general} \ref {lem-general-b}
    is that --- whenever it is applicable --- it means that there is a unique
    way to orient every unmarked edge of $ C=(\Gamma,x_1,\dots,x_{r+s},
    y_1,\dots,y_n,h) $ so that it points towards the unique unmarked non-fixed
    end of the component of $ \Gamma \backslash (x_1 \cup \cdots \cup x_{r+s})
    $ containing the edge. The picture on the right shows this for the curve of
    example \ref {ex-curve}. Note that the arrow will always point inwards on
    fixed ends, and outwards on non-fixed ends.
  \end {remark}
\end {sidepic}

To be able to talk about this concept in the future we will now introduce the
notion of oriented curves.

\begin {definition}[Oriented marked curves] \label {def-oriented}
  An \df {oriented $ (r,s) $-marked curve} is an $ (r,s) $-marked curve $
  C = (\Gamma,x_1,\dots,x_{r+s},y_1,\dots,y_n,h) $ as in definition \ref
  {def-curve} in which each unmarked edge of $ \Gamma $ is equipped with an
  orientation (which we will draw as arrows in our pictures). In accordance
  with our above idea, the subset $ F = F(C) \subset \{1,\dots,n\} $ of all $i$
  such that the unmarked end $ y_i $ is oriented inwards is called the \df {set
  of fixed ends} of $C$. The space of all oriented $ (r,s) $-marked curves with
  a given degree $ \Delta $ and set of fixed ends $F$ will be denoted $
  \Mrs^\ori(\Delta,F) $; for the case $ F=\emptyset $ of no fixed ends we write
  $ \Mrs^\ori(\Delta,\emptyset) $ also as $ \Mrs^\ori(\Delta) $. We denote by $
  \ft: \Mrs^\ori(\Delta,F) \to \Mrs(\Delta) $ the obvious \df {forgetful map}
  that disregards the information of the orientations.
\end {definition}

\begin {remark} \label {rem-oriented}
  Obviously, our constructions and results for non-oriented curves carry over
  immediately to the oriented case: $ \Mrs^\ori(\Delta,F) $ is a polyhedral
  complex with cells $ \Mrs^\alpha(\Delta,F) $ corresponding to the
  combinatorial types $ \alpha $ of the oriented curves (which now include the
  data of the orientations of all edges). The forgetful map $ \ft $ is a
  morphism of polyhedral complexes that is injective on each cell. There
  are evaluation maps on $ \Mrs^\ori(\Delta,F) $ as in definition \ref {def-ev}
  \ref {def-ev-a} that are morphisms of polyhedral complexes; by abuse of
  notation we will write them as in the unoriented case as $ \ev_F $.
\end {remark}

So far we have allowed any choice of orientations on the edges of our curves in
$ \Mrs^\ori(\Delta,F) $. To ensure that the orientations are actually as
explained in remark \ref {rem-general} we will now allow only certain types of
vertices. In the rest of the paper we will study various kinds of oriented
marked curves --- broccoli curves in section \ref {sec-broccoli}, Welschinger
curves in section \ref {sec-welschinger}, and bridge curves in section \ref
{sec-bridge} --- that differ mainly in their allowed vertex types. The
following definition gives a complete list of all vertex types that will occur
anywhere in this paper.

\begin {definition}[Vertex types and multiplicities] \label {def-vertex}
  We say that a vertex $V$ of an oriented $ (r+s) $-marked curve $C$ is of a
  certain type if the number, parity (even or odd), and orientation of its
  adjacent edges is as in the following table. In addition, two arrows pointing
  in the same direction (as in the types (6b) and (8)) require these odd edges
  to be two unmarked ends with the same direction, and an arc (as in the types
  (6a) and (9)) means that these two odd edges must \textsl {not} be two
  unmarked ends with the same direction. Hence the type (6) splits up into the
  two subtypes (6a) and (6b). All other types in the list are mutually
  exclusive.

  \begin {center} \input {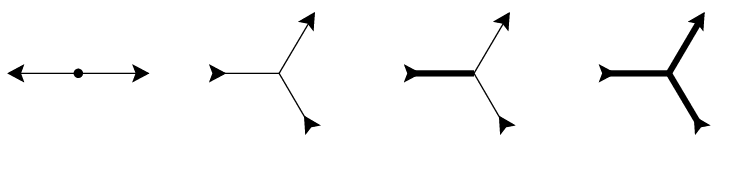} \end {center}

  \begin {center} \input {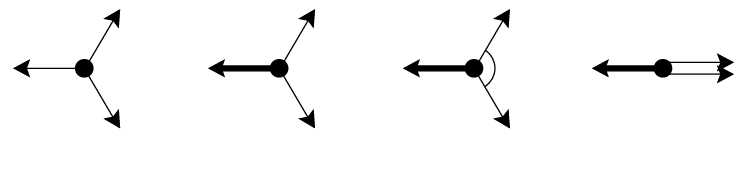} \end {center}

  \begin {center} \input {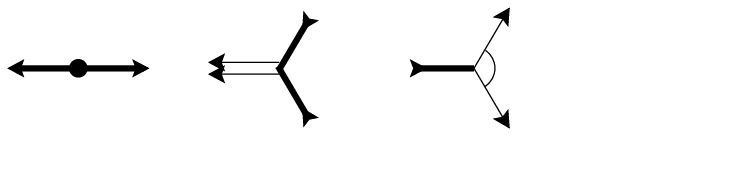} \end {center}

  In addition, each vertex $V$ of one of the above types is assigned a \df
  {multiplicity} $ m_V \in \CC $ that can also be read off from the table.
  Here, the number $a$ denotes the ``complex vertex multiplicity'' in the sense
  of Mikhalkin \cite {Mik03}, i.e.\ the absolute value of the determinant of
  two of the adjacent directions. For the type (8) it is the absolute value of
  the determinant of the two even adjacent directions.

  If $ C=(\Gamma,x_1,\dots,x_{r+s},y_1,\dots,y_n,h) $ consists only of vertices
  of the above types, we denote by $ n_\beta = n_\beta(C) $ the number of
  vertices in $C$ of a given type $ \beta $. In addition, we then define the
  \df {multiplicity} of $C$ to be
    \[ m_C := \prod_{k=1}^n \, i^{\omega(y_k)-1} \cdot
              \prod_V \, m_V, \]
  where the second product is taken over all vertices $V$ of $C$. Although some
  of the vertex multiplicities are complex numbers, the following lemma shows
  that the curve multiplicity $ m_C $ is always real. In fact, the complex
  vertex multiplicities are just a computational trick that makes the ``sign
  factor'', i.e.\ the power of $i$, the same for all the vertex types (2) to
  (6) (which will be the most important ones), leading to easier proofs in the
  rest of the paper.
\end {definition}

\begin {lemma} \label {lem-vertex}
  Every oriented marked curve that has only vertices of the types in definition
  \ref {def-vertex} has a real multiplicity.
\end {lemma}

\begin {proof}
  Let $V$ be a vertex of $C$, and denote by $ E_1,\dots,E_q $ the adjacent
  unmarked edges (so $ q \in \{ 2,3,4 \} $ depending on the type of the
  vertex). Pick's theorem implies that the complex vertex multiplicity $a$ as
  in definition \ref {def-vertex} satisfies $ a = \omega (E_1) + \cdots +
  \omega (E_q) \text{ mod }2$. By checking all vertex types we thus see that in
  each case
    \[ m_V \in \prod_{k=1}^q i^{\omega(E_k)-1} \cdot \RR. \]
  Now every unmarked edge is adjacent to exactly two vertices if it is bounded,
  and adjacent to exactly one vertex if it is unbounded. Hence
    \[ m_C \in \prod_E i^{2(\omega(E)-1)} \cdot \RR = \RR, \]
  where the sum is taken over all unmarked edges.
\end {proof}

\begin {example} \label {ex-vertex}
  The picture of example \ref{ex-curve} and remark \ref{rem-general} shows an
  oriented marked curve $C$ with $ F(C)=\emptyset $. Its vertices $ V_1 $, $
  V_2 $, $ V_3 $, labeled from left to right, are of the types (1), (3), and
  (6), respectively, so that e.g.\ $ \n6 = 1 $. The vertex $ V_3 $ is also of
  type (6a). The multiplicities of the vertices are $ m_{V_1} = 1 $, $ m_{V_2}
  = 2 \cdot i^{2-1} = 2i $, and $ m_{V_3} = i^{2-1} = i $. As all unmarked ends
  of $C$ have weight $1$ the multiplicity of $C$ is thus $ m_C = -2 $.
\end {example}

Let us now check that, with our list of allowed vertex types, in the situation
of lemma \ref {lem-general} \ref {lem-general-b} the only way to orient a given
curve is as explained in remark \ref {rem-general}.

\begin {lemma}[Uniqueness of the orientation of curves] \label {lem-oriented}
  Let the notations and assumptions be as in lemma \ref {lem-general} \ref
  {lem-general-b}. If there is a way to make $C$ into an oriented curve with
  vertices of the types (1) to (7) and so that the orientations of the unmarked
  ends are as given by $F$, this must be the orientation that lets each
  unmarked edge point towards the unique unmarked and non-fixed end in the
  component of $ \Gamma \backslash (x_1 \cup \cdots \cup x_{r+s}) $ containing
  it.
\end {lemma}

\begin {proof}
  By lemma \ref {lem-general} \ref {lem-general-b} there is a unique
  orientation on $C$ pointing on each unmarked edge towards the unmarked and
  non-fixed end in the component of $ \Gamma \backslash (x_1 \cup \cdots \cup
  x_{r+s}) $ containing the edge. Now assume that we have any orientation on
  $C$ with vertices of types (1) to (7). Denote by $ \Gamma' $ the subgraph of
  $ \Gamma $ where these two orientations differ; we have to show that $
  \Gamma' = \emptyset $.

  Note that $ \Gamma' $ is a bounded subgraph since the orientation on the ends
  is fixed by $F$. Moreover, $ \Gamma' $ cannot contain an edge adjacent to a
  marking since all possible vertex types (1), (5), (6), and (7) with markings
  require the orientation on the adjacent edges precisely as in remark \ref
  {rem-general}. So if $ \Gamma' $ is non-empty it must have a 1-valent
  vertex somewhere that is not adjacent to a marking. This can only be a vertex
  of the types (2), (3), or (4), and the condition of $ \Gamma' $ being
  1-valent means that the two orientations differ at exactly one adjacent edge.
  But this is impossible since both orientations have the property that they
  have one adjacent edge pointing outwards and two pointing inwards at this
  vertex.
\end {proof}

We will end this section by computing the dimensions of the cells of $
\Mrs^\ori(\Delta,F) $.

\begin {lemma} \label {lem-dim}
  Let $ C \in \Mrs^\ori(\Delta,F) $ be an oriented marked curve all of whose
  vertices are of the types listed in definition \ref {def-vertex}. Let $
  \alpha $ be the combinatorial type of $C$. Then the cell of $ \Mrs^\ori
  (\Delta,F) $ corresponding to $ \alpha $ has dimension
    \[ \dim \alpha = |\Delta|+r+\n7-\n8-1 = 2(r+s)+|F|+\n9. \]
\end {lemma}

\begin {proof}
  By remark \ref {rem-combtype} it suffices to show that the number of bounded
  edges of $C$ is equal to
    \[ \text {both} \quad |\Delta|+r+\n7(C)-\n8(C)-3
       \qquad
       \text {and} \quad 2(r+s)+|F|+\n9(C)-2. \]
  This is easily proven by induction on the number of vertices in $C$:
  if $C$ has only one vertex (and thus no bounded edge) it has to be one of the
  types in definition \ref {def-vertex}, and the statement is easily checked in
  all of these cases. If the curve $C$ has more than one vertex we cut it at
  any bounded edge into two parts $ C_1 $ and $ C_2 $, making the cut edge
  unbounded in both parts. Note that the cut edge points inward for one part,
  and thus becomes a fixed end for this part. If $ C_i \in M_{(r_i,s_i)}
  (\Delta_i,F_i) $ for $ i=1,2 $, then $ r=r_1+r_2 $, $ s=s_1+s_2 $, $
  |\Delta|=|\Delta_1|+|\Delta_2| -2 $, $ |F|=|F_1|+|F_2|-1 $, and $ n_\beta(C)
  = n_\beta(C_1)+n_\beta(C_2) $ for $ \beta \in \{ \text{(7)}, \text{(8)},
  \text{(9)} \} $. The number of bounded edges of $C$ is now just the number of
  bounded edges in $ C_1 $ and $ C_2 $ plus $1$, i.e.\ by induction equal to
  \begin {align*}
    &   |\Delta_1|+r_1+\n7(C_1)-\n8(C_1)-3
      + |\Delta_2|+r_2+\n7(C_2)-\n8(C_2)-3
      +1 \\
    & \qquad = |\Delta|+r+\n7(C)-\n8(C)-3
  \end {align*}
  as well as
  \begin {align*}
    &   2(r_1+s_1)+|F_1|+\n9(C_1)-2
      + 2(r_2+s_2)+|F_2|+\n9(C_2)-2
      +1 \\
    & \qquad = 2(r+s)+|F|+\n9(C)-2. \qedhere
  \end {align*}
\end {proof}

  \section {Broccoli curves} \label {sec-broccoli}

In this section we will introduce the most important type of curves considered
in this paper: the broccoli curves. We define corresponding numbers, and
show that they are independent of the chosen point conditions.

Broccoli curves can be defined with or without orientation. Both definitions
have their advantages: the oriented one is easier to state and local at the
vertices, whereas the unoriented one is easier to visualize (as one does not
need to worry about orientations at all). So let us give both definitions and
show that they agree for enumerative purposes.

\begin {definition}[Broccoli curves] \label {def-broccoli}
  Let $ r,s \ge 0 $, let $ \Delta = (v_1,\dots,v_n) $ be a collection of
  vectors in $ \ZZ^2 \backslash \{0\} $, and let $ F \subset \{1,\dots,n\} $.
  \begin {enumerate}
  \item \label {def-broccoli-a}
    An oriented curve $ C \in \Mrs^\ori(\Delta,F) $ all of whose vertices are
    of the types (1) to (6) of definition \ref {def-vertex} is called an \df
    {oriented broccoli curve}.
  \item \label {def-broccoli-b}
    Let $ C = (\Gamma,x_1,\dots,x_{r+s},y_1,\dots,y_n,h) \in \Mrs(\Delta) $.
    Consider the subgraph $ \Gamma_\even $ of $ \Gamma $ of all even edges
    (including the markings). The 1-valent vertices of $ \Gamma_\even $ as well
    as the $ y_i \subset \Gamma_\even $ with $ i \notin F $ are called the \df
    {stems} of $ \Gamma_\even $. We say that $C$ is an \df {unoriented broccoli
    curve} (with set of fixed ends $F$) if
    \begin {enumerate}
    \item [(i)]
      all complex markings are adjacent to 4-valent vertices;
    \item [(ii)]
      every connected component of $ \Gamma_\even $ has exactly one stem.
    \end {enumerate}
  \end {enumerate}
\end {definition}

\begin {example} \label {ex-broccoli}
  The picture below shows an oriented broccoli curve in which every allowed
  vertex type appears. We have labeled the vertices with their types. Note that
  by forgetting the orientations of the edges (and thus also disregarding the
  vertex types) one obtains an unoriented broccoli curve. Its subgraph $
  \Gamma_\even $ of even edges consists of all markings and thick edges. It has
  four connected components $ \Gamma_1,\dots,\Gamma_4 $, and each component has
  exactly one stem: the non-fixed unmarked end in $ \Gamma_1 $, the vertex of
  type (3) in $ \Gamma_2 $, and the unique vertices in $ \Gamma_3 $ and $
  \Gamma_4 $.

  \begin {center} \input {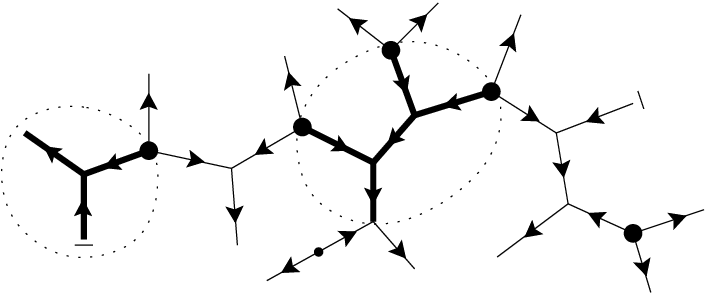} \end {center}
\end {example}

Of course, to count these curves we have to fix the right number of conditions
to get a finite answer. This dimension condition follows e.g.\ for oriented
broccoli curves from lemma \ref {lem-dim}: we must have $ r+2s+|F| = |\Delta|-1
$ since $ \n7=\n8=\n9=0 $.

\begin {proposition}[Equivalence of oriented and unoriented broccoli curves]
    \label {prop-broccoli}
  Let $ r,s \ge 0 $, let $ \Delta = (v_1,\dots,v_n) $ be a collection of
  vectors in $ \ZZ^2 \backslash \{0\} $, and let $ F \subset \{1,\dots,n\} $
  such that $ r+2s+|F| = |\Delta|-1 $. Moreover, let $ \calP \in
  \RR^{2(r+s)+|F|} $ be a collection of conditions in general position for $
  \ev_F: \Mrs(\Delta) \to \RR^{2(r+s)+|F|} $ (see example \ref {ex-general}).

  Then the forgetful map $ \ft $ of definition \ref {def-oriented} gives a
  bijection between oriented and unoriented $ (r,s) $-marked broccoli curves
  through $ \calP $ with degree $ \Delta $ and set of fixed ends $F$.
\end {proposition}

\begin {proof}
  We have to prove three statements.
  \begin {enumerate}
  \item \label {prop-broccoli-a}
    $ \ft $ maps oriented to unoriented broccoli curves through $ \calP $:
    Let $ C \in \Mrs^\ori(\Delta,F) $ be an oriented broccoli curve. The list
    of allowed vertex types for $C$ implies immediately that $C$ then satisfies
    condition (i) of definition \ref {def-broccoli}.

    To show (ii) let $ \Gamma' $ be a connected component of $ \Gamma_\even $.
    If $ \Gamma' $ contains no vertex of type (4) it can only be a single
    marking (types (1) or (5)) or a single unmarked edge with possibly attached
    markings (vertex types (3) together with (6), (3) with a fixed unmarked
    end, or (6) with a non-fixed unmarked end), and in each of these cases
    condition (ii) is satisfied. If there are vertices of type (4) they must
    form a tree in $ \Gamma' $, and obviously every such tree made up from type
    (4) vertices has exactly one outgoing end. This unique outgoing end must be
    a non-fixed end of $C$ or connected to a type (3) vertex, hence in any case
    it leads to a stem. On the other hand, the incoming ends of the tree must
    be fixed ends of $C$ or connected to a type (6) vertex, i.e.\ they never
    lead to a stem. Consequently, $ \Gamma' $ satisfies condition (ii).
  \item \label {prop-broccoli-b}
    $ \ft $ is injective on the set of curves through $ \calP $: Note that the
    conditions of lemma \ref {lem-general} \ref {lem-general-b} are satisfied
    by the dimension condition of lemma \ref {lem-dim} and our list of allowed
    vertex types. Hence lemma \ref {lem-oriented} implies that there is at most
    one possible orientation on $C$.
  \item \label {prop-broccoli-c}
    $ \ft $ is surjective on the set of curves through $ \calP $: Let $ C \in
    \Mrs(\Delta) $ be an unoriented broccoli curve through $ \calP $ with set
    of fixed ends $F$. Then by (i) the curve $C$ has $s$ 4-valent vertices at
    the complex markings, so by \cite {GM05b} proposition 2.11 the
    combinatorial type of $C$ has dimension $ |\Delta|-1+r-\sum_V (\val V-3) =
    2(r+s)+|F|-\sum_V (\val V-3) $, with the sum taken over all vertices $V$
    that are not adjacent to a complex marking. But as $ \calP $ is in general
    position this dimension cannot be less than $ 2(r+s)+|F| $. So we see
    that all vertices without adjacent complex marking are 3-valent, and that
    the combinatorial type of $C$ has dimension equal to $ 2(r+s)+|F| $. Hence
    we can apply lemma \ref {lem-general} \ref {lem-general-b} again to
    conclude that there is an orientation on $C$ that points on each edge
    towards the unique non-fixed unmarked end in $ \Gamma \backslash (x_1 \cup
    \cdots \cup x_{r+s}) $.

    It remains to be shown that with this orientation the only vertex types
    occurring in $C$ are (1) to (6). For this, note that for a vertex $V$
    \begin {itemize}
    \item as we have said above, $V$ is 4-valent if there is a complex marking
      at $V$, and 3-valent otherwise;
    \item by the construction of the orientation, all edges at $V$ are
      oriented outwards if there is a marking at $V$, and exactly one edge is
      oriented outwards otherwise;
    \item by the balancing condition, it is impossible that exactly one
      edge at $V$ is odd.
    \end {itemize}
    With these restrictions, the only possible vertex types besides (1) to (6)
    would be the ones in the picture below.
    \begin {center} \begin{picture}(0,0)%
\includegraphics{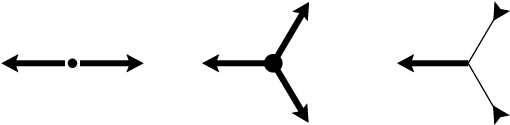}%
\end{picture}%
\setlength{\unitlength}{4144sp}%
\begingroup\makeatletter\ifx\SetFigFont\undefined%
\gdef\SetFigFont#1#2#3#4#5{%
  \reset@font\fontsize{#1}{#2pt}%
  \fontfamily{#3}\fontseries{#4}\fontshape{#5}%
  \selectfont}%
\fi\endgroup%
\begin{picture}(3886,948)(619,-621)
\end{picture}%
 \end {center}
    To exclude these three cases, note that in all of them $V$ would be
    contained in a connected component $ \Gamma' $ of $ \Gamma_\even $ that
    contains at least one unmarked edge. So let us consider such a component,
    and let $ W \in \Gamma' \cap (\Gamma \backslash \Gamma') $ be a vertex
    where $ \Gamma' $ meets the complement of $ \Gamma' $. Then there must be
    an odd as well as an unmarked even edge in $ \Gamma $ at $W$, so by the
    balancing condition as above there are exactly two odd edges and one even
    unmarked edge at $W$. Hence $W$ is a stem if and only if there
    is no marking at $W$. So a connection in $ \Gamma \backslash (x_1 \cup
    \cdots \cup x_{r+s}) $ from a point in the interior of $ \Gamma' $ to a
    non-fixed unmarked end can only be via a stem --- which is unique by (ii).
    This means that every point in the interior of $ \Gamma' $ must be
    connected in $ \Gamma \backslash (x_1 \cup \cdots \cup x_{r+s}) $ to the
    stem. In particular, the interior of $ \Gamma' $ can have no further
    markings, which rules out the first two vertex types in the picture above.
    The third vertex type is impossible since this would have to be the stem
    and thus the connection from $ \Gamma' $ to the non-fixed unmarked end,
    which does not match with the orientation of the even edge. \qedhere
  \end {enumerate}
\end {proof}

Let us now make the obvious definition of the enumerative invariants
corresponding to broccoli curves. Proposition \ref {prop-broccoli} tells us
that it does not matter whether we count oriented or unoriented broccoli
curves. We choose the oriented ones here as their definition is easier. So we
make the convention that from now on \textsl {a broccoli curve will always mean
an oriented broccoli curve}.

\begin {notation} \label {not-broccoli}
  We denote by $ \Mrs^B(\Delta,F) $ the closure of the space of all broccoli
  curves in $ \Mrs^\ori(\Delta,F) $; this is obviously a polyhedral subcomplex.
  By lemma \ref {lem-dim} it is non-empty only if the dimension condition $
  r+2s+|F| = |\Delta|-1 $ is satisfied. Moreover, in this case it is of pure
  dimension $ 2(r+s)+|F| $, and its maximal open cells correspond exactly to
  the broccoli curves in $ \Mrs^B(\Delta,F) $.
\end {notation}

\begin {definition}[Broccoli invariants] \label {def-broccoli-inv}
   As above, let $ r,s \ge 0 $, let $ \Delta = (v_1,\dots,v_n) $ be a
   collection of vectors in $ \ZZ^2 \backslash \{0\} $, and let $ F \subset
   \{1,\dots,n\} $ such that $ r+2s+|F| = |\Delta|-1 $. Moreover, let $ \calP
   \in \RR^{2(r+s)+|F|} $ be a collection of conditions in general position for
   broccoli curves, i.e.\ for the evaluation map $ \ev_F: \Mrs^B(\Delta,F) \to
   \RR^{2(r+s)+|F|} $. Then we define the \df {broccoli invariant}
     \[ \Nrs^B (\Delta,F,\calP) \;:=\; \frac 1{|G(\Delta,F)|} \cdot
          \sum_C m_C, \]
   where the sum is taken over all broccoli curves $C$ in $\Mrs^B(\Delta,F)$
   with degree $\Delta$, set of fixed ends $F$, and $ \ev(C) = \calP $. The
   group $ G(\Delta,F) $ as in definition \ref {def-ev} \ref {def-ev-b} takes
   care of the overcounting of curves due to relabeling the non-fixed unmarked
   ends. The sum is finite by the dimension statement of notation \ref
   {not-broccoli}, and the multiplicity $ m_C $ is as in definition \ref
   {def-vertex}.
\end {definition}

The main result of this section --- and in fact the most important point that
distinguishes our new invariants from the otherwise quite similar Welschinger
invariants that we will study in section \ref {sec-welschinger} --- is that
broccoli invariants are always independent of the choice of conditions $ \calP
$.

\begin {theorem} \label {thm-broccoliinvariant}
  The broccoli invariants $ \Nrs^B (\Delta,F,\calP) $ are independent of the
  collection of conditions $ \calP $. We will thus usually write them simply as
  $ \Nrs^B (\Delta,F) $ (or $ \Nrs^B(\Delta) $ for $ F=\emptyset $).
\end {theorem}

\begin {proof}
  The proof follows from a local study of the moduli space $ \Mrs^B(\Delta,F)
  $. Compared to the one for ordinary tropical curves in \cite {GM05a} theorem
  4.8 it is very similar in style and conceptually not more complicated; there
  are just (many) more cases to consider because we have to distinguish
  orientations as well as even and odd edges.

  By definition, the multiplicity of a curve depends only on its combinatorial
  type. So it is obvious that the function $ \calP \mapsto
  \Nrs^B(\Delta,F,\calP) $ is \textsl {locally} constant on the open subset of
  $ \RR^{2(r+s)+|F|} $ of conditions in general position for broccoli curves,
  and may jump only at the image under $ \ev_F $ of the boundary of
  top-dimensional cells of $ \Mrs^B(\Delta,F) $. This image is a union of
  polyhedra in $ \RR^{2(r+s)+|F|} $ of positive codimension. It suffices to
  show that the function $ \calP \mapsto \Nrs^B(\Delta,F,\calP) $ is locally
  constant around a cell in this image of codimension 1 in $ \RR^{2(r+s)+|F|} $
  since any two top-dimensional cells of $ \RR^{2(r+s)+|F|} $ can be connected
  to each other through codimension-1 cells.

  \begin {sidepic}{invariant}
    So let $ \alpha $ be a combinatorial type in $ \Mrs^B(\Delta,F) $ of
    dimension $ 2(r+s)+|F|-1 $ such that $ \ev_F $ is injective on $
    \Mrs^\alpha (\Delta,F) $ and thus maps this cell to a unique hyperplane $H$
    in $ \RR^{2(r+s)+|F|} $. As in the picture on the right let $ U_\alpha
    \subset \Mrs^B(\Delta,F) $ be the open subset consisting of $ \Mrs^\alpha
    (\Delta,F) $ together with all adjacent top-dimensional cells of $ \Mrs^B
    (\Delta,F) $. To prove the theorem we will show that for a point $ \calP $
    in a neighborhood of $ \ev_F (\Mrs^\alpha (\Delta,F)) $ the sum of the
    multiplicities of the curves in $ U_\alpha \cap \ev_F^{-1}(\calP) $ does
    not depend on $ \calP $, i.e.\ is the same on both sides of $H$. In our
    picture this would just mean that $ m_{\text {I}} + m_{\text {II}} =
    m_{\text {III}} $, where $ m_{\text {I}}, m_{\text {II}}, m_{\text {III}} $
    denote the multiplicities of $ C_{\text {I}}, C_{\text {II}}, C_{\text
    {III}} $, respectively.

  \end {sidepic}

  Actually, we will show this in a slightly different form: to each
  codimension-0 type $ \alpha_k $ in $ U_\alpha $ we will associate a so-called
  \df {$H$-sign} $ \sigma_k $ that is $1$ or $-1$ depending on the side of $H$
  on which $ \ev_F(\Mrs^{\alpha_k} (\Delta,F)) $ lies (it will be $0$ if $
  \ev_F(\Mrs^{\alpha_k} (\Delta,F)) \subset H $). So in the picture above on
  the right we could take $ \sigma_{\text {I}} = \sigma_{\text {II}} = 1 $ and
  $ \sigma_{\text {III}} = -1 $. We then obviously have to show that $ \sum_k
  \sigma_k \, m_k = 0 $, where the sum is taken over all top-dimensional cells
  adjacent to $ \alpha $.

  To prove this, we will start by listing all codimension-1 combinatorial
  types $ \alpha $ in $ \Mrs^B(\Delta,F) $. They are obtained by shrinking the
  length of a bounded edge in a broccoli curve to zero, thereby merging two
  vertices into one. Depending on the merging vertex types we distinguish the
  following cases:
  \begin {itemize}
  \item [(A)] a vertex (1) merging with a vertex (2)/(3), leading to a 4-valent
    vertex with one real marking, two outgoing edges, and one incoming edge.
  \item [(B)] a vertex (2)/(3)/(4) merging with a vertex (2)/(3)/(4), leading
    to a 4-valent vertex with no marking, one outgoing edge, and three incoming
    edges.
  \item [(C)] a vertex (5)/(6) merging with a vertex (2)/(3)/(4), leading to a
    5-valent vertex with one complex marking, three outgoing edges, and one
    incoming edge.
  \end {itemize}
  More precisely, noting that by the balancing condition it is impossible to
  have exactly one odd edge at a vertex, the cases (A), (B), and (C) split up
  into the following possibilities depending on the orientation and parity of
  the adjacent edges.

  \begin {center} \input {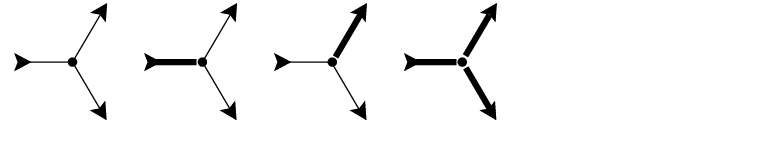} \end {center}
  \vspace {-3mm}
  \begin {center} \input {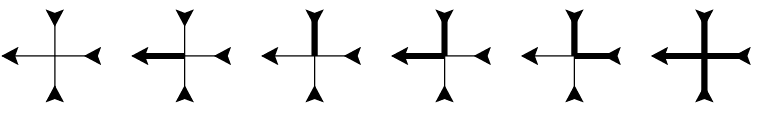} \end {center}
  \vspace {-3mm}
  \begin {center} \input {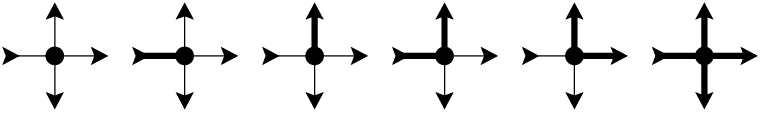} \end {center}

  Next, we will list the adjacent codimension-0 types in $ \Mrs^B(\Delta,F) $
  (called \df {resolutions}) that make up $ U_\alpha $ in the cases (A), (B),
  and (C). In this picture, the dashed lines can be even or odd depending on
  which of the subcases (A$\,\cdot\,$), (B$\,\cdot\,$), (C$\,\cdot\,$) we are
  in. The vectors $ v_1,\dots,v_4 $ will be used in the computations below;
  they are always meant to be oriented outwards (i.e.\ \textsl {not}
  necessarily in the direction of the orientation of the edge), so that $
  v_1+v_2+v_3=0 $ in case (A) and $ v_1+v_2+v_3+v_4=0 $ in the cases (B) and
  (C).

  \begin {center} \input {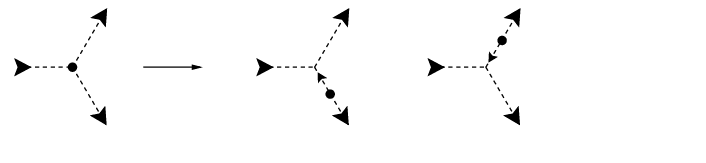} \end {center}
  \vspace {-2mm}
  \begin {center} \input {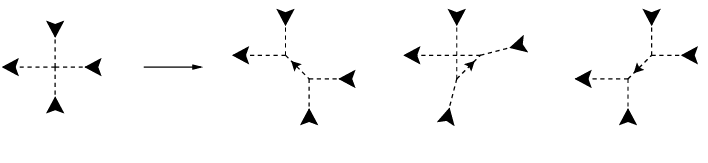} \end {center}
  \vspace {-2mm}
  \begin {center} \input {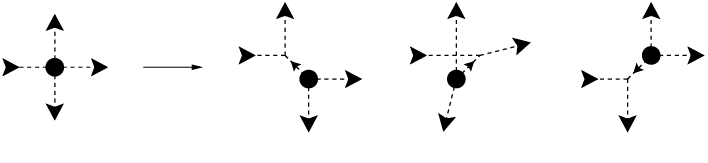} \end {center}

  Note that the allowed vertex types for broccoli curves fix the orientation of
  the newly inserted bounded edge in all these resolutions; it is already
  indicated in the picture above. Moreover, the requirement that there cannot
  be exactly one odd edge at a vertex fixes the parity of the new bounded edge
  in all cases except (B1) and (C1). In the (B1) and (C1) cases, there are two
  possibilities: the four vectors $ v_1,\dots,v_4 $ can either be all the same
  in $ (\ZZ_2)^2 $ (in which case the new bounded edge joining $V$ and $W$ is
  even in all three types I, II, III; we call this case (B$1_3$) and (C$1_3$),
  respectively), or they make up two non-zero equivalence classes in $
  (\ZZ_2)^2 $ (in which case the new bounded edge is even in exactly one of the
  types I, II, III; we call this case (B$1_1$) and (C$1_1$), respectively). In
  the (B$1_1$) and (C$1_1$) cases, we can assume by symmetry that the even
  bounded edge occurs in type I. So in total we now have 18 codimension-1 cases
  (A1), \dots, (A4), (B$1_1$), (B$1_3$), (B2),\dots (B6), (C$1_1$), (C$1_3$),
  (C2),\dots (C6) to consider, and in each of these cases we know the
  resolutions together with all parities and orientations of all edges of the
  curves --- in particular, with the vertex types of $V$ and $W$ (as in the
  picture above). For example, in case (B6) the new bounded edge must be even
  in all three resolutions. Hence in all three resolutions all edges are even,
  and thus both vertices $V$ and $W$ are of type (4).
  
  The following table lists the vertex types for $V$ and $W$ for all
  resolutions I, II, III of all codimension-1 cases. The symbol ``---'' means
  that the required vertex type is not allowed in broccoli curves and thus that
  a corresponding codimension-0 cell does not exist. The columns labeled $ m_*
  $ and $ \mu_*/\mu_* $ will be explained below.

  \begin {center} \tabcolsep 1mm \begin {tabular}{|c|ccc|cccc|} \hline
    codim-1 & \multicolumn 3 {c|}{resolution I}
            & \multicolumn 4 {c|}{resolution II} \\
    case    & $V$ & $W$ & $m_{\text {I}}$
            & $V$ & $W$ & $\mu_{\text {II}}/\mu_{\text {I}}$
            & $m_{\text {II}}$ \\ \hline
    A1      & (2) & (1) & 1
            & (2) & (1) & $-1$ & 1 \\
    A2      & (3) & (1) & $(v_1,v_2)$
            & (3) & (1) & 1 & $(v_1,v_3)$ \\
    A3      & --- & (1) & 0
            & (3) & --- & 1 & 0 \\
    A4      & (4) & --- & 0
            & (4) & --- & 1 & 0 \\ \hline
  \end {tabular} \end {center}
  \begin {center} \tabcolsep 1mm \begin {tabular}{|c|ccc|ccc@{}c|ccc@{}c|}
    \hline
    codim-1 & \multicolumn 3 {c|}{resolution I}
            & \multicolumn 4 {c|}{resolution II}
            & \multicolumn 4 {c|}{resolution III} \\
    case    & $V$ & $W$ & $m_{\text {I}}$
            & $V$ & $W$ & $\mu_{\text {II}}/\mu_{\text {I}}$
	    & $m_{\text {II}}$
            & $V$ & $W$ & $\mu_{\text {III}}/\mu_{\text {I}}$
	    & $m_{\text {III}}$ \\ \hline
    B$1_1$  & (3) & --- & 0
            & (2) & (2) & 1 & 1
            & (2) & (2) & $-1$ & 1 \\
    B$1_3$  & (3) & --- & 0
            & (3) & --- & 1 & 0
            & (3) & --- & 1 & 0 \\
    B2      & --- & (2) & 0
            & --- & (2) & 1 & 0
            & --- & (2) & 1 & 0 \\
    B3      & (3) & (2) & $(v_1,v_2)$
            & (2) & (3) & 1 & $(v_4,v_2)$
            & (2) & (3) & $-1$ & $(v_2,v_3)$ \\
    B4      & (4) & --- & 0
            & --- & (3) & 1 & 0
            & --- & (3) & 1 & 0 \\
    B5      & (3) & (3) & $(v_1,v_2)(v_3,v_4)$
            & (3) & (3) & 1 & $(v_1,v_3)(v_4,v_2)$
            & (3) & (4) & 1 & $(v_1,v_4)(v_2,v_3)$ \\
    B6      & (4) & (4) & $(v_1,v_2)(v_3,v_4)$
            & (4) & (4) & 1 & $(v_1,v_3)(v_4,v_2)$
            & (4) & (4) & 1 & $(v_1,v_4)(v_2,v_3)$ \\ \hline
  \end {tabular} \end {center}
  \begin {center} \tabcolsep 1mm \begin {tabular}{|c|ccc|ccc@{}c|ccc@{}c|}
    \hline
    codim-1 & \multicolumn 3 {c|}{resolution I}
            & \multicolumn 4 {c|}{resolution II}
            & \multicolumn 4 {c|}{resolution III} \\
    case    & $V$ & $W$ & $m_{\text {I}}$
            & $V$ & $W$ & $\mu_{\text {II}}/\mu_{\text {I}}$
            & $m_{\text {II}}$
            & $V$ & $W$ & $\mu_{\text {III}}/\mu_{\text {I}}$
            & $m_{\text {III}}$ \\ \hline
    C$1_1$  & (3) & (6) & $(v_1,v_2)$
            & (2) & (5) & 1 & $(v_4,v_2)$
            & (2) & (5) & $-1$ & $(v_2,v_3)$ \\
    C$1_3$  & (3) & (6) & $(v_1,v_2)$
            & (3) & (6) & 1 & $(v_1,v_3)$
            & (3) & (6) & 1 & $(v_1,v_4)$ \\
    C2      & (3) & (5) & $(v_1,v_2)(v_3,v_4)$
            & (3) & (5) & 1 & $(v_1,v_3)(v_4,v_2)$
            & (3) & (5) & 1 & $(v_1,v_4)(v_2,v_3)$ \\
    C3      & --- & (5) & 0
            & (2) & (6) & 1 & 1
            & (2) & (6) & $-1$ & 1 \\
    C4      & (4) & (6) & $(v_1,v_2)$
            & (3) & (6) & 1 & $(v_1,v_3)$
            & (3) & (6) & 1 & $(v_1,v_4)$ \\
    C5      & --- & (6) & 0
            & --- & (6) & 1 & 0
            & (3) & --- & 1 & 0 \\
    C6      & (4) & --- & 0
            & (4) & --- & 1 & 0
            & (4) & --- & 1 & 0 \\ \hline
  \end {tabular} \end {center}
  Let us now determine the $H$-sign of the resolutions above, i.e.\ figure out
  which of them occur on which side of $H$. To do this we set up the system of
  linear equations determining the lengths of the bounded edges of the curve in
  terms of the positions of the markings in $ \RR^2 $. For such a given
  position of the markings (on the one or on the other side of $H$), a given
  resolution type is then possible if and only if the required length of the
  new bounded edge is positive.

  More concretely, let $a$ be the length of the newly created bounded edge, and
  denote by $ P \in \RR^2 $ in the cases (A) and (C) the required image point
  for the marking. In the cases (A) and (C) the end $ v_1 $ is fixed, so to
  determine the existing resolutions we may assume that there is another
  marking on the $ v_1 $ end at a distance of $ l_1 $ on the graph that is
  required to map to a point $ P_1 \in \RR^2 $. In the case (B) the ends $ v_2
  $, $ v_3 $, and $ v_4 $ are fixed, so we do the same then with lengths $
  l_2,l_3,l_4 $ and points $ P_2,P_3,P_4 \in \RR^2 $. As an example, these
  notions are illustrated for the resolution I in the following picture.
  \begin {center} \input {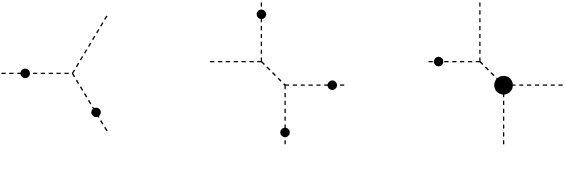} \end {center}
  The systems of linear equations that determine the relative positions of $
  P,P_1,\dots,P_4 $ in terms of $ a,l_1,\dots,l_4 $ are then as follows (where
  all entries are in $ \RR^2 $ and thus each row stands for two equations).
  \begin {center} \tabcolsep 1mm
    \begin {tabular}{|cc|c|}
      \multicolumn 3c{(A)-I} \\ \hline
      $l_1$ & $a$ & \\ \hline
      $-v_1$ & $v_3$ & $P-P_1$ \\ \hline
    \end {tabular}
    \hspace {3mm}
    \begin {tabular}{|cc|c|}
      \multicolumn 3c{(A)-II} \\ \hline
      $l_1$ & $a$ & \\ \hline
      $-v_1$ & $v_2$ & $P-P_1$ \\ \hline
    \end {tabular}
  \end {center}
  \begin {center} \tabcolsep 1mm
    \begin {tabular}{|cccc|c|}
      \multicolumn 5c{(B)-I} \\ \hline
      $l_2$ & $l_3$ & $l_4$ & $a$ & \\ \hline
      $-v_2$ & $v_3$ & $0$ & $v_3+v_4$ & $P_3-P_2$ \\
      $-v_2$ & $0$ & $v_4$ & $v_3+v_4$ & $P_4-P_2$ \\ \hline
    \end {tabular}
    \hspace {3mm}
    \begin {tabular}{|cccc|c|}
      \multicolumn 5c{(B)-II} \\ \hline
      $l_2$ & $l_3$ & $l_4$ & $a$ & \\ \hline
      $-v_2$ & $v_3$ & $0$ & $v_1+v_3$ & $P_3-P_2$ \\
      $-v_2$ & $0$ & $v_4$ & $0$ & $P_4-P_2$ \\ \hline
    \end {tabular}
    \hspace {3mm}
    \begin {tabular}{|cccc|c|}
      \multicolumn 5c{(B)-III} \\ \hline
      $l_2$ & $l_3$ & $l_4$ & $a$ & \\ \hline
      $-v_2$ & $v_3$ & $0$ & $0$ & $P_3-P_2$ \\
      $-v_2$ & $0$ & $v_4$ & $v_1+v_4$ & $P_4-P_2$ \\ \hline
    \end {tabular}
  \end {center}
  \begin {center} \tabcolsep 1mm
    \begin {tabular}{|cc|c|}
      \multicolumn 3c{(C)-I} \\ \hline
      $l_1$ & $a$ & \\ \hline
      $-v_1$ & $v_3+v_4$ & $P-P_1$ \\ \hline
    \end {tabular}
    \hspace {3mm}
    \begin {tabular}{|cc|c|}
      \multicolumn 3c{(C)-II} \\ \hline
      $l_1$ & $a$ & \\ \hline
      $-v_1$ & $v_2+v_4$ & $P-P_1$ \\ \hline
    \end {tabular}
    \hspace {3mm}
    \begin {tabular}{|cc|c|}
      \multicolumn 3c{(C)-III} \\ \hline
      $l_1$ & $a$ & \\ \hline
      $-v_1$ & $v_2+v_3$ & $P-P_1$ \\ \hline
    \end {tabular}
  \end {center}
  To determine $a$ in terms of $ P,P_1,\dots,P_4 $ we use Cramer's rule: if $M$
  is the (quadratic) matrix of a system of linear equations as above and $M'$
  the matrix obtained from $M$ by replacing the $a$-column by the right hand
  side of the equation, then $ a = \det M' / \det M $. But within a case (A),
  (B), (C) the matrix $M'$ does not depend on the resolution I, II, III, and
  thus it is simply the sign of $ \det M $ that tells us whether $a$ is
  positive or negative, i.e.\ whether this resolution exists for the chosen
  points $ P,P_1,\dots,P_4 $. We can therefore take the $H$-sign to be the sign
  of $ \det M $ (note that this will be $0$ if and only if the relative
  position of $ P,P_1,\dots,P_4 $ is not determined uniquely by the equations
  and thus if and only if the codimension-0 cell maps to $H$). An elementary
  computation of the determinants shows that these $H$-signs are as in the
  following table, where $ (v_i,v_j) $ stands for the determinant of the
  $ 2 \times 2 $ matrix with columns $ v_i,v_j $ (and where we have used $
  v_1+v_2+v_3=0 $ in case (A) as well as $ v_1+v_2+v_3+v_4=0 $ in the cases (B)
  and (C)).
  \begin {center} \begin {tabular}{|c|ccc|} \hline
        & $H$-sign for I & $H$-sign for II & $H$-sign for III \\ \hline
    (A) & $ \sign (v_1,v_2) $ & $ \sign (v_1,v_3) $ & \\
    (B) & $ \sign \big((v_1,v_2) (v_3,v_4)\big) $ &
          $ \sign \big((v_1,v_3) (v_4,v_2)\big) $ &
          $ \sign \big((v_1,v_4) (v_2,v_3)\big) $ \\
    (C) & $ \sign (v_1,v_2) $ & $ \sign (v_1,v_3) $ & $ \sign (v_1,v_4) $
          \\ \hline
  \end {tabular} \end {center}
  Note that these $H$-signs follow a special pattern: for each of the vertices
  $V$ and $W$ that is of type (2), (3), or (4) we get a factor of $ \sign
  (v_i,v_j) $ in the $H$-sign of the resolution, where $ (i,j) \in \{ (1,2),
  (1,3), (1,4), (3,4), (4,2), (2,3) \} $ is the unique pair such that the $ v_i
  $ and $ v_j $ edges are adjacent to the vertex. On the other hand, by
  definition \ref {def-vertex} the multiplicity of such a vertex is $1$ in type
  (1), $ i^{|(v_i,v_j)|-1} $ in types (2) and (6), and $ |(v_i,v_j)| \cdot
  i^{|(v_i,v_j)|-1} $ in types (3), (4), and (5). If one replaces $ |(v_i,v_j)|
  $ by $ -|(v_i,v_j)| $ in these expressions, the vertex multiplicities remain
  the same for the types (1), (5) and (6), and are replaced by their negatives
  for the types (2), (3), and (4). It follows that the $H$-sign can be taken
  care of by replacing $ a=|(v_i,v_j)| $ by $ (v_i,v_j) $ in the vertex
  multiplicities of definition \ref {def-vertex} for $V$ and $W$.

  More precisely, if $\sigma$ denotes the $H$-sign and $m$ the multiplicity
  of a curve in a given resolution, then $ \sigma \, m = \lambda \, \tilde m_V
  \, \tilde m_W $, where $ \tilde m_V $ and $ \tilde m_W $ are the
  multiplicities of the vertices $V$ and $W$ as in definition \ref {def-vertex}
  with $a$ replaced by $ (v_i,v_j) $ as above, and $ \lambda $ is the product
  of the vertex multiplicities of all other vertices. To show that the sum of
  these numbers over all resolutions is zero we can obviously divide by the
  constant $ \lambda $ (which is the same for the resolutions I, II, III) and
  only consider $ \tilde m_V \, \tilde m_W $. Let us split this number as $
  \tilde m_V \, \tilde m_W = \mu \, m $, where $ \mu $ collects all factors $
  i^{(v_i,v_j)-1} $ and $m$ all factors $ (v_i,v_j) $ for $V$ and $W$. The
  values for $ m=m_{\text {I}},m_{\text {II}},m_{\text {III}} $ are listed in
  the table of resolutions above. As for $ \mu $, note that this number is
  \begin {itemize}
  \item in case (A): $ \mu_{\text {I}} := i^{(v_1,v_2)-1} $ for I and $
    \mu_{\text {II}} := i^{(v_1,v_3)-1} $ for II;
  \item in cases (B) and (C): $ \mu_{\text {I}} := i^{(v_1,v_2)+(v_3,v_4)-2} $
    for I, $ \mu_{\text {II}} := i^{(v_1,v_3)+(v_4,v_2)-2} $ for II, and $
    \mu_{\text {III}} := i^{(v_1,v_4)+(v_2,v_3)-2} $ for III.
  \end {itemize}
  To simplify these expressions we divide them by $ \mu_{\text {I}} $ and get
  (using $ v_1+v_2+v_3=0 $ in (A) and $ v_1+v_2+v_3+v_4=0 $ in (B) and (C))
  \begin {itemize}
  \item in case (A): $ \mu_{\text {II}} / \mu_{\text {I}} =
  i^{2(v_2,v_1)} = (-1)^{(v_2,v_1)} $;
  \item in cases (B) and (C): $ \mu_{\text {II}} / \mu_{\text {I}} =
    i^{2(v_2,v_1)} = (-1)^{(v_2,v_1)} $ and $ \mu_{\text {III}} / \mu_{\text
    {I}} = i^{2(v_1,v_4)} = (-1)^{(v_1,v_4)} $.
  \end {itemize}
  The values for these quotients are also listed in the table of resolutions
  above. Using these values for the quotients and $ m_{\text {I}}, m_{\text
  {II}}, m_{\text {III}} $, one can now easily check the required statement
    \[   \mu_{\text {I}} \cdot m_{\text {I}}
       + \mu_{\text {II}} \cdot m_{\text {II}}
       + \mu_{\text {III}} \cdot m_{\text {III}}
       = \mu_{\text {I}} \cdot \big( m_{\text {I}}
       + \mu_{\text {II}} / \mu_{\text {I}} \cdot m_{\text {II}}
       + \mu_{\text {III}} / \mu_{\text {I}} \cdot m_{\text {III}} \big)
       = 0 \]
  in all 18 codimension-1 cases, using the identities
  \begin {itemize}
  \item $ (v_1,v_2)+(v_1,v_3)=0 $ for (A),
  \item $ (v_1,v_2)+(v_4,v_2)+(v_3,v_2)=0 $, $ (v_1,v_2)(v_3,v_4)+(v_1,v_3)
    (v_4,v_2)+(v_1,v_4)(v_2,v_3) = 0 $, and $ (v_1,v_2)+(v_1,v_3)+(v_1,v_4) = 0
    $ for (B) and (C),
  \end {itemize}
  that follow from $ v_1+v_2+v_3=0 $ and $ v_1+v_2+v_3+v_4=0 $, respectively.
\end {proof}

  \section {Welschinger curves} \label {sec-welschinger}

In this section we define tropical curves that we call Welschinger curves.
Their count (for certain choices of $\Delta$) yields Welschinger invariants,
i.e.\ numbers of real rational curves on a toric Del Pezzo surface $\Sigma$
belonging to an ample linear system $D$ and passing through a generic
conjugation invariant set of $-K_{\Sigma}\cdot D-1$ points, weighted with $\pm
1$, depending on the nodes of the curve.

As we have mentioned already in the introduction, we will parametrize even
non-fixed unmarked ends of Welschinger curves as two ends of half the weight
--- this way we can avoid this kind of splitting on the bridges of section \ref
{sec-bridge}. We will refer to such ends, i.e.\ pairs of non-fixed ends of the
same odd direction adjacent to the same 4-valent vertex, as double ends. In the
following, we will first settle how to deal with these double ends. Then we
define oriented and unoriented Welschinger curves and prove that they are
equivalent. We relate unoriented Welschinger curves to tropical curves in other
literature that are counted to determine Welschinger invariants, cite the
Correspondence Theorem, and discuss some invariance and non-invariance
properties of tropical Welschinger numbers. 

\begin{definition}[Double ends and end-gluing] \label{def-endgluing}
  Let $\alpha$ be a combinatorial type of $\Mrs(\Delta)$ with $\Delta =
  (v_1,\dots,v_n)$, and let $F\subset \{1,\ldots,n\}$ be a set of fixed ends.
  Assume that there are exactly $k$ pairs $i_1<j_1,\ldots,i_k<j_k$ in
  $\{1,\ldots,n\}\setminus F$ such that the unmarked ends $y_{i_l}$ and
  $y_{j_l}$ have the same odd direction and are adjacent to the same $4$-valent
  vertex, for all $l=1,\ldots,k$. We refer in the following to such a pair of
  ends as a \df {double end}. We then set
    \[ \Delta' = \big( ( v(y_i) : i \neq i_1,j_1,\dots,i_k,j_k),
       (2\cdot v(y_{i_1}),\ldots, 2\cdot v(y_{i_k})) \big). \]
  Moreover, we define $\alpha'$ by gluing each pair of double ends $y_{i_l}$
  and $y_{j_l}$ to one unmarked end of direction $2\cdot v(y_{i_l})$, and
  denote by $ F' \subset \{1,\dots,n-k\} $ the set of entries of $ \Delta' $
  corresponding to the fixed ends $F$ in $ \Delta $. There is then an
  associated map $\Mrs^\alpha(\Delta)\rightarrow \Mrs^{\alpha'}(\Delta')$
  which we call the \df {end-gluing map}.
\end{definition}

The analogous end-gluing map $ (\Mrs^\ori)^\alpha (\Delta,F) \rightarrow
(\Mrs^\ori)^{\alpha'} (\Delta',F') $  also exists for oriented curves. The map
sending a combinatorial type $\alpha$ of $\Mrs(\Delta)$ as above to $\alpha'$
is injective, because if we want to produce a preimage $\alpha$ from $\alpha'$,
we just have to split the last $k$ ends of $ \Delta' $, producing $4$-valent
vertices.

\begin{example}\label{ex-endgluing}
  The following picture shows a curve $C$ and its image $C'$ under the
  end-gluing map. Although mainly following convention \ref{conv-curve}, we
  draw double ends separately even though this is actually a feature of the
  graph $\Gamma$ and cannot be seen in $h(\Gamma)$.
  \begin {center} \input {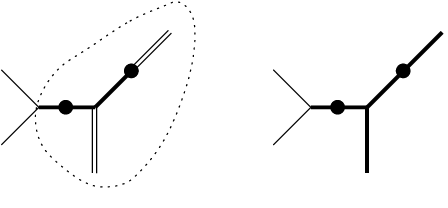} \end {center}
\end{example}

\begin{remark}\label{rem-dimaftergluing}
  It follows from example \ref{ex-general} that if a collection of conditions $
  \calP \in \RR^{2(r+s)+|F|} $ as in remark \ref {rem-ev} is in general
  position for $ \ev_F :\Mrs^\alpha(\Delta) \to \RR^{2(r+s)+|F|} $ then it
  is also in general position after end-gluing for $\ev_{F'} :\Mrs^{\alpha'}
  (\Delta') \to \RR^{2(r+s)+|F|} $, and vice versa. Notice also that 
  $ \dim \Mrs^\alpha(\Delta) = \dim \Mrs^{\alpha'}(\Delta') $: by \cite
  {GM05b} proposition 2.11 a combinatorial type has dimension $
  |\Delta|-1+r+s-\sum_V(\val(V)-3)$ where the sum goes over all vertices $V$ of
  $\Gamma$, and the end-gluing map decreases the number of entries of $\Delta$
  by the same number as it decreases the number of $4$-valent vertices. As
  orienting the edges does not change dimensions we conclude that the
  end-gluing map does not change the dimension of combinatorial types of
  oriented curves either.
\end{remark}

\begin{definition}[$ \Gamma_\even $ and roots] \label{def-gammaeven}
  Let $ C = (\Gamma,x_1,\dots,x_{r+s},y_1,\dots,y_n,h) \in \Mrs(\Delta) $.
  Let $C'$ be the image of $C$ under the end-gluing map of definition \ref
  {def-endgluing} and call the underlying graph $\Gamma'$. Consider the
  subgraph $ \Gamma'_\even $ of $ \Gamma' $ of all even edges (including the
  markings), and its preimage $\Gamma_\even$. That is, $\Gamma_\even$ consists
  of all even edges and all double ends of $\Gamma$. Vertices of $ \Gamma_\even
  \cap \overline {\Gamma \backslash \Gamma_\even} $ as well as unmarked
  non-fixed \emph{even} ends of $\Gamma_\even$ are called the \df {roots} of
  $\Gamma_\even$.
\end{definition}

\begin{example}
  For the curve of example \ref{ex-endgluing}, the part $\Gamma_\even$ is
  encircled. It has one root, namely the vertex denoted by $V$.
\end{example}

\begin {definition}[Welschinger curves] \label {def-welschinger}
  Let $ r,s \ge 0 $, let $ \Delta = (v_1,\dots,v_n) $ be a collection of
  vectors in $ \ZZ^2 \backslash \{0\} $, and let $ F \subset \{1,\dots,n\} $.
  \begin {enumerate}
  \item \label {def-welschinger-a}
    An oriented curve $ C \in \Mrs^\ori(\Delta,F) $ all of whose vertices are
    of the types (1) to (5), (6b), (7), or (8) of definition \ref {def-vertex}
    is called an \df {oriented Welschinger curve}. 
  \item \label {def-welschinger-b}
    Let $ C = (\Gamma,x_1,\dots,x_{r+s},y_1,\dots,y_n,h) \in \Mrs(\Delta) $,
    and let $\Gamma_\even$ be as in definition \ref{def-gammaeven}. We say that
    $C$ is an \df {unoriented Welschinger curve} (with set of fixed ends $F$)
    if
    \begin {enumerate}
    \item [(i)] complex markings are adjacent to 4-valent vertices, or
      non-isolated in $\Gamma_\even$; 
    \item [(ii)] each connected component of $\Gamma_\even$ has a unique root.
    \end {enumerate}
  \end {enumerate}
\end {definition}

\begin{example}
  The following picture shows an oriented Welschinger curve with an even and an
  odd fixed end. As in example \ref{ex-endgluing}, we indicate double ends in
  the picture while otherwise following convention \ref{conv-curve}. Each
  vertex is labeled with its type, every allowed vertex type occurs. If we
  forget the orientations of the edges, we get an unoriented Welschinger curve.
  There are four connected components of $\Gamma_\even$. The subgraph
  $\Gamma_3$ consists of a complex marking and $\Gamma_4$ of a real marking.
  $\Gamma_1$ and $\Gamma_2$ both have one root, namely the vertex of type (3).
  Three complex markings are adjacent to $4$-valent vertices, four are
  non-isolated in $\Gamma_\even$. 

  \begin {center} \input {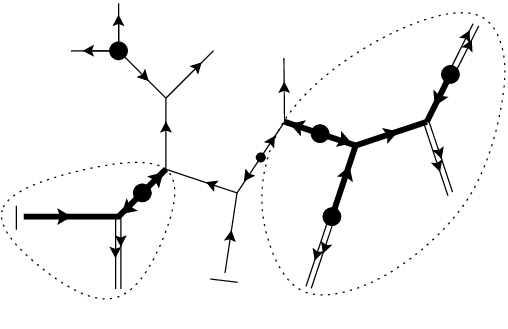} \end {center}
\end{example}

As for broccoli curves, we want to show that oriented and unoriented
Welschinger curves are equivalent for enumerative purposes. The following
remark and lemma are needed as preparation.  

\begin {remark} \label {rem-outends}
  Let $ C \in \Mrs^\ori(\Delta,F) $ be an oriented Welschinger curve. 
  \begin {enumerate}
  \item \label {rem-outends-a}
    By lemma \ref {lem-dim}, the curve $C$ has $ |\Delta|-|F| = r+2s+1-\n7+\n8
    $ outward pointing ends. In particular, if $ |\Delta|-1 = r+2s+|F|$ then $
    \n7=\n8 $.
  \item \label {rem-outends-b}
    If $C$ consists only of vertices of types (4), (6b), (7) and (8), then
    we have $ r=0 $, $ s=\n{6b}+\n7 $, and the number of odd outward pointing
    ends is $ 2\n{6b}+2\n8 $. Hence in this case it follows from \ref
    {rem-outends-a} that $C$ has exactly $ 1+\n7-\n8 $ even outward pointing
    ends.
  \end {enumerate}
\end {remark}

\begin {lemma} \label{lem-forest}
  Let $|\Delta|-1=r+2s+|F|$, let $ C\in \Mrs^\ori(\Delta,F) $ be an oriented
  Welschinger curve, and let $\Gamma_\even$ be as in definition \ref
  {def-gammaeven}. Then every connected component of $ \, \Gamma_\even$ has
  exactly one root.
\end {lemma}

\begin {proof}
  If $\Gamma_\even=\Gamma$ then $\Gamma$ has only vertices of type (4), (6b),
  (7), and (8). By remark \ref{rem-outends} \ref{rem-outends-a} we have
  $ \n7=\n8 $, so from remark \ref{rem-outends} \ref{rem-outends-b} it then
  follows that $\Gamma$ has exactly one even outward pointing end, which is the
  unique root.
  
  If $\Gamma_\even \neq \Gamma$, every connected component $\tilde{\Gamma}$ of
  $\Gamma_\even$ needs to be adjacent to odd edges which are not double ends.
  The only allowed vertex type for oriented Welschinger curves to which both
  even edges (resp.\ double edges) and odd edges (which are not double ends)
  are adjacent is type (3). Each vertex of type (3) yields a $1$-valent vertex
  in $\Gamma_\even$. Remove these $1$-valent vertices from the component
  $\tilde{\Gamma}$, and call the resulting graph $\tilde{\Gamma}^\circ$. A
  vertex of type (3) leads to an outward pointing end of $ \tilde{\Gamma}^\circ
  $. Note that $ \tilde\Gamma^\circ $ has vertices of types (4), (6b), (7), and
  (8). Thus by remark \ref{rem-outends} \ref{rem-outends-b} we have $
  \n8^{\tilde\Gamma^\circ} \leq \n7^{\tilde\Gamma^\circ} $, where the
  superscripts indicate that we refer to numbers of vertices of $ \tilde
  \Gamma^\circ $. By remark \ref{rem-outends} \ref{rem-outends-a} we have $
  \n7^C=\n8^C $. Since any vertex of type (7) or (8) belongs to exactly one
  graph $ \tilde\Gamma^\circ $ associated to a connected component $
  \tilde\Gamma $ of $ \Gamma_\even $, and since the inequality $ \n8^{\tilde
  \Gamma^\circ} \leq \n7^{\tilde\Gamma^\circ} $ holds for any such $
  \tilde\Gamma $, we conclude that it is an equality. Then by remark \ref
  {rem-outends} \ref {rem-outends-b} every $ \tilde\Gamma^\circ $ has exactly
  one even outward pointing end. It follows that every $ \tilde\Gamma $ has
  exactly one root.
\end{proof}

With this preparation we can prove the following statement analogously to
proposition \ref {prop-broccoli}.

\begin {proposition}[Equivalence of oriented and unoriented Welschinger curves]
    \label {prop-welschingeroriented}
  Let $ r,s \ge 0 $, let $ \Delta = (v_1,\dots,v_n) $ be a collection of
  vectors in $ \ZZ^2 \backslash \{0\} $ , and let $ F \subset \{1,\dots,n\} $
  such that $ r+2s+|F| = |\Delta|-1 $. Moreover, let $ \calP \in
  \RR^{2(r+s)+|F|} $ be a collection of conditions in general position for $
  \ev_F: \Mrs(\Delta) \to \RR^{2(r+s)+|F|} $ (see example \ref {ex-general}).

  Then the forgetful map $ \ft $ of definition \ref {def-oriented} gives a
  bijection between oriented and unoriented $ (r,s) $-marked Welschinger curves
  through $ \calP $ with degree $ \Delta $ and set of fixed ends $F$.
\end {proposition}

\begin {proof}
  As in proposition \ref {prop-broccoli}, we have to prove three statements.
  \begin {enumerate}
  \item 
    $ \ft $ maps oriented to unoriented Welschinger curves through $ \calP $:
    Let $ C \in \Mrs^\ori(\Delta,F) $ be an oriented Welschinger curve. The
    list of allowed vertex types for $C$ implies that $C$ satisfies condition
    (i) of definition \ref {def-welschinger}. Condition (ii) follows from lemma
    \ref{lem-forest}.
  \item 
    $ \ft $ is injective on the set of curves through $ \calP $: Notice that
    under the end-gluing map of definition \ref{def-endgluing}, a vertex of
    type (8) becomes a vertex of type (4), and type (6b) becomes (7). Thus the
    image $C'$ under the end-gluing map satisfies the conditions of lemma \ref
    {lem-general} \ref {lem-general-b} by lemma \ref {lem-dim} and remark \ref
    {rem-dimaftergluing}. Lemma \ref {lem-oriented} implies that there is at
    most one possible orientation on $C'$, and it follows immediately that
    there is only one possible orientation on $C$, since double ends have to
    point outwards (types (6b) and (8)).
  \item 
    $ \ft $ is surjective on the set of curves through $ \calP $: Let $ C \in
    \Mrs(\Delta) $ be an unoriented Welschinger curve through $ \calP $ with
    set of fixed ends $F$. Let $\alpha$ be the combinatorial type of $C$ and $
    \Mrs^\alpha(\Delta)$ its corresponding cell in $\Mrs(\Delta) $. Denote by
    $s_1$ the number of isolated complex markings in $\Gamma_\even$,
    and by $k$ the number of double ends. As this means by definition \ref
    {def-endgluing} and condition (i) that there are at least $ s_1+k $
    vertices of valence 4 it follows from \cite {GM05b} proposition 2.11 that
    the dimension of $\Mrs^\alpha(\Delta)$ is at most $|\Delta|+r+s-1-s_1-k =
    2r+3s+|F|-s_1-k $. On the other hand, $C$ passes through a collection of
    conditions in general position, so $\dim (\Mrs^\alpha(\Delta)) \geq
    2r+2s+|F|$. It follows that
      \[ \qquad \quad
         s-s_1-k \ge 0. \tag {$*$} \]
    In fact, we want to show that we always have equality here. For this
    let $ \tilde\Gamma $ be a connected component of $ \Gamma_\even \backslash
    \big( \overline {\Gamma \backslash \Gamma_\even} \big) $ --- i.e.\ we
    remove from $ \Gamma_\even $ all attachment vertices to its complement ---
    which is not an isolated marked end. Denote by $ \tilde\Gamma' $ its image
    under the end-gluing map. Let $ \tilde s $ be the number of complex
    markings belonging to $ \tilde \Gamma $, and let $ \tilde k $ be the number
    of its double ends. Then $ \tilde \Gamma' $ contains possibly fixed even
    ends, the $ \tilde k $ ends coming from the double ends, and one extra end
    (which is either the root itself or the edge with which it is adjacent to $
    \Gamma \backslash \Gamma_\even $). If $\tilde s > \tilde k $ it follows
    that there is a component of $ \tilde \Gamma' $ minus the $ \tilde s $
    complex markings which does not contain a non-fixed end, which would be a
    contradiction to lemma \ref {lem-general} \ref {lem-general-a}. Thus $
    \tilde s \le \tilde k $. Summing this up over all such components $
    \tilde \Gamma $ it follows that the number $ s-s_1 $ of complex markings
    which are non-isolated in $ \Gamma_\even $ satisfies $ s-s_1 \le k $.
    Together with $(*)$ this yields $s-s_1=k$, as desired.

    Hence equality holds in all our estimates above. There are various
    consequences of this: first of all, we have $ \dim (\Mrs^\alpha(\Delta)) =
    2r+2s+|F| $, and $C$ has exactly $s$ vertices of valence $4$, namely $s_1$
    adjacent to complex markings which are isolated in $\Gamma_\even$, and $
    s-s_1$ adjacent to double ends. All other vertices have valence $3$. In
    particular, if the root of a connected component of $\Gamma_\even $ is not
    an end, it has to be at a $3$-valent vertex. Also, since we have $ \tilde
    s = \tilde k $ complex markings on the components $ \tilde \Gamma $ above,
    it follows that there cannot be additional real markings on these
    components, since otherwise there would be a connected component of $
    \tilde\Gamma' $ without the markings again which does not contain a
    non-fixed end. Thus there are no real markings which are non-isolated in $
    \Gamma_\even $.

    The combinatorial type of the image $C'$ of $C$ under the end-gluing map is
    of dimension $\dim (\Mrs^\alpha(\Delta))= 2r+2s+|F|$ by remark
    \ref{rem-dimaftergluing}. Since $C$ has $4$-valent vertices only at complex
    markings resp.\ double ends, it follows that $C'$ has $4$-valent vertices
    only at complex markings, and so we can apply lemma \ref {lem-general} to
    $C'$ to see that there is an orientation on $C'$ that points on each edge
    towards the unique non-fixed unmarked end in $ \Gamma' \backslash (x_1 \cup
    \cdots \cup x_{r+s}) $. We can define an orientation on $C$ by orienting
    double ends just as the end they map to under the end-gluing map.

    It remains to be shown that, for this orientation of $C$, we only have the
    vertex types (1) to (5), (6b), (7) or (8). As in the proof of proposition
    \ref{prop-broccoli} \ref{prop-broccoli-c}, all edges adjacent to a vertex
    $V$ point outwards if there is a marking at $V$, and exactly one points
    outwards otherwise. It is impossible that exactly one edge at $V$ is odd.
    We have seen that $V$ is $4$-valent if it is adjacent to a double end, or
    to a complex marking, and $3$-valent otherwise. The only vertex types
    compatible with all these restrictions are the types (1) to (8), and the
    three special ones in the picture of the proof of proposition \ref
    {prop-broccoli} \ref {prop-broccoli-c}. Type (6a) cannot appear since each
    root has to be $3$-valent by the above. The left picture in the proof of
    proposition \ref {prop-broccoli} \ref {prop-broccoli-c} is excluded since
    there are no non-isolated real markings in $\Gamma_\even$. The middle
    picture is excluded since we have $4$-valent vertices only at isolated
    complex markings or double ends. The right picture would be a root of a
    component $ \tilde\Gamma $ as above. But because of the orientation there
    is no connection from this vertex via one of the odd edges to a non-fixed
    unmarked end without passing a marking. With $ \tilde k $ non-fixed ends
    and $ \tilde k $ complex markings in $ \tilde\Gamma $ this would again
    lead to a connected component of $ \Gamma $ minus the markings with no
    non-fixed end, a contradiction to lemma \ref {lem-general} \ref
    {lem-general-a}. \qedhere
  \end {enumerate}
\end {proof}

\begin {remark}[Unoriented Welschinger curves after end-gluing]
    \label {rem-welschingerunoriented}
  In addition to definition \ref {def-welschinger} \ref {def-welschinger-b} we
  can also describe unoriented Welschinger curves after the end-gluing: fix a
  degree $\Delta = (v_1,\dots,v_n)$ and $F\subset \{1,\ldots,n\}$. We
  then allow curves of any degree $\Delta' = \big( ( v(y_i) : i \neq
  i_1,j_1,\dots,i_k,j_k), (2\cdot v(y_{i_1}),\ldots, 2\cdot v(y_{i_k})) \big)$ for
  some $i_1<j_1,\ldots,i_k<j_k$ in $\{1,\ldots,n\}\setminus F$ such that the
  unmarked ends $y_{i_l}$ and $y_{j_l}$ have the same odd direction. For a
  curve $ C = (\Gamma,x_1,\dots,x_{r+s},y_1,\dots,y_{n-k},h) \in \Mrs(\Delta')
  $, we define $\Gamma_\even$ as in definition \ref{def-broccoli} as the
  subgraph of all even edges. We then require that complex markings are
  adjacent to 4-valent vertices, or non-isolated in $\Gamma_\even$; and that
  each connected component of $\Gamma_\even$ has a unique root. An example of
  such an unoriented Welschinger curve after end-gluing is the top left curve
  in the bridge picture in the introduction.
\end{remark}

Now we define enumerative numbers of Welschinger curves. As for broccoli
curves, we work with oriented Welschinger curves from now on, keeping in mind
that it does not matter whether we count oriented or unoriented Welschinger
curves by proposition \ref {prop-welschingeroriented}.

\begin {notation} \label {not-welschinger}
  Let $ r+2s+|F| = |\Delta|-1 $, and denote by $ \Mrs^W(\Delta,F) $ the
  closure of the space of all Welschinger curves in $ \Mrs^\ori(\Delta,F) $.
  This is obviously a polyhedral subcomplex. By lemma \ref {lem-dim} it is of
  pure dimension $ 2(r+s)+|F| $, and its maximal open cells correspond exactly
  to the Welschinger curves in $ \Mrs^W(\Delta,F) $. For $ F=\emptyset $ we
  write $ \Mrs^W(\Delta,F) $ also as $ \Mrs^W(\Delta) $.
\end {notation}

\begin {definition}[Welschinger numbers] \label {def-welschinger-inv}
  Let $ r,s \ge 0 $, let $ \Delta = (v_1,\dots,v_n) $ be a collection of
  vectors in $ \ZZ^2 \backslash \{0\} $, and let $ F \subset \{1,\dots,n\} $
  such that $ r+2s+|F| = |\Delta|-1 $. Moreover, let $ \calP \in
  \RR^{2(r+s)+|F|} $ be a collection of conditions in general position for
  Welschinger curves, i.e.\ for the evaluation map $ \ev_F: \Mrs^W(\Delta,F)
  \to \RR^{2(r+s)+|F|} $. Then we define the \df {Welschinger number}
    \[ \Nrs^W (\Delta,F,\calP) \;:=\; \frac 1{|G(\Delta,F)|} \cdot
         \sum_C m_C, \]
  where the sum is taken over all Welschinger curves $C$ in with degree $
  \Delta $, set of fixed ends $F$, and $ \ev(C) = \calP $. As in the case of
  broccoli invariants, the group $ G(\Delta,F) $ compensates for the
  overcounting of curves due to relabeling the non-fixed unmarked ends (see
  remark \ref {rem-endlabeling}). The sum is finite by the dimension statement
  of notation \ref {not-welschinger}, and the multiplicity $ m_C $ is as in
  definition \ref {def-vertex}. For $ F=\emptyset $ we abbreviate the numbers
  as $ \Nrs^W (\Delta,\calP) $.
\end {definition}

In contrast to the broccoli invariants of definition \ref {def-broccoli-inv} we
will see in remark \ref {rem-notinvariant} that these Welschinger numbers will
in general depend on the choice of conditions $ \calP $. For $ F=\emptyset $
and certain choices of the degree $ \Delta $ however, there exist well-known
Welschinger invariants in the literature that count real rational algebraic
curves through given points in the plane, and that do not depend on the choice
of point conditions. We want to show now that they agree with our Welschinger
numbers in these cases.

\begin{remark}[Welschinger curves compared to \cite {Shu06}]
    \label{rem-shustin}
  Notice that (unoriented) Welschinger curves where all unmarked ends are
  non-fixed and odd correspond precisely to the curves considered by Shustin in
  \cite {Shu06} (in the way described in remark \ref
  {rem-welschingerunoriented}). There, unparametrized tropical curves are
  considered, i.e.\ the images $ h(\Gamma) $ without the parametrizing graph $
  \Gamma $, and it is required that the point conditions are general enough so
  that the Newton subdivision dual to $h(\Gamma)$ (see \cite{Mik03} proposition
  3.11) consists only of triangles and parallelograms. In this case each such
  unparametrized curve can uniquely (up to the labeling of the unmarked ends)
  be parametrized by a graph $ \Gamma' $ such that the map to $ \RR^2 $
  identifies only finitely many points. Adding an end for each marking and
  reversing the end-gluing by splitting each even unmarked end into a double
  end then gives a graph $ \Gamma $ together with a map $ h: \Gamma \to \RR^2 $
  satisfying the conditions of definition \ref {def-welschinger} \ref
  {def-welschinger-b}. The part $h(\Gamma_\even)$ coincides with the subgraph
  $G$ in \cite{Shu06} consisting of all the even edges; their components are
  connected to odd edges at exactly one vertex, the root. (Other authors
  consider parametrized curves and the even part $G$ as the non-fixed locus of
  a certain involution on the tropical curve, from which it also follows that
  each connected component has one root \cite{BM08}.)
  
  The definition of the multiplicities of these curves in \cite {Shu06} looks
  at first a little different compared to our definition \ref {def-vertex}. We
  will recall it here and then show that it in fact coincides with ours.
\end{remark}

\begin{definition}[W-multiplicity, see \cite {Shu06} section 2.5]
    \label{def-multshu}
  Let $C=(\Gamma,x_1,\ldots,x_{r+s},y_1,\ldots,y_n,h)$ be a Welschinger curve,
  and assume that the Newton subdivision dual to $h(\Gamma)$ (see \cite{Mik03}
  proposition 3.11) consists of only triangles and parallelograms. Denote by
  $\tilde{a}$ the number of lattice points inside triangles of this
  subdivision, by $\tilde{b}$ the number of triangles such that all sides have
  even lattice length, and by $\tilde{c}$ the number of triangles whose lattice
  area is even. Then we define the \df {W-multiplicity} of $C$ to be
    \[ \tilde m_C := (-1)^{\tilde{a}+\tilde{b}} \cdot 2^{-\tilde{c}} \cdot
       \prod_V \mult(V), \]
  where the product goes over all triangles with even lattice area or dual to
  vertices with a complex marking, and where $ \mult(V) $ denotes the integer
  area of this triangle, i.e.\ the complex vertex multiplicity as in definition
  \ref {def-vertex}.

  For an unparametrized curve $h(\Gamma)$, this coincides with the definition
  of multiplicity in \cite{Shu06} section 2.5.
\end{definition}

\begin {example}
  The following picture shows a Welschinger curve (without orientation) and its
  dual Newton subdivision. The triangles $V$ contributing to $ \tilde m_C $ are
  shaded and labeled with their integer area; we have $ \tilde m_C = (-1)^{1+1}
  \cdot 2^{-2} \cdot 4 \cdot 2 \cdot 3 \cdot 1 = 6 $.

  \begin {center} \input {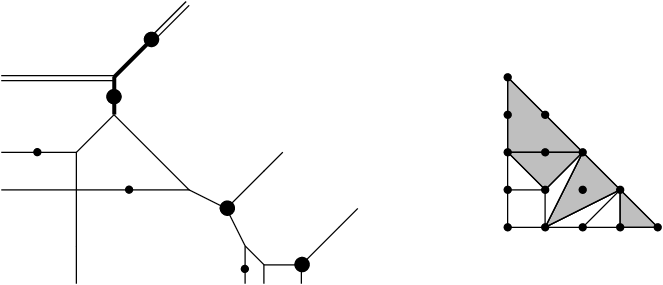} \end {center}
\end {example}

\begin {remark}[Labeled and unlabeled curves] \label {rem-endlabeling}
  Note that we consider curves with labeled unmarked ends, whereas the
  unparametrized curves in \cite{Shu06} come without this data. Thus we
  overcount each unparametrized curve by a factor that records the different
  ways to label the (non-fixed) unmarked ends so that we get different
  parametrized curves. If $k$ denotes the number of double ends then this
  overcounting factor is $ |G(\Delta)| \cdot 2^{-k} $, where the $ 2^{-k} $
  term arises because exchanging the two labels of a double end does not change
  the parametrized curve.
\end {remark}

\begin{lemma}[Multiplicity and W-multiplicity] \label{lem-multequal}
  Let $C=(\Gamma,x_1,\ldots,x_{r+s},y_1,\ldots,y_n,h)$ be a Welschinger curve
  of degree $\Delta$ with no fixed ends, satisfying $\omega(y_i)=1$ for all
  $i=1,\ldots,n$, and passing through points in general position as in example
  \ref{ex-general}. Then the multiplicity $ m_C $ and the W-multiplicity $
  \tilde m_C $ of $C$ are related by $ m_C = 2^k \cdot \tilde m_C $, where $k$
  is the number of double ends of $C$.
\end{lemma}

\begin{proof}
  It follows from the list of allowed vertex types and their multiplicities
  that a vertex $V$ contributes a factor of $ \mult(V) $ to $m_C$ if and only
  if $V$ is adjacent to a complex marking or dual to a triangle with even
  lattice area.

  The number $\tilde{c}$ of triangles with even lattice area equals
  $\n3+\n4+\n8$. Let $\tilde{\Gamma}$ be a connected component of
  $\Gamma_\even$. We know that $\tilde{\Gamma}$ has a unique root. Since
  $\omega(y_i)=1$ for all $ i=1,\ldots,n$, this root cannot be an end of
  $\Gamma$, so it has to be a vertex of type (3) in $\Gamma$, i.e.\ a
  $1$-valent vertex in $\Gamma_\even$. Remove the $1$-valent vertex from
  $\tilde{\Gamma}$, thus producing an end, apply the end-gluing map of
  definition \ref{def-endgluing}, and forget all markings (straightening the
  $2$-valent vertices). Call the resulting graph $\tilde{\Gamma}^\circ$. This
  graph is $3$-valent and has $ 1+\n{6b}^{
  \tilde{\Gamma}}+\n8^{\tilde{\Gamma}}$ ends, and thus it has $
  \n{6b}^{\tilde{\Gamma}}+\n8^{\tilde{\Gamma}}-1$ vertices. But this number of
  $3$-valent vertices also equals $\n4^{\tilde{\Gamma}}+\n8^{\tilde{\Gamma}}$,
  and so $\n{6b}^{\tilde{\Gamma}}+\n8^{\tilde{\Gamma}}=\n4^{\tilde{\Gamma}}+
  \n8^{\tilde{\Gamma}}+1=\n4^{\tilde{\Gamma}}+\n8^{\tilde{\Gamma}}+
  \n3^{\tilde{\Gamma}}$. Since this holds for any $\tilde{\Gamma}$, it follows
  that $\n{6b}+\n8=\n3+\n4+\n8$. Thus $k=\tilde{c}$, where $k$ denotes the
  number of double ends. The factor $ 2^k $ in the lemma thus corresponds
  exactly to the factor $ 2^{-\tilde{c}} $ in the definition \ref {def-multshu}
  of $ \tilde m_C $.

  Hence it only remains to show that $(-1)^{\tilde{a}+\tilde{b}}$ equals the
  sign contribution coming from factors of $i$ in the definition \ref
  {def-vertex} of $m_C$, where $\tilde{a}$ denotes the number of lattice points
  in the interior of triangles and $\tilde{b}$ denotes the number of triangles
  such that all sides have even lattice length. We refer to the power of $i$ in
  the vertex multiplicity $m_V$ of definition \ref{def-vertex} as the sign.

  Consider a vertex $V$ and let $A=\mult(V)$. If $V$ is of type (2) to (5),
  assume the three adjacent (non-marked) edges have weights $\omega_1$,
  $\omega_2$ and $\omega_3$. By Pick's formula, $A=2I+B-2$, where $I$ denotes
  the number of lattice points in the interior of the triangle dual to $V$ and
  $B$ denotes the number of lattice points on the boundary. By our assumptions,
  $B=\omega_1+\omega_2+\omega_3$. If $V$ is of type (2) or (5), then its sign
  is
    \[ i^{A-1} = (-1)^{\frac{A-1}{2}}
               = (-1)^{\frac{2I+\omega_1+\omega_2+\omega_3-2-1}{2}}
               = (-1)^I\cdot (-1)^{\frac{\omega_1-1}{2}}\cdot
                 (-1)^{\frac {\omega_2-1}2}\cdot (-1)^{\frac{\omega_3-1}2}. \]
  If $V$ is of type (3), its sign is
  \begin{align*}
    i^{A-1} &= i^{-1}\cdot i^A
             = i^{-1}\cdot (-1)^{\frac{A}{2}}
             = i^{-1}\cdot (-1)^{\frac{2I+\omega_1+\omega_2+\omega_3-2}{2}} \\
            &= i^{-1}\cdot (-1)^I \cdot (-1)^{\frac{\omega_1-1}{2}} \cdot
               (-1)^{\frac{\omega_2-1}{2}}\cdot
               (-1)^{\frac{\omega_3}{2}},
  \end{align*}
  where we assume that $\omega_3$ is the even weight. For type (4), we get
  \begin{align*}
    i^{A-1} &= i^{-1}\cdot i^A
             = i^{-1}\cdot (-1)^{\frac{A}{2}}
             = i^{-1}\cdot (-1)\cdot
               (-1)^{\frac{2I+\omega_1+\omega_2+\omega_3}{2}} \\
            &= i^{-1}\cdot (-1)\cdot (-1)^I \cdot
               (-1)^{\frac{\omega_1}{2}}\cdot (-1)^{\frac{\omega_2}{2}}\cdot
               (-1)^{\frac{\omega_3}{2}}.
  \end{align*}
  We write the sign of type (6b) as $i^{-1}=i\cdot (-1)=i\cdot
  (-1)^{\frac{2}{2}}$, and $2$ is the weight of the even adjacent edge (since
  the double ends are of weight $1$ by assumption). The sign of (8) is 
    \[ -1 = (-1)\cdot i^{A}
          = (-1)\cdot (-1)^I \cdot (-1)^{\frac{\omega_1}{2}}\cdot
            (-1)^{\frac{\omega_2}{2}}, \]
  where $\omega_1$ and $\omega_2$ are the weights of the two adjacent even
  edges. This is true since the two edges of the same direction which are
  adjacent to (8) are ends and thus their weight is $1$ by assumption. The sign
  of (1) can be written as $1=(-1)^{\frac{\omega_1-1}{2}}\cdot
  (-1)^{\frac{\omega_2-1}{2}}$, where $\omega_1=\omega_2$ is the odd weight of
  the adjacent edges. Analogously, we can write the sign of (7) as
  $1=(-1)^{\frac{\omega_1}{2}}\cdot (-1)^{\frac{\omega_2}{2}}$, where now
  $\omega_1=\omega_2$ is the even weight of the adjacent edges.

  Notice that the product of the factors $(-1)^I$ which appear for each vertex
  dual to a triangle is $(-1)^{\tilde{a}}$. Also, for each vertex of type (4)
  and (8) --- which are the vertices dual to triangles such that all sides have
  even lattice length --- we have a factor of $(-1)$ which yields
  $(-1)^{\tilde{b}}$ as product. In addition, we have extra factors of $i^{-1}$
  for each vertex of type (3) and (4), and $i$ for each vertex of type (6b).
  But since $\n4+\n3=\n{6b}$ as we have seen above, these extra factors cancel.
  Furthermore, we have factors of $(-1)^{\frac{\omega-1}{2}}$ for each edge of
  odd weight ending at a vertex, and $(-1)^{\frac{\omega}{2}}$ for each even
  edge. Every bounded edge ends at two vertices, so these contributions cancel.
  Since we require that the weights of all ends are $1$, the corresponding
  factors for the ends are just $1$. Thus all the factors $ (-1)^{
  \frac{\omega-1}{2}}$ resp.\ $(-1)^{\frac{\omega}{2}}$ cancel, and it follows
  that the sign equals $(-1)^{\tilde{a}+\tilde{b}}$, as required.
\end{proof}

\begin{remark}[Welschinger numbers compared to \cite {Shu06}]
    \label{rem-unparametrizednumbers}
  It follows from remark \ref {rem-shustin}, remark \ref {rem-endlabeling}, and
  lemma \ref {lem-multequal} that for $F=\emptyset$ and $\Delta$ consisting of
  primitive vectors (i.e.\ of directions of weight one) our Welschinger number
  $\Nrs^W (\Delta,\calP)$ of definition \ref {def-welschinger-inv} equals the
  number of unparametrized curves as in \cite{Shu06}, counted with their
  W-multiplicities as in definition \ref{def-multshu}.
\end{remark}

\begin {example}[Welschinger numbers for degrees with non-fixed even ends]
    \label {ex-welschinger-even}
  In two special cases when the degree $ \Delta=(v_1,\dots,v_n) $ contains one
  or several non-fixed even ends we can actually compute the Welschinger
  numbers directly:
  \begin {enumerate}
  \item \label {ex-welschinger-even-a}
    Assume that $\Delta$ contains more than one non-fixed even end.

    Consider a Welschinger curve $ C = (\Gamma,x_1,\dots,x_{r+s},y_1,\dots,y_n,
    h) $ contributing to the number $ \Nrs^W (\Delta,F,\calP)$. Every even
    non-fixed end belongs to a connected component of $\Gamma_\even$ and is a
    root. Since every connected component has a unique root by definition \ref
    {def-welschinger} \ref {def-welschinger-b} (ii) it follows that such a
    component cannot meet the remaining part $ \overline{\Gamma \setminus
    \Gamma_\even}$. But as the curve is connected this means that $
    \Gamma_\even $ can have only one connected component and thus only one
    root. This is a contradiction, showing that there is no Welschinger curve
    with more than one non-fixed even end, and thus that in this case
      \[ \qquad \quad \Nrs^W (\Delta,F,\calP)=0. \]
  \item \label {ex-welschinger-even-b}
    Assume now that $\Delta$ contains exactly one non-fixed end of weight $2$,
    of direction $v_1$, and only non-fixed edges of weight $1$ otherwise.

    Assume that $ \Nrs^W (\Delta,\calP) \neq 0$. By the same argument as in
    \ref {ex-welschinger-even-a} each curve contributing to $ \Nrs^W
    (\Delta,\calP) $ is totally even (containing one even and $
    \frac{|\Delta|-1}{2}$ double ends). Hence $|\Delta|$ must be odd and must
    contain each vector $v_i$ ($i\neq 1$) twice. Without restriction we
    can assume that $v_i=v_{i+\frac{|\Delta|-1}{2}}$ for $ 1 < i \leq \frac
    {|\Delta|-1}{2}+1 $. Furthermore, it then follows that $r=0$ and $ s =
    \frac{|\Delta|-1}{2}$.

    In other words, each curve contributing to $N_{(0,s)}^W (\Delta,\calP) $
    contains only vertices of type (4), (6b), (7), and (8). We can thus
    interpret the number $N_{(0,s)}^W (\Delta,\calP) $ as a ``double complex
    enumerative number'' in the following sense: let $ \Delta' =
    (\frac{1}{2}v_1,v_2,\ldots,v_{\frac{|\Delta|-1}{2}+1})$ and denote by $
    N_s^C (\Delta',\calP) $ the number of ($3$-valent) tropical curves (without
    labeled ends) passing through $\calP$ as e.g.\ in \cite {GM05a}, i.e.\ each
    curve is counted with its usual complex multiplicity as in \cite{Mik03}. If
    we forget the labels of the non-marked ends, the set of curves contributing
    to $ N_{(0,s)}^W (\Delta,\calP) $ is then obviously in bijection to the set
    of curves contributing to $ N_s^C (\Delta',\calP) $ by multiplying each
    direction vector (after end-gluing) with $\frac{1}{2}$. However, $
    N_{(0,s)}^W (\Delta,\calP) $ is not quite equal to $ N_s^C (\Delta',\calP)
    $ since the multiplicities of the curves are slightly different:
    \begin {itemize}
    \item If the vector $ \frac 12 v_1 $ occurs $d$ times in $ \Delta' $ then
      there are $d$ choices in the count of $ N_{(0,s)}^W (\Delta,\calP) $
      which of the ends of the ``double complex curve'' is the weight-2 end of
      the Welschinger curve.
    \item As we count Welschinger curves with labeled ends to get the
      number $N_{(0,s)}^W (\Delta,\calP) $, we overcount each curve without
      labeled ends by a factor of $ |G(\Delta)| \cdot 2^{-\frac{|\Delta|-1}{2}}
      $ (see remark \ref{rem-endlabeling}), since $\frac{|\Delta|-1}{2}$ is the
      number of double ends.
    \item Under the bijection, each vertex of type of type (4) and (8) maps to
      a vertex of complex multiplicity $\frac{a}{4}$. Denote by $\Gamma'$ the
      graph after end-gluing and forgetting the marked points. This graph has
      $\frac{|\Delta|-1}{2}+1$ ends and is $3$-valent, thus we have $ \n4+\n8
      = \frac{|\Delta|-1}{2}-1$. Therefore we overcount each Welschinger curve
      by an additional factor of $4^{\frac{|\Delta|-1}{2}-1}$.
    \item In addition, we count each Welschinger curve with a sign, namely $ 
      i\cdot (-1)^{\n8} \cdot i^{-\n4-\n{6b}}$, where the factor of $i$ arises
      because of the end of weight $2$ and the other factors arise because of
      the vertex multiplicities. The number of ends of the graph $\Gamma'$
      equals $ \n{6b}+\n7+1 = \frac{|\Delta|-1}{2}+1$, thus we have $\n4+\n8+1
      = \n{6b}+\n7$. Since $\n7=\n8$ by \ref{rem-outends}, we can conclude $
      \n4+1 = \n{6b}$, thus the sign above equals $(-1)^{\n8}\cdot i^{-2\n4} =
      (-1)^{\n4+\n8} = (-1)^{\frac{|\Delta|-1}{2}-1} $.
    \end {itemize}
    Taking all these factors together, it follows that
    \begin {align*} \qquad \quad
      N_{(0,s)}^W (\Delta,\calP)
      &= d \cdot
         (-1)^{\frac{|\Delta|-1}{2}-1} \cdot 2^{-\frac{|\Delta|-1}{2}} \cdot
           4^{\frac{|\Delta|-1}{2}-1} \cdot  N_s^C (\Delta',\calP) \\
      &= d \cdot
         (-1)^{\frac{|\Delta|-1}{2}-1} \cdot 2^{\frac{|\Delta|-1}{2}-2} \cdot
           N_s^C (\Delta',\calP).
    \end {align*}
    In particular, in this case $N_{(0,s)}^W (\Delta,\calP)$ does not depend on
    the exact position of the points $\calP$.
  \end {enumerate}
  We will see in example \ref{ex-broccoli-even} that in some cases these
  results hold for broccoli invariants as well.
\end{example}

\begin{remark}[Algebraic Welschinger invariants]
    \label {rem-toricdelpezzo}
  To see the enumerative meaning of the Welschinger numbers let us now discuss
  a Correspondence Theorem stating that our tropical count determines the
  algebraic Welschinger invariants, i.e.\ numbers of real rational curves
  passing through a set of conjugation invariant points, counted with weight
  $\pm 1$ according to the nodes. More precisely, let $\Sigma$ be a real toric
  unnodal Del Pezzo surface with the tautological real structure, and $D$ a
  real ample linear system on $\Sigma$. There are five such surfaces, namely
  $\PP^2$, $\PP^1\times \PP^1$ or $\PP^2$ blown up at $k\leq 3$ generic real
  points (denoted by $\PP^2_k$), equipped with the standard real structure. The
  linear system $D$ is in suitable toric coordinates generated by monomials
  $x^i y^j$, where $(i, j)$ ranges over all lattice points of a polygon $Q(D)$
  of the following form. If $ \Sigma=\PP^2 $ and $D$ is the class of $d$ times
  a line, then $Q(D)$ is the triangle with vertices $(0,0)$, $(d,0)$, and
  $(0,d)$. If $ \Sigma=\PP^1\times \PP^1 $ and $D$ is of bidegree $(d_1,d_2)$
  then $Q(D)$ is the rectangle with vertices $(0,0)$, $(d_1,0)$, $(d_1,d_2)$,
  and $(0,d_2)$. If $ \Sigma=\PP^2_k $ and $ D = d\cdot L- \sum_{i=1}^k d_i E_i
  $, (where $L$ denotes the class of the pull-back of a line, and $E_i$ denote
  the exceptional divisors of $ \PP^2_k \rightarrow \PP^2$), then $Q(D)$ is the
  trapezoid with vertices $(0,0)$, $(d-d_1,0)$, $(d-d_1,d_1)$, $(0,d)$ if $ k=1
  $, the pentagon with vertices $(d_2,0)$, $(d-d_1,0)$, $(d-d_1,d_1)$, $(0,d)$,
  $(0,d_2)$ if $ k=2 $, and the hexagon with vertices $(d_2,0)$, $(d-d_1,0)$,
  $(d-d_1,d_1), (d_3,d-d_3)$, $(0,d-d_3)$, $(0,d_2)$ if $ k=3 $.

  \begin {center} \input {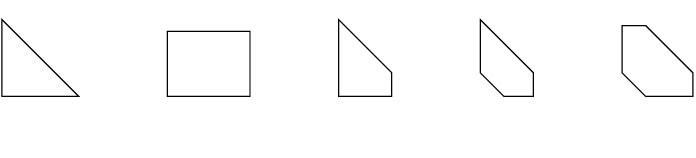} \end {center}

  Let $r$ and $s$ be non-negative integers satisfying $ \# (\partial Q(D) \cap
  \ZZ^2)-1 = r+2s $, and let $\mathcal{P}$ be a generic conjugation invariant
  set of $r+2s$ points of which exactly $r$ points are conjugation invariant
  themselves. By the Welschinger theorem (\cite{Wel03}, \cite{Wel05}), the set
  $\mathcal{R}(\Sigma,D,\mathcal{P})$ of algebraic real rational curves $
  C\in D $ passing through $\mathcal{P}$ is finite, consists only of nodal and
  irreducible curves, and the number
    \[ W_\Sigma(D,r,s)
         := \sum_{C\in \mathcal{R}(\Sigma,D,\mathcal{P})}(-1)^{s(C)} \]
  called \df {Welschinger invariant} does not depend on the special choice of
  $\mathcal{P}$, where $s(C)$ denotes the number of solitary nodes of $C$,
  i.e.\ real points where the curve is locally given by the equation
  $x^2+y^2=0$.
\end{remark}

\begin{definition}[Toric Del Pezzo degrees] \label{def-degree}
  We say that a degree $\Delta$ is \emph{toric Del Pezzo} if it consists of
  the primitive normal directions of facets of one of the polytopes $Q(D)$ of
  remark \ref{rem-toricdelpezzo}, where each direction appears $l$ times if $l$
  is the lattice length of the corresponding facet. If $Q(D)$ is the triangle
  with endpoints $(0,0)$, $(d,0)$ and $(0,d)$ (corresponding to the class of
  $d$ times a line in $\PP^2$), then we call curves of degree $\Delta$
  consisting of the normal directions $(-1,0)$, $(0,-1)$ and $(1,1)$ each $d$
  times \emph{curves of degree $d$}.
\end{definition}

Notice that a toric Del Pezzo degree consists of directions of weight one, so
the requirements of lemma \ref{lem-multequal} are satisfied.

\begin{theorem}[Correspondence Theorem] \label {thm-correspondence}
  Let $\Sigma$ be a toric Del Pezzo surface, $D$ a real ample linear system,
  $Q(D)$ the corresponding polytope as in remark \ref{rem-toricdelpezzo}, and
  $\Delta$ the corresponding degree. Let $r$ and $s$ satisfy $|\Delta|-1 =
  \# (\partial Q(D)\cap \ZZ^2)-1 = r+2s$. Then $ \Nrs^W (\Delta,\calP) =
  W_\Sigma(D,r,s) $ for any choice of points $\calP$ in general position.
  In particular, the Welschinger numbers $ \Nrs^W (\Delta,\calP) $ are
  independent of $ \calP $ in this case.
\end{theorem}

\begin {proof}
  Using remark \ref {rem-unparametrizednumbers}, this is theorem 3.1 of \cite
  {Shu06}. Note that the proof establishes not only an equality of numbers,
  but also a finite-to-one map between algebraic and tropical curves reflecting
  the tropical multiplicity.
\end {proof}

\begin{remark}[Welschinger numbers are not locally invariant in the moduli
    space] \label{rem-notinvariant}
  It is a striking feature of the Welschinger numbers $\Nrs^W (\Delta,\calP)$
  that, although they are invariant under $ \calP $ in the cases mentioned
  in theorem \ref {thm-correspondence}, one cannot show this by a local study
  of the moduli space as in the proof of theorem \ref {thm-broccoliinvariant}.
  In short, the reason for this is that the absence of the vertex type (6a)
  breaks the local invariance argument in the codimension-1 case (C$1_1$) (see
  the proof of theorem \ref {thm-broccoliinvariant}, in particular the table of
  codimension-1 cases and their resolutions).

  For example, consider a combinatorial type corresponding to a cell of $
  \Mrs^W(\Delta) $ of codimension one which locally contains the left picture
  $C$ below:

  \begin {center} \input {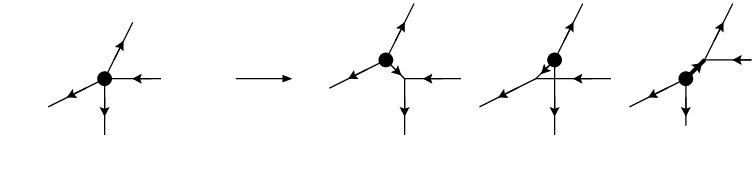} \end {center}

  Curves of this type pass through conditions which are not in general
  position, since the horizontal edge is fixed and the complex point is exactly
  on this horizontal line. There are two Welschinger curves $C_1$ and $C_2$
  as in the picture above such that this type appears in their boundary. Their
  multiplicities are $m_{C_1}=i^0\cdot 3i^2=-3$ and $m_{C_2}=i^0\cdot 1\cdot
  i^0=1$. We can see that they both satisfy the conditions when we move the
  complex point above the horizontal line. In contrast, no Welschinger curve
  satisfies the conditions if we move the point below the line: the third
  resolution $ C_3 $ would require a vertex of type (6a), which is not allowed
  for Welschinger curves. Thus locally around this codimension-1 cone, the
  number of Welschinger curves is not invariant.

  Of course, this leads to choices of $\Delta$ for which the Welschinger
  numbers are not invariant. For example, we can pick $\Delta=((1,0), (0,-1),
  (-2,-1), (1,2))$ such that the picture above is actually a global picture.
  Then this example shows that $\Nrs^W (\Delta,\calP)=-2$ if we pick $\calP$
  with the complex point above the horizontal line, and $\Nrs^W
  (\Delta,\calP)=0$ if we pick $\calP$ with the complex point below the line.
  Thus the numbers depend on the choice of $\calP$ and are not invariant.

  However, if $\Delta$ is a toric Del Pezzo degree as in definition
  \ref{def-degree}, then it follows from the Correspondence Theorem \ref
  {thm-correspondence} (and the Welschinger theorem) that the numbers $\Nrs^W
  (\Delta,\calP)$ are invariant. 

  Since this is true in spite of the missing local invariance around
  codimension-1 cones we can observe the following interesting fact about the
  moduli space $ \Mrs^W(\Delta) $ and the map $\ev$: given a collection of
  points $\calP$ not in general position such that a curve of a codimension-1
  type is in the preimage $\ev^{-1}(\calP)$ for which we do not have local
  invariance (as for the example above), there must be another curve in
  $\ev^{-1}(\calP)$ which is also of a codimension-1 type not satisfying
  local invariance, such that the differences to the invariance cancel exactly.
  For example, if we consider the above example as a local picture of the curve
  of degree $3$ below, then there is a second curve of codimension 1 such that
  the two differences cancel. The following picture shows these two
  codimension-1 curves passing through $\calP$ not in general position:

  \begin {center} \begin{picture}(0,0)%
\includegraphics{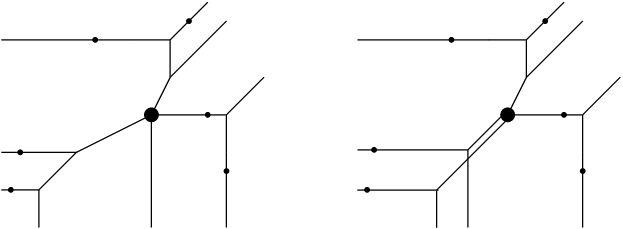}%
\end{picture}%
\setlength{\unitlength}{3947sp}%
\begingroup\makeatletter\ifx\SetFigFont\undefined%
\gdef\SetFigFont#1#2#3#4#5{%
  \reset@font\fontsize{#1}{#2pt}%
  \fontfamily{#3}\fontseries{#4}\fontshape{#5}%
  \selectfont}%
\fi\endgroup%
\begin{picture}(4974,1826)(2539,-7125)
\end{picture}%
 \end {center}

  We have seen already that the left picture produces a local difference of
  $-2$: locally, the difference between the numbers of curves passing through
  the configuration where we move the complex point up and down is $-2$. The
  right picture now produces a local difference of $+2$:

  \begin {center} \begin{picture}(0,0)%
\includegraphics{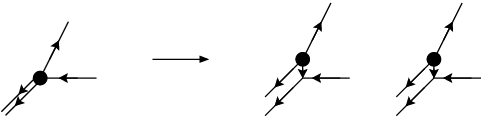}%
\end{picture}%
\setlength{\unitlength}{3947sp}%
\begingroup\makeatletter\ifx\SetFigFont\undefined%
\gdef\SetFigFont#1#2#3#4#5{%
  \reset@font\fontsize{#1}{#2pt}%
  \fontfamily{#3}\fontseries{#4}\fontshape{#5}%
  \selectfont}%
\fi\endgroup%
\begin{picture}(3859,924)(2679,-8173)
\end{picture}%
 \end {center}

  There are again two Welschinger curves which have this codimension one curve
  in their boundary (notice that the two edges pointing to the bottom-left are
  distinguishable in the big picture). They both satisfy the conditions when
  the complex point is moved up. No Welschinger curve satisfies the conditions
  if the complex point is moved down. Their multiplicity is $i^0\cdot 1\cdot
  i^0=1$ each.
\end{remark}

If the degree $\Delta$ is not a toric Del Pezzo degree, in particular if $
\Delta $ contains non-primitive vectors (i.e.\ we consider \df{relative}
Welschinger numbers), it may happen that these numbers are not even globally
invariant. This has already been observed in \cite{ABM08} with the following
example.

\begin{example}[Welschinger numbers are in general not invariant, see
    \cite {ABM08} section 7.2] \label {ex-notinvariant}
  The following picture shows the three Welschinger curves $ C_1 $, $ C_2 $, $
  C_3 $ (up to relabeling of the unmarked ends) of degree
    \[ ((-3,0),(0,-1),(0,-1),(0,-1),(1,1),(1,1),(1,1)) \]
  passing through some given configuration $ \calP $ of points. Each counts
  with multiplicity $3$, so for this configuration we have $ \Nrs^W (\Delta,
  \calP) = 9 $. For the configuration on the bottom right however, there is
  only one Welschinger curve $ C'$ passing through it, and it is of
  multiplicity one. So in this case $ \Nrs^W (\Delta,\calP')=1 $, i.e.\ the
  number depends on the choice of $\calP$.

  \begin {center} \input {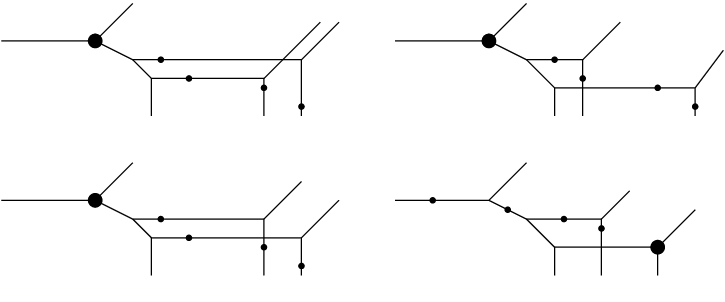} \end {center}
\end{example}

  \section {Bridge curves} \label {sec-bridge}

The aim of the following section is to prove that for toric Del Pezzo degrees
$\Delta$ (see definition \ref{def-degree}) the Welschinger numbers $\Nrs^W
(\Delta,\calP)$ coincide with the broccoli invariants $\Nrs^B (\Delta,\calP)$
(see corollary \ref{cor-bridge}). Since broccoli invariants are independent of
the chosen conditions, this result provides a tropical proof of the invariance
of Welschinger numbers, without having to use the detour via the Correspondence
and the Welschinger theorem. When considering degrees $\Delta$ that are not
toric Del Pezzo, the equivalence of Welschinger numbers and broccoli invariants
no longer holds, and consequently the Welschinger numbers may actually not be
invariant. 

We start with the definition of the class of bridge curves. It is a special
case of the class of oriented marked curves and includes oriented broccoli and
Welschinger curves. When a bridge curve is a broccoli curve having vertices of
type (6a) or a Welschinger curve having vertices of type (8), this curve allows
to start a so called \df{bridge}, that is, a 1-dimensional family of bridge
curves connecting broccoli and Welschinger curves. We show the invariance of
the curve multiplicities $m_C$ along these bridges, which then leads to the
equality of broccoli and Welschinger numbers mentioned above.

Throughout this section let $r,s \geq 0$, let $\Delta=(v_1,\ldots,v_n)$ be a
collection of vectors in $\ZZ^2 \setminus \{ 0 \}$, and let $F \subset \{
1,\ldots, n \}$ such that $|\Delta| -1 =r+2s+ |F|$. Moreover, fix conditions $
\calP \in \RR^{2(r+s)+|F|} $ in general position for $ \ev_F: \Mrs^\ori
(\Delta,F) \to \RR^{2(r+s)+|F|} $ as in definition \ref{def-general} and
example \ref{ex-general}, and consider only curves satisfying these conditions.

\begin{remark} \label{rem-vertexnumbers}
  Note that by lemma \ref{lem-dim} an oriented curve $C \in \Mrs^\ori
  (\Delta,F)$ all of whose vertices are of the types (1) to (9) of definition
  \ref{def-vertex} satisfies $ \n7 = \n8+\n9 $ (similarly to remark
  \ref{rem-outends} \ref{rem-outends-a} for Welschinger curves).
\end{remark}

\begin{definition}[Bridge curves] \label{def-bridgecurve}
  Let $r$, $s$, $ \Delta $, and $F$ be as in remark \ref{rem-vertexnumbers}. A
  \df {bridge curve} consists of the data of:
  \begin{itemize}
  \item an oriented curve $ C \in \Mrs^\ori (\Delta,F) $ all of whose vertices
    are of the types (1) to (9) of definition \ref{def-vertex}, and
  \item a bijection between its vertices of type (7) and those of types (8) or
    (9) (see remark \ref{rem-vertexnumbers}),
  \end{itemize}
  such that the following conditions hold:
  \begin{enumerate}
  \item \label{def-bridgecurve-a}
    There is at most one vertex of type (9).
  \item \label{def-bridgecurve-b}
    Each vertex of type (8) or (9) is connected to its corresponding vertex
    of type (7) (under the given bijection) starting with one of its even edges
    by a sequence of edges with no markings on them.
  \item \label{def-bridgecurve-c}
    Consider the set $M$ of vertices of type (6a) and (7); by abuse of notation
    we will sometimes also think of it as the set of all complex markings at
    these vertices. We split this set as $ M= \M8 \stackrel{\cdot}{\cup} \M9
    \stackrel{\cdot}{\cup} \M{6a} $, where
    \begin{itemize} 
    \item $ \M8 $ contains the vertices of type (7) corresponding to vertices
      of type (8) under the given bijection,
    \item $ \M9 $ contains the vertices of type (7) corresponding to vertices
      of type (9) under the given bijection,
    \item $ \M{6a} $ contains the vertices of type (6a).  
    \end{itemize}
    We define a partial order on $M$ by considering each vertex in $M$ with one
    even adjacent edge --- in the case of a vertex of type (7) we take the edge
    that does not connect this vertex to its corresponding vertex of type (8)
    or (9). For complex markings $ x_i \neq x_j $ in $M$ we say $ x_i < x_j $
    if the unique path connecting $x_i$ and $x_j$ does not pass through the
    even edge of $x_i$, but does pass through the even edge of $x_j$. Refine
    this partial order to a total order by considering vertices which are
    minimal under the partial order and comparing the (numerical) value of
    their markings. Choose the numerically minimal one and repeat the procedure
    without the chosen vertex until all vertices are ordered. We require now
    that the labeling of the complex markings is chosen such that vertices in $
    \M8 $ are smaller than vertices in $ \M9 $, and vertices in $ \M9 $ are
    smaller than vertices in $ \M{6a} $. 
  \end{enumerate}
  The multiplicity $ m_C $ of a bridge curve $C$ is given as usual by
  definition \ref{def-vertex}.
\end{definition}

\begin{example} \label{ex-order}
  For an example of the partial order in definition \ref{def-bridgecurve}
  \ref{def-bridgecurve-c} consider the picture below on the left, in which $
  x_2 $, $ x_3 $, and $ x_5 $ are the complex markings of type (6a) or (7). We
  have $ x_5 < x_2 < x_3 $, where dotted lines stand for parts of the graph
  between the distinguished edges and vertices. In this case, the total order
  on $M$ of definition \ref{def-bridgecurve} \ref{def-bridgecurve-c} agrees
  with this partial order. In the picture on the right however we get the
  partially ordered sets $ x_7 < x_8 < x_5 < x_1 $, $ x_7 < x_8 < x_2 < x_3 $,
  $ x_6 < x_4 $, and the total order $ x_6 < x_4 < x_7 < x_8 < x_2 < x_3 < x_5
  < x_1 $.
  \begin {center} \input {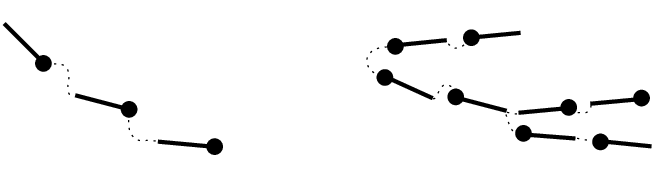} \end {center}
\end{example}

\begin{example} \label{ex-bridgecurve}
  An example of a bridge curve (containing a vertex of type (9)) is given in
  the following picture; the bijection between the vertices of type (7) and
  those of types (8) and (9) is indicated by the dotted arrows. We have labeled
  the vertices by their types only in the cases (6), (7), (8), and (9)
  since these are the most relevant ones for our study of bridge curves. In
  this example we have $ M=\{x_3,x_5,x_6\} $ and $ \M8 = \{x_5\} $, $ \M9 =
  \{x_6\} $, $ \M{6a} = \{x_3\} $. The partial order on $M$ is given by $ x_6 <
  x_3 $ and the total order by $ x_5 < x_6 < x_3 $. The dashed edges are
  ordinary odd edges (they form a string as explained in definition
  \ref{def-string} and remark \ref{rem-string}).

  \begin {center} \input {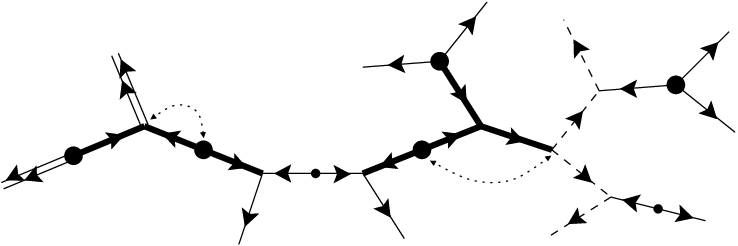} \end {center}
\end{example}

\begin{remark} \label{rem-fatpath}
  From the allowed vertex types of definition \ref{def-vertex} it follows that
  the sequence of edges of definition \ref{def-bridgecurve}
  \ref{def-bridgecurve-b} connecting each vertex of type (7) to its
  corresponding vertex of type (8) or (9) just contains even edges which are
  then adjacent to vertices of type (4).
\end{remark}

\begin{remark} \label{rem-totalorder}
  The choice of the total order refining the partial order in definition
  \ref{def-bridgecurve} \ref{def-bridgecurve-c} is not important. While the
  definition of bridge curves depends on this choice, the result of invariance
  in theorem \ref{thm-bridgeinv} does not.
\end{remark}

\begin{remark}[Dimension of the space of bridge curves] \label{rem-genpos}
  These (oriented) bridge curves can be constructed with the bridge algorithm
  \ref{algo-bridge} from oriented broccoli or Welschinger curves without
  changing the conditions $ \calP $. In particular, bridge curves are curves
  passing through conditions in general position. In fact, since the number of
  our conditions is $ 2(r+s)+|F| $ it follows from lemma \ref{lem-dim}
  that the space of bridge curves of a given combinatorial type through $ \calP
  $ is $0$-dimensional if there is no vertex of type (9) (i.e.\
  if $ \M9 = \emptyset $), and $1$-dimensional otherwise. If we even have $
  \M8 = \M9 = \emptyset $ or $ \M9 = \M{6a} = \emptyset $, the bridge curves
  specialize to the broccoli and Welschinger curves that we already know:
\end{remark}

\begin{lemma}[Broccoli and Welschinger curves as bridge curves]
    \label{lem-bridgecurveex}
  For fixed $r$, $s$, $\Delta$, $F$ the operation of forgetting the
  correspondence between the vertices of type (7) and those of types (8) or (9)
  of definition \ref{def-bridgecurve} induces bijections between curves through
  $ \calP $
  \begin{align*}
    \{\text{bridge curves with $ \M8=\M9=\emptyset $}\}
      & \quad \stackrel{1:1}{\longleftrightarrow} \quad
    \{\text{oriented broccoli curves}\} \\
    \text{and} \qquad \quad
    \{\text{bridge curves with $ \M9=\M{6a}=\emptyset $}\}
      & \quad \stackrel{1:1}{\longleftrightarrow} \quad
    \{\text{oriented Welschinger curves}\}.
  \end{align*}
\end{lemma}

\begin{proof}
  First of all, given a bridge curve with $ \M8=\M9=\emptyset$, it follows
  directly $ \n7 = \n8 = \n9 = 0 $. Hence the curve consists only of vertices
  of types (1) to (6) and is therefore a broccoli curve. In the same way, $ \M9
  = \M{6a} = \emptyset $ for a bridge curve implies $ \n9=0 $ and $ \n{6a}=0 $
  by definition \ref{def-bridgecurve} \ref{def-bridgecurve-c}. So we obtain a
  Welschinger curve. Hence the two maps of the lemma (from left to right) are
  well-defined.

  Conversely, an oriented broccoli curve has only vertices of type (1) to (6).
  Hence $ \M8=\M9=\emptyset$, and the correspondence between vertices of types
  (7), (8), and (9) is trivial. So the statement of the lemma about broccoli
  curves is obvious.

  Analogously, we have $ \M9=\M{6a}=\emptyset$ for each oriented Welschinger
  curve as we just allow vertices of types (1) to (5), (6b), (7), and (8).
  Conditions \ref{def-bridgecurve-a} and \ref{def-bridgecurve-c} of definition
  \ref{def-bridgecurve} are clear. So we have to prove the existence and
  uniqueness of a correspondence between the vertices of type (7) and (8) that
  satisfies \ref{def-bridgecurve-b}. To do this, we perform an induction over
  the number $ \n7 $ of vertices of type (7) in the underlying graph $\Gamma$.
  For $ \n7=0 $ there is nothing to show. Let $V$ be such a vertex of type (7)
  in a connected component $\Gamma'$ of $\Gamma_\even $ such that the part of
  $\Gamma' \setminus \{V\}$ not containing the root of $\Gamma'$ (see
  definitions \ref{def-gammaeven} and \ref{def-welschinger}
  \ref{def-welschinger-b} and the equivalence of oriented and unoriented
  Welschinger curves through conditions in general position in proposition
  \ref{prop-welschingeroriented}) contains no other vertices of type (7). Using
  remark \ref{rem-outends} \ref{rem-outends-b} for the encircled part $R$ in
  the picture below, we know that it has exactly one vertex $W$ of type (8).
  Now $V$ and $W$ are obviously connected by a sequence of even edges as
  required by definition \ref{def-bridgecurve} \ref{def-bridgecurve-b}, and
  moreover $V$ is the only vertex of type (7) that $W$ can be connected to
  without passing through other markings. Cut off $R$ and replace $V$ by a
  vertex of type (6b). Applying the induction hypothesis to the rest of
  $\Gamma$, we obtain the required existence and uniqueness of the bijection
  between the vertices of type (7) and (8).

  \begin {center} \input {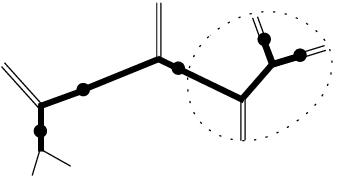} \end {center}

  \vspace{-4ex}
\end{proof}

We will now study the $1$-dimensional types of bridge curves through $ \calP $
and the boundary cases to which they can degenerate.

\begin{definition}[Strings] \label{def-string}
  Let $C=(\Gamma, x_1,\ldots,x_{r+s},y_1,\ldots,y_n,h) \in \Mrs^\ori (\Delta,F)
  $ be an oriented marked curve. As in definition 3.5 (a) of \cite{GM05b}, a
  \df{string} of $C$ is a subgraph of $\Gamma$ (after the end-gluing of
  definition \ref{def-endgluing}) homeomorphic to $\RR$ which does not
  intersect the closures $\overline{x_i}$ of the marked points and whose two
  ends are not fixed.
\end{definition}

\begin{remark} \label{rem-string}
  A bridge curve with a vertex of type (9) contains a unique string (containing
  this vertex) since the orientation of the two odd edges prescribes that they
  both lead in a unique way to a non-fixed unbounded end without passing
  through any markings (see example \ref{ex-bridgecurve}). As an example, the
  dashed edges in example \ref{ex-bridgecurve} are ordinary odd edges; they
  form a string.

  Note that the allowed vertex types require that these paths to the non-fixed
  unbounded ends go only through vertices of types (2) and (3). In particular,
  the string then contains only odd edges. On the other hand, a curve without
  vertex of type (9) does not contain a string.

  By remark \ref{rem-genpos}, a bridge curve through conditions in general
  position that has a vertex of type (9) (and thus a string) moves in a
  1-dimensional family --- namely by moving this string, as already observed
  in remark 3.6 of \cite{GM05b}. Let us now figure out what boundary cases can
  occur at the end of such 1-dimensional families.
\end{remark}

\begin{lemma}[Codimension-$1$ cases for bridge curves] \label{lem-dockcases}
  Let $C$ be a bridge curve through $ \calP $ with a vertex of type (9), thus
  having a string as in remark \ref{rem-string}. This string can be moved until
  two vertices of $C$ merge. The possible resulting vertices are as follows; we
  call them \df{codimension-$1$ cases for bridge curves}. As before, the arc in
  type (D2) means that the two odd edges must not be ends of the same
  direction.

  \begin {center} \input {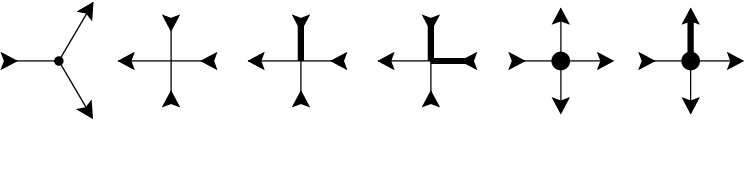} \end {center}

  \begin {center} \input {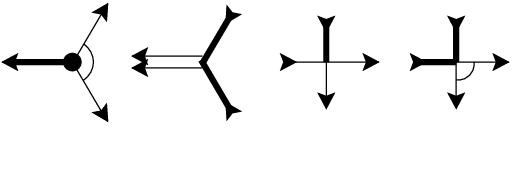} \end {center}

\end{lemma}

\begin{proof}
  For the terminology used in the following, we refer to the proof of
  theorem \ref{thm-broccoliinvariant}. Note that, when moving the string, two
  vertices on the string can merge, or a vertex on the string can merge with a
  vertex not on the string (if the two vertices are connected by a bounded
  edge).
  \begin{description} \itemsep0.5ex plus 0.5ex
  \item[Case 1] Assume the two vertices merging are of types (1) to (6). Then
    $V$ is a vertex of type (A$\,\cdot\,$), (B$\,\cdot\,$), or (C$\,\cdot\,$).
    The bridge curve we started with has already a vertex $W$ of type (9).
    Hence, just resolutions that do not create a vertex of type (9) are
    allowed. As $C$ originates from a bridge curve with a string, two of the
    edges adjacent to $V$ are contained in the string; more precisely by remark
    \ref{rem-string} there must be one incoming and one outgoing odd edge.
    If we just consider vertices with allowed bridge curve resolutions, the
    only possible vertices which are left then are (A1),
    ($\textnormal{B1}_\textnormal{1}$), (B3), (B5),
    ($\textnormal{C1}_\textnormal{1}$), ($\textnormal{C1}_\textnormal{3}$), and
    (C3).
  \item[Case 2] One vertex is of type (1) to (8) and the other one of type (7)
    or (8). Note that the string has to pass through one of the merging
    vertices in order to create the codimension-$1$ case. So we cannot have two
    vertices of type (7) and/or (8) as they do not allow the existence of the
    string. We thus need one vertex of type (1) to (6) which has one incoming
    and one outgoing odd edge, i.e.\ a vertex of type (3) merging with a vertex
    of type (7). But in this case, this vertex of type (7) (which necessarily
    lies in $ \M8 $) is bigger than the type (7) vertex in $ \M9 $
    corresponding to the type (9) vertex at which the string starts --- in
    contradiction to part \ref{def-bridgecurve-c} of the definition
    \ref{def-bridgecurve} of a bridge curve. And indeed, the vertex arising
    from merging type (3) with (7) has no other legal resolution, so such a
    case does not appear. Case 2 is thus impossible.
  \item[Case 3] One of the vertices is of type (9). Then the other vertex must
    be of type (2) to (4) or (7) as the other vertices of type (1), (5), (6),
    (8) do not fit together with the parity and the direction of the edges
    adjacent to the vertex of type (9).

    \begin{itemize}
    \item If $V$ arises from merging a vertex of type (9) with a vertex of type
      (7) we obtain a bridge curve with a vertex of type (6a), but without
      vertex of type (9).
    \item Merging a vertex of type (9) with a vertex of type (3) gives a bridge
      curve with a vertex of type (8) or (D2), depending on whether the
      resulting two odd edges are ends of the same direction or not.
    \item If the second vertex is of type (2) or (4), we obtain a vertex of
      type (D1) resp.\ (D2). \qedhere
    \end{itemize}
  \end{description}
\end{proof}

\begin{remark}[Bridge graphs and bridges] \label{rem-bridge}
  We are now able to explain the idea of bridges connecting broccoli to
  Welschinger curves more precisely. For this let us construct a so-called
  \df{bridge graph} as follows: the edges are the $1$-dimensional types of
  bridge curves through $ \calP $ (i.e.\ those containing a vertex of type (9)
  and thus a string), and the vertices are their $0$-dimensional boundary
  degenerations as described in lemma \ref{lem-dockcases} (we will see in lemma
  \ref{lem-end} that in the toric Del Pezzo case the string movement actually
  ends at both sides and thus leads to two vertices for each edge in the bridge
  graph). Note that the bijection between vertices of type (7) and those of
  types (8) and (9) that we have for the $1$-dimensional types can be extended
  to a map between vertices in the $0$-dimensional boundary types. We identify
  two such $0$-dimensional boundary types, i.e.\ represent them by the same
  vertex in the bridge graph, if they have the same underlying oriented curve
  and this map between vertices agrees, where we discard any mapping of a
  vertex to itself (which can occur if a type (7) vertex merges with a type (9)
  vertex to one of type (6a)).

  Note that some vertices in the bridge graph correspond to bridge curves with
  no type (9) vertex, whereas others (corresponding to codimension-$1$ cases
  (A$\,\cdot\,$), (B$\,\cdot\,$), (C$\,\cdot\,$), (D$\,\cdot\,$)) are not
  bridge curves in the sense of our definition. Included are however (as we
  will see in theorem \ref{thm-bridgeinv}) all broccoli and Welschinger curves
  through $ \calP $, so that we can think of the bridge graph as connecting
  broccoli and Welschinger curves. We will call a connected component of the
  bridge graph a \df{bridge}.
  
  The following picture shows a schematic example of a bridge graph. Its
  vertices corresponding to broccoli and Welschinger curves are drawn as big
  dots (on the left resp.\ right hand side of the diagram), the other ones as
  small dots. The dashed line indicates a curve which is both broccoli and
  Welschinger (i.e.\ has $ \M8=\M9=\M{6a}=\emptyset $), so it does not
  correspond to an edge in the bridge graph. The broccoli and Welschinger
  curves, as well as the $1$-dimensional types of bridge curves, are labeled
  with their multiplicities as in definition \ref{def-vertex}. 

  \begin {center} \input {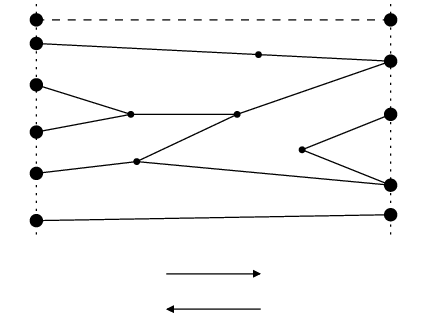} \end {center}

  The idea to prove the equality of broccoli and Welschinger numbers is now
  that there is a \df{local balancing condition} on the bridge graph, i.e.\
  that (as in the picture above) at each vertex the sum of the incoming equals
  the sum of the outgoing curve multiplicities when we move from the broccoli
  to the Welschinger side. To make this idea work, we first of all have to see
  that the edges of the bridge graph have a natural orientation so that it is
  well-defined which direction leads to the broccoli and which to the
  Welschinger side.
\end{remark}

\begin{definition}[Direction of string movement] \label{def-stringdir}
  For a given bridge curve $C$ with a vertex $V$ of type (9) consider the
  even edge $E$ adjacent to $V$. Changing the length of $E$ induces the
  movement of the string in $C$. Namely, making this edge longer makes the
  curve ``more Welschinger''; we want to call this the \df{positive direction}
  ($+$) of the string movement. Making $E$ shorter leads to a ``more broccoli''
  like curve; we want to call this the \df{negative direction} ($-$) of the
  string movement.
\end{definition}

\begin{theorem}[Invariance along bridges] \label{thm-bridgeinv}
  Let $C$ be an oriented curve containing a vertex $V$ of one of the
  codimension-$1$ types (A$\,\cdot\,$), (B$\,\cdot\,$), (C$\,\cdot\,$),
  (6a)/(8), or (D$\,\cdot\,$) as in lemma \ref{lem-dockcases}, and only
  vertices of types (1) to (9) otherwise. Assume as in lemma
  \ref{lem-dockcases} that $C$ arises from moving a string in a bridge curve
  with a vertex of type (9). Consider all bridge curves $C'$ that resolve $C$
  and that have matching bijections between their vertices of type (7) and
  those of type (8) and (9). (In the language of remark \ref{rem-bridge} this
  means that $C$ corresponds to a vertex and $C'$ to its adjacent edges in the
  bridge graph.)
  
  The curves $C'$ all contain a string and thus we can define $\sign_{C'}$ as
  the direction of the movement of the string away from $C$. Then $ \sum_{C'}
  \sign_{C'} \cdot m_{C'} $ equals\dots
  \begin{enumerate}
  \item \label{thm-bridgeinv-a}
    $ m_C $ if $C$ is a broccoli curve (i.e.\ we are on the left side of the
    bridge graph in remark \ref{rem-bridge});
  \item \label{thm-bridgeinv-b}
    $ -m_C $ if $C$ is a Welschinger curve (i.e.\ we are on the right side of
    the bridge graph);
  \item \label{thm-bridgeinv-c}
    $0$ in all other cases.
  \end{enumerate}
\end{theorem}

\begin{proof}
  For the terminology used in the following, we refer to the proof of theorem
  \ref{thm-broccoliinvariant}. We consider the resolving bridge curves $C'$
  and distinguish the types of $V$ as in lemma \ref{lem-dockcases}.

  \textbf {Case 1:} $V$ is a vertex of type (A$\,\cdot\,$), (B$\,\cdot\,$), or
  (C$\,\cdot\,$) (we are then in case \ref{thm-bridgeinv-c} of the theorem).
  Imagine to put a marking $m$ on the bounded edge adjacent to $V$ that
  connects this vertex on the string to the vertex $W$ of type (9). We then
  compare the resulting $H$-sign as in the proof of theorem
  \ref{thm-broccoliinvariant} with the direction of the string movement for
  $C'$. We know from \ref{lem-dockcases} that $V$ can be resolved into two
  vertices of types (1) to (6). As the two odd edges adjacent to $W$ are
  contained in the string, the $1$-dimensional movement of the marking $m$
  generated by resolving $V$ is reflected by the $1$-dimensional movement of
  the string and hence by varying the length of the even edge at $W$:

  \begin {center} \input {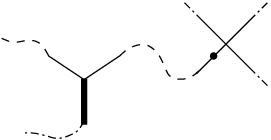} \end {center}

  Thus the $H$-sign equals the sign defined by the direction of the string
  movement (up to the same sign for all resolutions). Since we proved $
  \sum_{C'} \text {($H$-sign)} \cdot m_{C'}=0 $ in theorem
  \ref{thm-broccoliinvariant} already, it only remains to be shown in each case
  that all resolving curves are actually bridge curves, i.e.\ satisfy the
  conditions \ref{def-bridgecurve-a} to \ref{def-bridgecurve-c} of definition
  \ref{def-bridgecurve}. Condition \ref{def-bridgecurve-a} is always satisfied
  as we do not create a vertex of type (9).

  Concerning condition \ref{def-bridgecurve-b} of the definition of a bridge
  curve, note that in the cases (B$\,\cdot\,$) the connection between vertices
  of type (7), (8), and (9) are not modified as no vertices of type (7), (8),
  and (9) and no markings are involved. Hence, condition
  \ref{def-bridgecurve-b} is satisfied in all resolutions in this case. In the
  resolutions of vertices of type (A$\,\cdot\,$) and (C$\,\cdot\,$), no
  vertices of type (4) are involved, which are however necessary by remark
  \ref{rem-fatpath} to connect vertices of type (7) and (8), (9). Hence, also
  in these cases condition \ref{def-bridgecurve-b} is satisfied in all
  resolutions.

  Looking at condition \ref{def-bridgecurve-c} of definition
  \ref{def-bridgecurve}, the cases (A$\,\cdot\,$) and (B$\,\cdot\,$) are easy
  to manage as no vertices of type (6a) and (7) are involved (the partition of
  $M$ and the total order are not changed). For the case (C$\,\cdot\,$) we have
  to go into more details.
  \begin{enumerate}
  \item[($\textnormal{C1}_\textnormal{1}$)] Resolution (I) has a supplementary
    vertex $V$ of type (6). If the supplementary vertex is of type (6b), it is
    not contained in $M$ and need not be considered, so let us assume that
    $V$ is of type (6a). Then the set $M$ contains one more element (lying in
    $ \M{6a} $) compared to the resolutions (II) and (III). The string contains
    the edge $v_1$ and therefore, the vertex contained in $ \M9 $ also lies
    behind $v_1$. Hence, $V$ is bigger than the vertex of $ \M9 $ under the
    partial order. As the total order refines the partial order condition
    \ref{def-bridgecurve-c} is still satisfied. 
  \item[($\textnormal{C1}_\textnormal{3}$)] All three resolutions contain one
    more vertex of type (6a) in $ \M{6a} $ than $C$. But also in this case,
    this new vertex is bigger than the already existing vertex in $ \M9 $.
    Condition \ref{def-bridgecurve-c} is thus satisfied for all three
    resolutions simultaneously.
  \item[(C3)] Here, there are just two resolutions with a vertex of type (6a),
    where each time the new bounded edge is odd. The edge $v_2$ is even as
    before, the vertex in $ \M9 $ lies behind $v_1$, so the vertex in $ \M9 $
    and this vertex can be compared under the total order but not under the
    partial order. Hence, condition \ref{def-bridgecurve-c} satisfied in both
    cases simultaneously.
  \end{enumerate}
  In total, we can conclude that conditions \ref{def-bridgecurve-b} and
  \ref{def-bridgecurve-c} are fulfilled for all resolutions (if for any).

  \textbf {Case 2:} $V$ is a vertex of type (6a) or (8) (note that $V$ is a
  priori not unique then since $C$ has in general several vertices of type (6a)
  or (8)). We want to resolve vertices in this curve such that the resolutions
  are bridge curves with a vertex of type (9). The other way around we can ask
  ourselves which vertices in a bridge curve with vertex of type (9) can be
  merged in order to create $C$. After testing all possibilities we obtain two
  cases:
  \begin{itemize}
  \item[(A)] the vertex of type (9) can melt with a vertex of type (7) into a
    vertex of type (6a);
  \item[(B)] the vertex of type (9) can melt with a vertex of type (3) into a
    vertex of type (8), if the odd outgoing edge of the vertex of type (3) is
    an end and if one of the odd outgoing edges of the vertex of type (9) is
    also an end of the same direction.
  \end{itemize}
  Hence if we want to go the other way around, we can resolve
  \begin{itemize}
  \item[(A)] a vertex of type (6a) into a vertex of type (7) and a vertex of
    type (9);
  \item[(B)] a vertex of type (8) into a vertex of type (3) and a vertex of
    type (9). The so newly created bounded edge can have both orientations, due
    to the symmetric situation at the vertex of type (8). The question is just
    which of the vertices will become the vertex of type (3) and which one the
    vertex of type (9).   
  \end{itemize}
  For these two types of resolutions we have to check if the conditions
  \ref{def-bridgecurve-b} and \ref{def-bridgecurve-c} of the definition
  \ref{def-bridgecurve} of a bridge curve are satisfied.
  \begin{itemize}
  \item[(A)] The set $M$ remains the same as before resolving. The connections
    between vertices considered in condition \ref{def-bridgecurve-b} also
    remain the same. Before resolving the marking is at a vertex in $ \M{6a} $,
    but after resolving it becomes a vertex in $ \M9 $. This is just allowed if
    the marking was the smallest element in $ \M{6a} $, which is the case for
    exactly one marking if we assume $ \M{6a} \neq \emptyset $. Then the
    partial and the total order on $M$ also remain the same and condition
    \ref{def-bridgecurve-c} is satisfied.
  \item[(B)] The set $M$ is conserved also in this case. Consider the marking
    $x_i$ which corresponds to the vertex of type (8). In order to satisfy
    condition \ref{def-bridgecurve-b} of the definition we have to meet the
    vertex of type (9) at its even edge if we start at the marking. This means
    that we must choose the orientation of the inserted bounded edge such that
    this holds. To satisfy condition \ref{def-bridgecurve-c} the marking $x_i$
    has to be the biggest point in $ \M8 $ (assuming $ \M8 \neq \emptyset $).
    We need this since, after resolving the vertex, the marking lies in $ \M9 $
    and not anymore in $ \M8 $. But note that we still have two resolutions as
    we have two possibilities to enumerate the two odd edges at the vertex of
    type (8) that we resolve.       
  \end{itemize}
  Observe that both the multiplicity of the curve in (A) and the sum of the
  multiplicities of the two resolutions from (B) equal the multiplicity of $C$
  --- due to the fact that the multiplicity of the vertex of type (8) resolved
  in (B) is the double of the multiplicity of the vertex of type (3) after the
  resolution. Thus, as the even edge $E$ adjacent to the type (9) vertex
  becomes longer in (A) and shorter in the resolutions (B), the invariance
  holds if $ \M8 \neq \emptyset \neq \M{6a} $ so that both cases (A) and (B)
  exist. If $ \M8 $ is empty, the bridge curve we are looking at is a broccoli
  curve by lemma \ref{lem-bridgecurveex}. We then resolve a vertex of type (6a)
  by making $E$ longer. Hence $\sign_{C'} \cdot m_{C'}$ is plus the broccoli
  multiplicity. In the same way, if $ \M{6a} $ is empty, the considered bridge
  curve is a Welschinger curve by lemma \ref{lem-bridgecurveex}. As we then
  resolve a vertex of type (8), $E$ becomes shorter, so $\sign_{C'} \cdot
  m_{C'}$ is minus the Welschinger multiplicity.

  \textbf {Case 3:} $V$ is a vertex of type (D1) or (D2) (we are then in case
  \ref{thm-bridgeinv-c} of the theorem). Remember from lemma
  \ref{lem-dockcases} that $V$ can then be resolved into a vertex of type (2)
  to (4) and a vertex of type (9). The vertex of type (7) corresponding to the
  vertex of type (9) has to lie behind one of the even edges at the $4$-valent
  vertex by definition \ref{def-bridgecurve} \ref{def-bridgecurve-b}; we choose
  it to be behind the edge with direction $v_2$. The orientation and the parity
  of the bounded edge which appears when resolving are determined.  

  \begin {center} \input {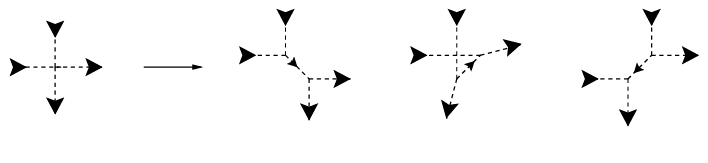} \end {center}

  Observe that resolution I does not exist for the vertex of type (D1) as the
  $3$-valent vertices that appear then are not allowed for bridge curves.
  The vertices appearing are listed in the table below. The last column
  $m_{\textnormal{I/II/III}}$ shows the absolute value of the product of the
  two vertex multiplicities in the resolutions I, II, and III.

  \begin {center} \tabcolsep 1mm \begin {tabular}{|c|ccc|ccc|ccc|} \hline
    codim-1 & \multicolumn 3 {c|}{resolution I}
            & \multicolumn 3 {c|}{resolution II}
            & \multicolumn 3 {c|}{resolution III} \\
    case    & $V$ & $W$ & $m_\textnormal{I}$
            & $V$ & $W$ & $m_\textnormal{II}$
            & $V$ & $W$ & $m_\textnormal{III}$\\ 
    \hline
    D1 & & & & $(2)$ & $(9)$ & 1 & $(2)$ & $(9)$& 1\\
    D2 & $(4)$ & $(9)$ & $|(v_1,v_2)|$ & $(3)$ & $(9)$ & $|(v_1,v_3)|$ &
      $(3)$ & $(9)$ & $|(v_1,v_4)|$ \\
    \hline 
  \end {tabular} \end {center}

  We have to check if conditions \ref{def-bridgecurve-b} and
  \ref{def-bridgecurve-c} of definition \ref{def-bridgecurve} are satisfied.
  Connections between vertices of type (7) to vertices of type (8) are not
  modified as no vertices of type (7), (8) and markings are involved in the
  resolutions. Similarly, the connection between the vertex of type (9) and the
  corresponding vertex of type (7) is not modified as the vertex of type (7)
  lies behind the edge of direction $v_2$. Hence, condition
  \ref{def-bridgecurve-b} is satisfied in all resolutions or in none of them.
  As no markings are involved in the resolutions, the set $M$, the splitting of
  $M$, and the total order are also preserved. So condition
  \ref{def-bridgecurve-c} holds in all three resolutions or in none of them.

  In order to prove the local invariance we also have to compute the direction
  of the string movement as in definition \ref{def-stringdir}. In resolution I
  we create a vertex of type (9), so the edge $E$ of definition
  \ref{def-stringdir} becomes longer.

  As in the proof of theorem \ref{thm-broccoliinvariant} we can imagine to have
  for the other resolutions II and III two other markings $P_1,P_2\in \RR^2$ on
  the edges $v_1,v_2$ as these are fixed. Hence we have two bounded edges of
  lengths $l_1$ and $l_2$, in addition to the (by resolving) new inserted
  bounded edge of length $a$. The direction of the string movement as in
  definition \ref{def-stringdir} is positive if and only if $l_2$ becomes
  longer when $a$ becomes longer. We can describe the condition that the curve
  has to pass through the given point conditions by the following linear
  systems of equations in the variables $l_1,l_2,a$.

  \begin {center} \tabcolsep 1mm
    \begin {tabular}{|ccc|c|}
      \multicolumn 4c{II} \\ \hline
      $l_1$ & $l_2$ & $a$ & \\ \hline
      $-v_1$ & $v_2$ & $-v_1-v_3$ & $P_2-P_1$ \\\hline
    \end {tabular}
    \hspace {6mm}
    \begin {tabular}{|ccc|c|}
      \multicolumn 4c{III} \\ \hline
      $l_1$ & $l_2$ &  $a$ & \\ \hline
      $-v_1$ & $v_2$ & $-v_1-v_4$ & $P_2-P_1$ \\ \hline
    \end {tabular}
  \end{center} 

  Obviously, these systems both have a one-dimensional space of solutions. In
  case II the homogeneous solution vector $(l_1,l_2,a)$ has the following
  entries:
    \[ l_1=(v_2,-v_1-v_3),\qquad l_2=-(-v_1,-v_1-v_3), \qquad a=(-v_1,v_2), \]
  where as above $(v_i,v_j)$ is the determinant of the matrix consisting of the
  column vectors $v_i$, $v_j$. So in order to determine the direction of the
  string movement we have to multiply the signs of $l_2$ and $a$, that is
  $\sign(v_1,v_3) \sign(v_1,v_2) $. In case III we just have to
  substitute the vector $v_3$ by $v_4$ and obtain therefore as sign $
  \sign(v_1,v_4)\sign(v_1,v_2)$. So in total the sign for the directions of the
  string movements are given by the following table.

  \begin{center}
    \begin {tabular}{|c|ccc|} \hline
      & sign for I & sign for II & sign for III \\ \hline
      (D) & $1$ & $\sign((v_1,v_3)(v_1,v_2))$ & $\sign((v_1,v_4)(v_1,v_2))$
        \\ \hline
    \end {tabular}
  \end{center}

  We are now able to verify the local invariance. We will use the same
  identities to deal with vertex multiplicities and signs as in the proof
  of theorem \ref{thm-broccoliinvariant}. Mainly, we use the formulas
  $\sign(v_i,v_j)i^{|(v_i,v_j)|-1}=i^{(v_i,v_j)-1}$ if $|(v_i,v_j)|$ is odd and
  $i^{|(v_i,v_j)|-1}=i^{(v_i,v_j)-1}$ if $|(v_i,v_j)|$ is even.

  In case (D1), we then obtain for the product of the vertex multiplicities
  together with the direction of the string movement in the resolutions II and
  III:
  \begin{align*}
    \textnormal{(II)} &= \sign((v_1,v_3)(v_1,v_2))\cdot i^{|(v_1,v_3)|-1}\cdot
         i^{|(v_2,v_4)|-1}=\sign(v_1,v_2)\cdot i^{(v_1,v_3)+(v_4,v_2)-2}, \\
    \textnormal{(III)} &= \sign((v_1,v_4)(v_1,v_2))\cdot i^{|(v_1,v_4)|-1}\cdot
         i^{|(v_2,v_3)|-1}=\sign(v_1,v_2)\cdot i^{(v_1,v_4)+(v_2,v_3)-2}.
  \end{align*}
  We have $\sign(v_1,v_2)\neq 0$ since $v_1$ and $v_2$ cannot be parallel as
  our curves pass through conditions in general position. Dividing equation
  (III) by (II) yields $i^{2(v_3,v_1)}=(-1)^{(v_3,v_1)}=-1$ as $(v_3,v_1)$ is
  odd. Hence (II)$+$(III)$=0$.
  
  Similarly, for (D2) we obtain:
  \begin{align*}
    \textnormal{(I)} &=
      |(v_1,v_2)|\cdot i^{|(v_1,v_2)|-1}\cdot
      i^{|(v_3,v_4)|-1}=\sign(v_1,v_2)\cdot(v_1,v_2)\,
      i^{(v_1,v_2)+(v_3,v_4)-2}, \\
    \textnormal{(II)} &=
      \sign((v_1,v_3)(v_1,v_2))\cdot
      |(v_1,v_3)|\cdot i^{|(v_1,v_3)|-1}\cdot i^{|(v_2,v_4)|}
      =\sign(v_1,v_2)\cdot(v_1,v_3)\, i^{(v_1,v_3)+(v_4,v_2)-2}, \\
    \textnormal{(III)} &=
      \sign((v_1,v_4)(v_1,v_2))\cdot|(v_1,v_4)|\cdot
      i^{|(v_1,v_4)|-1}\cdot i^{|(v_2,v_3)|-1}=\sign(v_1,v_2)\cdot
      (v_1,v_4)\, i^{(v_1,v_4)+(v_2,v_3)-2}.
  \end{align*}
  Let us divide all three terms by $ \sign(v_1,v_2) \,i^{(v_1,v_2)+(v_3,v_4)-2}
  $. For (I) we then get $ (v_1,v_2) $. In term (II) we obtain $ i^{2(v_2,v_1)}
  \cdot (v_1,v_3)=(-1)^{(v_2,v_1)}\cdot(v_1,v_3)=(v_1,v_3) $ as $(v_2,v_1)$ is
  even. Finally, for (III) we get $ i^{2(v_1,v_4)} \cdot (v_1,v_4) =
  (-1)^{(v_1,v_4)} \cdot (v_1,v_4) = (v_1,v_4) $ as $(v_1,v_4)$ is also even.
  So we have (I)$+$(II)$+$(III)$=(v_1,v_2)+(v_1,v_3)+(v_1,v_4)=0$.

  Hence we have shown the invariance for all codimension-$1$ cases for bridge
  curves.
\end{proof}

In order to prove the equality of broccoli and Welschinger numbers with the
idea of remark \ref{rem-bridge} we need one more final ingredient: that each
edge in the bridge graph is actually bounded, i.e.\ that the string movement in
each $1$-dimensional type of bridge curves is bounded in both directions by a
codimension-$1$ case. It is actually only this last step that requires a toric
Del Pezzo degree and thus spoils the equality of broccoli and Welschinger
numbers (as well as the invariance of Welschinger numbers, see example
\ref{ex-notinvariant}) in other cases.

\begin{lemma}[Boundedness of bridges] \label{lem-end}
  Assume that $ \Delta $ is a toric Del Pezzo degree (see definition
  \ref{def-degree}). Let $C$ be a bridge curve through $ \calP $ with a vertex
  of type (9), thus having a string as in remark \ref{rem-string}. Then the
  movement of the string within this combinatorial type is bounded in both
  directions.
\end{lemma}

\begin{proof}
  Assume that we have a bridge curve through $ \calP $ with a string that can
  be moved infinitely far. By the proof of proposition 5.1 in \cite{GM05b} such
  a string then has to consist of two edges which are both ends of the curve.
  Let us brief\/ly repeat the arguments for the sake of completeness. 

  \begin {center} \input {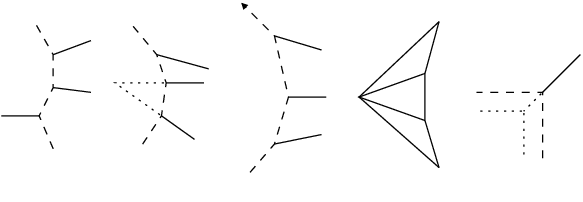} \end {center}

  If there are bounded edges adjacent to the string on both sides, the deformations of
  the string are bounded on both sides of the string, see case (a) in the
  figure above (where the string is marked with dashed lines). If
  there are only bounded edges adjacent to the string on one side as in (b), their
  extensions must not meet on the other side of the string, because otherwise
  the string is bounded on both sides as in (a). So the edges adjacent to the
  string look as in (c). This implies that the corresponding local dual
  subdivision has a concave side as depicted in (d). As the edges with
  direction vectors $v_1$ and $v_2$ in case (d) are dual to the two unbounded
  edges of the string, they must be $(\pm 1,0)$, $(0,\pm 1)$ or $\pm(1,1)$
  depending on the chosen toric Del Pezzo degree $\Delta$ (see definition
  \ref{def-degree}). Considering the lattice area of the triangle spanned by
  any two of these vectors, which is at most $1$, it is obvious that this
  triangle has no interior lattice point. Hence, there cannot be a vertex in
  the curve dual to this triangle. It follows that the string only consists of
  two unbounded edges as shown in (e).
  
  \vspace {1ex}

  \begin{sidepic}{vector}
    As we are dealing with bridge curves the string must then consist
    of the two odd edges adjacent to the vertex of type (9). From the
    definition of the vertex type (9) we know that the two ends cannot have the
    same direction. We thus see that these ends have two of the directions shown in the
    picture on the right. But in all these cases the third direction at the
    vertex of type (9) would be odd (in contradiction to the definition of
    type (9)) or $0$ (which is impossible for curves through conditions in
    general position). Hence the string movement cannot be unbounded.
  \end{sidepic}

  \vspace{-3.3ex}
\end{proof}

\begin{corollary}[Welschinger numbers $=$ broccoli invariants in the toric Del
    Pezzo case] \label{cor-bridge}
  Let $r,s \geq 0$, let $\Delta=(v_1,\ldots,v_n)$ be a toric Del Pezzo degree,
  and let $F \subset \{ 1,\ldots, n \}$ such that $|\Delta| -1 =r+2s+ |F|$. Fix
  a configuration $\calP$ of conditions in general position. Then $ \Nrs^W
  (\Delta,F,\calP) = \Nrs^B (\Delta,F,\calP) $.
\end{corollary}

\begin{proof}
  By theorem \ref{thm-bridgeinv} and definitions \ref{def-broccoli-inv} and
  \ref{def-welschinger-inv} we have
    \[ |G(\Delta,F)| \cdot
       \big( \Nrs^B (\Delta,F,\calP) - \Nrs^W (\Delta,F,\calP) \big)
       = \sum_C \, \sum_{C'} \sign_{C'} \cdot m_{C'}, \]
  where the sum is taken over all $C$ as in theorem \ref{thm-bridgeinv} and all
  resolutions $ C' $ of $C$ (i.e.\ over all vertices and adjacent edges in the
  bridge graph of remark \ref{rem-bridge}). Note that this in fact a finite sum
  since there are only finitely many types of bridge curves. Now by lemma
  \ref{lem-end} each $1$-dimensional type $C'$ of bridge curves occurs in this
  sum exactly twice with the same multiplicity, once with a positive and once
  with a negative sign. Hence the sum is $0$, proving the corollary.
\end{proof}

\begin{corollary}[Invariance of Welschinger numbers in the toric Del Pezzo
    case] \label{cor-welinv}
  With the assumptions and notations as in corollary \ref{cor-bridge}, the
  Welschinger numbers $ \Nrs^W (\Delta,F,\calP) $ are independent of the
  conditions $ \calP $.
\end{corollary}

\begin{proof}
  This follows from corollary \ref{cor-bridge} and theorem
  \ref{thm-broccoliinvariant}.
\end{proof}

In the remaining part of this section we want to construct bridges explicitly
and give some examples. The following algorithm, which follows from the proof
of theorem \ref{thm-bridgeinv}, shows how to construct a bridge from a given
starting point.

\begin{algo}[Bridge algorithm] \label{algo-bridge}
  Let $r,s\geq 0$, let $\Delta=(v_1,\ldots,v_n)$ be a toric Del Pezzo degree,
  and let $F\subset \{1,\ldots,n\}$ be such that $|\Delta|-1=r+2s+|F|$. Fix a
  configuration $\mathcal{P}$ of conditions in general position. Consider a
  bridge curve $C$ passing through $\mathcal{P}$; we want to construct the
  bridge that contains $C$.
  \begin{itemize}
  \item[(1)] If $C$ is a broccoli and Welschinger curve simultaneously (hence
    $ \M8=\M9=\M{6a}=\emptyset $), do nothing.
  \item[(2)] Given a bridge curve $C$ with $ \M9 \neq \emptyset$ (hence with a
    string) together with a direction for the movement of the string, move the
    string in the direction until we hit a codimension-$1$ type $C'$ as in
    lemma \ref{lem-dockcases}. Go to (2) with each new resolution in the
    direction away from $C'$.
  \item[(3)] If the curve is a broccoli curve, that is $ \M8=\M9=\emptyset$,
    choose the smallest vertex in $ \M{6a} $ under the total order defined in
    \ref{def-bridgecurve} \ref{def-bridgecurve-c}. Pull out an even edge of
    this vertex of type (6a) in order to create a vertex of type (7) and a
    vertex of type (9), thus producing a bridge curve with a string and a
    direction for the movement. Go to (2).
 \item[(4)] If the curve is a Welschinger curve, that is $ \M9=\M{6a}=\emptyset
    $, choose the vertex of type (8) corresponding to the biggest vertex in $
    \M8 $ under the total order defined in \ref{def-bridgecurve}
    \ref{def-bridgecurve-c}. Pull apart the two odd edges in order to create a
    string between the two even edges and a direction for the movement. We thus
    transform the vertex of type (8) into a vertex of type (3) and a vertex of
    type (9). Go to (2).
  \item[(5)] If the curve is a bridge curve with $ \M9 = \emptyset$, but $ \M8
    \neq \emptyset \neq \M{6a} $, we can choose the biggest vertex (under the
    total order) in $ \M8 $ or the smallest in $ \M{6a} $ in order to construct
    the bridge in direction ``broccoli'' or in direction ``Welschinger''.
    Transform the vertex as described in the two last items, respectively, thus
    producing a bridge curve with a string and a direction. Go to (2).
  \end{itemize}
\end{algo}

\begin{example}[A bridge connecting only broccoli curves] \label{ex-brocbridge}
  Following algorithm \ref{algo-bridge}, the following picture shows a bridge
  connecting one broccoli curve (a) to another broccoli curve (e) (and to no
  Welschinger curve). In curve (c) we resolve a $4$-valent vertex of type (D1).
  The types (b) and (d) are $1$-dimensional, the other three $0$-dimensional.

  \begin {center} \input {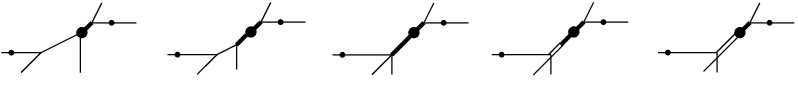} \end {center}
\end{example}

An example of a bridge connecting a broccoli curve with a Welschinger curve can
be found in section \ref{sec-content} of the introduction.

\begin{example}[Two cases that are not toric Del Pezzo] \label{ex-nopezzo}
  The boundedness of bridges of lemma \ref{lem-end}, and consequently the
  equality of broccoli and Welschinger numbers as well as the invariance of
  Welschinger numbers, are false in general for degrees that are not toric Del
  Pezzo:
  \begin{enumerate}
  \item \label{ex-nopezzo-a}
    Consider the following Newton polytope and its subdivision. It is obviously
    not toric Del Pezzo. A broccoli curve having this Newton subdivision is
    depicted on the right hand side. Starting the bridge as in algorithm
    \ref{algo-bridge} yields a string going to infinity (very right hand side),
    so the broccoli curve is not connected to a Welschinger curve by a bridge.

    \begin {center} \input {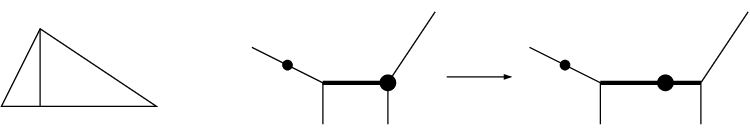} \end {center}
  \item \label{ex-nopezzo-b}
    Recall example \ref{ex-notinvariant} where we have shown that Welschinger
    numbers are not invariant if we do not have a toric Del Pezzo degree. If we
    choose the point configuration $\mathcal{P}$ as in example
    \ref{ex-notinvariant}, the Welschinger curves $C_1$, $C_2$, $C_3$ with
    multiplicity $3$ shown there are also broccoli curves, and in addition
    there are $4$ more broccoli curves passing through $\mathcal{P}$ as
    depicted below.

    \begin {center} \input {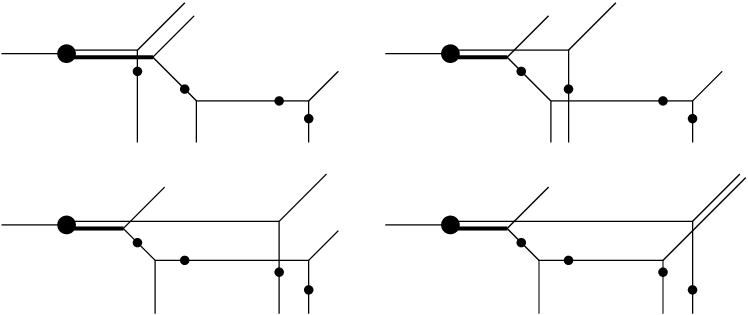} \end {center}

    Each of them has multiplicity $ -2 $, so the broccoli invariant is
    $ \Nrs^B (\Delta,\calP) = 3 \cdot 3 + 4 \cdot (-2) = 1 $. In particular,
    it is not equal to $ \Nrs^W (\Delta,\calP) = 9 $. Indeed, starting a bridge
    at the complex marking of each of the four curves above gives a curve
    having a string going to infinity as in \ref{ex-nopezzo-a}, so the
    contribution of $ -8 $ to the broccoli invariant is not seen on the
    Welschinger side.
  \end{enumerate}
\end{example}

\begin{example}[Broccoli invariants for degrees with non-fixed even ends]
    \label{ex-broccoli-even}
  By remark \ref{rem-string} the ends of a string are always unfixed and odd.
  In particular, this means that the proof of lemma \ref{lem-end} (and thus
  also of the equality of broccoli and Welschinger numbers) only requires that
  the \textsl{unfixed odd} ends in $ \Delta $ are those occurring in a toric
  Del Pezzo degree.

  Let us review example \ref{ex-welschinger-even} from this point of view.
  \begin{enumerate}
  \item \label{ex-broccoli-even-a}
    If $ \Delta $ has more than one non-fixed even end, and all other non-fixed
    ends are only those occurring in a toric Del Pezzo degree, then the
    result $ \Nrs^W (\Delta,F,\calP) = 0 $ of example \ref{ex-welschinger-even}
    \ref{ex-welschinger-even-a} implies that also $ \Nrs^B (\Delta,F) = 0 $.
  \item \label{ex-broccoli-even-b}
    If $ \Delta $ has one non-fixed even end, and all other ends are non-fixed
    and among those occurring in a toric Del Pezzo degree, then the formula for
    $ \Nrs^W (\Delta,\calP) $ of example \ref{ex-welschinger-even}
    \ref{ex-welschinger-even-b} holds in the same way for $ \Nrs^B (\Delta) $.
  \end{enumerate}
\end{example}

  \section {The Caporaso-Harris formula for broccoli curves} \label {sec-ch}

In this section, we establish a Caporaso-Harris formula for broccoli curves of
degree dual to the triangle with endpoints $(0,0)$, $(d,0)$ and $(0,d)$. This
is a recursive formula computing all broccoli invariants with weight conditions
on fixed and non-fixed left ends in addition to the usual point conditions. As
usual for Caporaso-Harris type formulas, the idea to obtain these relations is
to move one of the point conditions to the far left so that the curve splits
into a left part (passing through the moved point) and a right part (passing
through the remaining points). Since broccoli invariants of curves with ends of
weight one (i.e.\ of degree $d$) equal Welschinger numbers $\Nrs^W (d)$ by
corollary \ref{cor-bridge} and the latter equal Welschinger invariants
$W_{\PP^2}(d,r,s)$ by the Correspondence Theorem \ref{thm-correspondence}, our
formula then computes all Welschinger invariants of the plane recursively.

It is also possible to use Welschinger curves directly to establish a similar
formula. However, since the numbers of Welschinger curves of degree dual to the
triangle with endpoints $(0,0)$, $(d,0)$, and $(0,d)$ and with ends of higher
weight are not invariant (as we have seen in example \ref {ex-notinvariant}),
the arguments are then getting significantly more complicated as one always has
to pick special configurations of points. This is the content of \cite{ABM08}.
There, the authors pick a configuration of points such that the Welschinger
curves passing through these points decompose totally into floors (see
proposition \ref {prop-shape}), and count them by means of floor diagrams. This
yields a recursive formula for floor diagrams which also computes all
Welschinger invariants of the plane.

Let us first fix some notation.

\begin {notation}
  Let $ \alpha = (\alpha_1,\ldots,\alpha_m) $, $ \beta = (\beta_1,\ldots,
  \beta_{m'}) $, $ \alpha^1 = (\alpha^1_1,\ldots,\alpha^1_{m_1}) $, \ldots, $
  \alpha^k = (\alpha^k_1,\ldots,\alpha^k_{m_k}) $ be finite sequences with $
  \alpha_i,\beta_i,\alpha_i^j \in \N $. For simplicity, we will usually
  consider them to be infinite sequences by setting the remaining entries to
  $0$. We then define:
  \begin {enumerate}
  \item $ | \alpha |:= \sum_{i=1}^m \alpha_i$,
  \item $ I \alpha := \sum_{i=1}^m i \cdot \alpha_i$,
  \item $ \alpha + \beta := (\alpha_1 + \beta_1,\alpha_2 + \beta_2, \dots) $,
  \item $ \alpha \leq \beta :\Leftrightarrow \alpha_i \leq \beta_i \text { for
    all $i$} $,
  \item $ \alpha < \beta :\Leftrightarrow \alpha \le \beta \text { and } \alpha
    \neq \beta $,
  \item $\binom{n}{\alpha_1,\ldots,\alpha_m}:=\frac{n !}{\alpha_1! \cdot \ldots
    \cdot \alpha_m! (n - \alpha_1-\ldots -\alpha_m)!}$ for $ |\alpha| \le n $,
  \item $\binom{\alpha}{\alpha^1,\ldots, \alpha^k} :=\prod_i
    \binom{\alpha_i}{\alpha_i^1,\ldots,\alpha_i^k}$.
  \end {enumerate}
  Furthermore, we define $e_k$ to be the sequence having only $0$ as entries
  except a $1$ in the $k$-th entry.
\end{notation}

\begin{definition}[Broccoli curves of type $ (\alpha, \beta) $]
    \label{def-delta}
  Let $ d>0 $, and let $\alpha$ and $\beta$ be two sequences satisfying $
  I\alpha+I\beta=d$. We define $ \Delta(\alpha,\beta) $ to be the degree
  consisting of $d$ times the vectors $ (0,-1) $ and $ (1,1) $ each, and $
  \alpha_i+\beta_i $ times $ (-i,0) $ for all $i$ (in any fixed order). Let $
  F(\alpha,\beta) \subset \{1,\dots,|\Delta(\alpha,\beta)|\} $ be a fixed
  subset with $ |\alpha| $ elements such that the entries of $
  \Delta(\alpha,\beta) $ with index in $F$ are $ \alpha_i $ times $ (-i,0) $
  for all $i$. If no confusion can result we will often abbreviate $
  \Delta(\alpha,\beta) $ as $ \Delta $ and $ F(\alpha,\beta) $ as $F$.

  Broccoli curves in $ \Mrs^B(\Delta,F) $ will be called \df {curves of type $
  (\alpha,\beta) $}. We speak of their unmarked ends with directions $ (-i,0) $
  as the \df {left ends}. So $ \alpha_i $ and $ \beta_i $ are the numbers of
  fixed and non-fixed left ends of weight $i$, respectively.
\end{definition}

\begin{definition}[Relative broccoli invariants]
  Let $ \Delta=\Delta(\alpha,\beta) $ and $ F=F(\alpha,\beta) $ be as in
  definition \ref{def-delta}, and $ r,s $ such that the dimension condition
  $ |\Delta|-1-|F| = 2d+|\beta|-1=r+2s $ is satisfied. To simplify notation, we
  define the \df {relative broccoli invariant}
    $$ N^d(\alpha,\beta,s) := \Nrs^B (\Delta(\alpha,\beta),F(\alpha,\beta)). $$
\end{definition}

\begin{remark}[Unlabeled non-fixed ends] \label{rem-unlabeled}
  Notice that by remark \ref{rem-endlabeling} a broccoli curve without labels
  on the unmarked ends yields $2^{-k}\cdot |G(\Delta,F)|$ labeled curves
  contributing to the broccoli invariant, where $|G(\Delta,F)|$ as in
  definition \ref {def-ev} \ref {def-ev-b} denotes the number of ways to
  relabel the non-fixed unmarked ends without changing the degree, and $
  k=\n{6b}+\n8 $ is the number of double ends. In contrast, in the definition
  \ref {def-broccoli-inv} of broccoli invariants we multiply the number of
  broccoli curves with $\frac{1}{|G(\Delta,F)|}$. Thus a curve without labels
  contributes $2^{-k}$ to the count. Hence, when counting broccoli curves whose
  non-fixed unmarked ends are not labeled, we have to change the multiplicity
  of vertices of type (6b) to $ \frac 12 \cdot i^{-1} $. In the following, we
  will drop the labels of the non-fixed ends and change the multiplicity
  accordingly. Note that for the degree $ \Delta $ and $F$ as above we have
  $|G(\Delta,F)|= d! \cdot d! \cdot \beta_1! \cdot \beta_2! \cdot \; \cdots $.
\end{remark}

\begin{remark}\label{rem-relbroccoliwelschinger}
  It follows from theorem \ref{thm-broccoliinvariant} that $N^d(\alpha,\beta,s)
  $ is invariant, i.e.\ does not depend on the choice of the conditions.
  Note that if $\alpha=(0)$ and $\beta=(d)$ then
    $$ N^d((0),(d),s)= \Nrs^B (d)=\Nrs^W (d)=W_{\PP^2}(d,3d-2s-1,s), $$
  where the second equality follows from theorem \ref{thm-bridgeinv} and the
  last equality from theorem \ref{thm-correspondence}.
\end{remark}

Now we describe the properties of configurations $\calP$ of points that we
obtain by moving one of the point conditions (w.l.o.g.\ $P_1$) to the left. Let
us show first that then curves satisfying these conditions decompose into a
left and a right part.

\begin{lemma}[Decomposing curves into a left and right part]
    \label{lem-pointleft}
  Let $\Delta$ and $F$ be as in definition \ref{def-delta}, and let $
  2d+|\beta|-1=r+2s $. Fix a small real number $\epsilon > 0$ and a large one
  $N>0$. Choose $r+s$ (real and complex) points $P_1,\ldots,P_{r+s}$ and
  $|\alpha|$ $y$-coordinates for the fixed left ends in general position such
  that
  \begin{itemize}
  \item the $y$-coordinates of all $P_i$ and the fixed ends are in the open
    interval $(-\epsilon,\epsilon)$,
  \item the $x$-coordinates of $P_2,\ldots,P_{r+s}$ are in $(-\epsilon,
    \epsilon)$,
  \item the $x$-coordinate of $P_1$ is smaller than $-N$.
  \end{itemize}
  Let $ C=(\Gamma,x_1,\dots,x_{r+s},y_1,\dots,y_n,h) \in \Mrs^B(\Delta,F)$ be a
  broccoli curve satisfying these conditions. Then no vertex of $C$ can have
  its $y$-coordinate below $-\epsilon$ or above $\epsilon$. There is a
  rectangle $R=[a,b]\times [-\epsilon,\epsilon]$ (with $a\geq -N$, $b\leq
  -\epsilon$ only depending on $d$) such that $R \cap h(\Gamma)$ contains only
  horizontal edges of $C$.
\end{lemma}

\begin{proof}
  Notice that it follows from lemma \ref{lem-general} that each connected
  component of $C$ minus the marked points contains exactly one non-fixed
  unmarked end, a statement analogous to remark 2.10 of \cite{GM05c}. The fact
  that the $y$-coordinates of the vertices of $C$ cannot be above $\epsilon$ or
  below $-\epsilon$ and the existence of the rectangle $R$ follow analogously
  to the first part of the proof of theorem 4.3 of \cite{GM05c}.
\end{proof}

A configuration of points and $y$-coordinates for the fixed left ends as in
lemma \ref {lem-pointleft} can be obtained from any other by moving $P_1$ far
to the left. So in this situation the curves decompose into a left and a right
part connected by only horizontal edges in the rectangle $R$. A picture showing
this can be found in example \ref{ex-leftright}. In the following, we study the
possibilities for the shapes of the left and right part.

\begin{notation}[Left and right parts] \label{not-decompose}
  With notations as in lemma \ref {lem-pointleft}, cut $C$ at each bounded
  edge $e$ such that $ h(e)\cap R\neq \emptyset $. Denote the component passing
  through $P_1$ by $C_0$ (the left part), and the union of the other connected
  components by $\tilde{C}$ (the right part).
\end{notation}

\begin{proposition}[Possible shapes of the left and right part]
    \label{prop-shape}
  Let $ C_0 $ and $ \tilde C $ be the left and right part of a broccoli curve
  as in lemma \ref {lem-pointleft} and notation \ref {not-decompose}.
  \begin {enumerate}
  \item \label {prop-shape-a}
    If $C_0$ has no bounded edges, it looks like (A), (B), or (C) in the
    picture below (in which the edges are labeled with their weights).
    Moreover:
    \begin {itemize}
    \item In case (A), $ \tilde C $ is an irreducible curve of
      type $ (\alpha+e_k,\beta-e_k)$.
    \item In case (B), $ \tilde C $ is an irreducible curve of type $
      (\alpha+e_{k_1+k_2},\beta-e_{k_1}-e_{k_2}) $.
    \item In case (C), $ \tilde C $ decomposes into two connected components
      $C_1$ and $C_2$ of types $(\alpha^1,\beta^1)$ resp.\ $(\alpha^2,\beta^2)$
      with $I\alpha^j+I\beta^j=d_j$ for $ j=1,2 $, $d_1+d_2=d$, $\alpha^1+
      \alpha^2=\alpha+e_{k_1}+e_{k_2}$, and $\beta^1+\beta^2=\beta-e_{k_1+k_2}
      $. The curve $C_j$ for $ j=1,2 $ passes through $r_j$ real and $s_j$
      complex given points, where $ 2d_j+|\beta^j|-1=r_j+2s_j $.
    \end {itemize}
    In case (A) (for real $ P_1 $) the left end is odd, in the cases (B) and
    (C) (for complex $ P_1 $) exactly one of the three edges adjacent to $ P_1
    $ is even.
  \item \label {prop-shape-b}
    If $C_0$ has bounded edges (it is then called a \df{floor}), it looks like
    (D), (E), or (F) in the picture below, and has one end of direction
    $(0,-1)$ and one of direction $(1,1)$. We call the ends of $ C_0 $ of
    direction $ (i,0) $ for $ i>0 $ the \df {right ends}. Moreover:
    \begin {itemize}
    \item In case (D) (for real $ P_1 $), $ C_0 $ has only fixed left and right
      ends.
    \item In case (E) (for complex $ P_1 $), $ P_1 $ is adjacent to a left
      non-fixed end of $ C_0 $, and all other left and right ends of $ C_0 $
      are fixed.
    \item In case (F) (for complex $ P_1 $), $ P_1 $ is adjacent to a right
      non-fixed end of $ C_0 $, and all other left and right ends of $ C_0 $
      are fixed.
    \end {itemize}
    In any case, $ \tilde C $ consists of some number $l$ of connected
    components $C_1,\ldots,C_l$. Each $C_j$ is a curve of some type
    $(\alpha^j,\beta^j)$ with $I\alpha^j+I\beta^j=d_j$ and $ \sum_{j=1}^l
    d_j=d-1 $. The curve $ C_j $ for $ j=1,\dots,l $ passes through $r_j$ real
    and $s_j$ complex given points, where $ 2d_j+|\beta^j|-1=r_j+2s_j $. Note
    that (D), (E), and (F) are meant to be schematic pictures in which the
    thin and thick horizontal edges are just examples. The non-horizontal edges
    are always odd however.
  \end {enumerate}

  \begin {center} \input {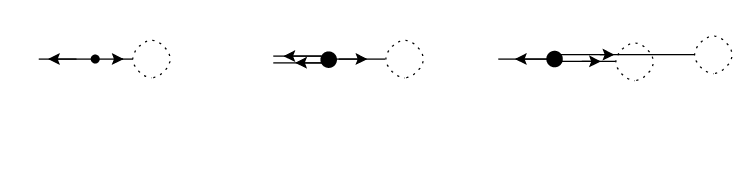} \end {center}

  \begin {center} \input {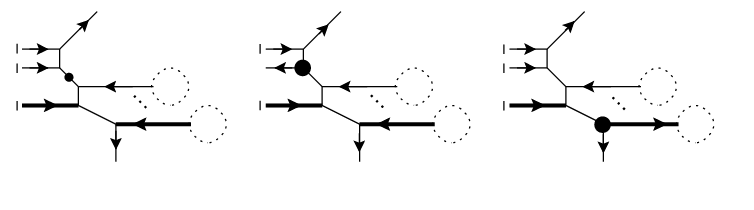} \end {center}
\end{proposition}

\begin{proof}
  \ref {prop-shape-a} Assume $C_0$ contains no bounded edge and $P_1$ is real.
  Then $C_0$ contains exactly one vertex, of type (1). Both adjacent edges are
  ends of $C_0$. Since $C$ is connected, one of the ends of $C_0$ results from
  cutting a bounded horizontal edge of $C$. Because of the balancing condition,
  it follows that the other end is of direction $(-k,0)$ for some $k>0$, which
  has to be odd since $P_1$ is of vertex type (1). Hence we are then in case
  (A).
  
  Assume now that $P_1$ is complex. Then $C_0$ consists of a vertex of type (5)
  or (6). At least one of the adjacent edges is of direction $(k,0)$ for some
  $k>0$ since it results from cutting a horizontal bounded edge. The other
  adjacent edges are ends of $C$. It follows from the balancing condition that
  all three adjacent edges are horizontal, and so we have type (B) or (C).
  Exactly one of the adjacent edges is even (and so vertex type (5) is
  impossible). In (A) and (B), we just cut one edge, so it follows that
  $\tilde{C}$ is irreducible and of the degree as claimed above. In (C), we cut
  two edges, so $\tilde{C}$ consists of two connected components $C_1$ and
  $C_2$. Ends of $C_1$ and $C_2$ are either ends of $C$ or the two cut edges.
  Denote their weights by $k_1$ resp.\ $k_2$, then it follows that $C_j$ is of
  a type $(\alpha^j,\beta^j)$ for $ j=1,2 $ with $\alpha^1+\alpha^2 =
  \alpha+e_{k_1}+e_{k_2}$ and $\beta^1+\beta^2=\beta-e_{k_1+k_2}$. If $
  2d_j+|\beta^j|-1 <r_j+2s_j$ for $j=1$ or $j=2$, then it follows that there is
  a connected component of $\Gamma$ minus the marked ends which does not
  contain a non-fixed unmarked end, a contradiction to lemma \ref{lem-general}.
  Thus we have $2d_j+|\beta^j|-1 \geq r_j+2s_j$, and since $
  2d_1+|\beta^1|-1+2d_2+|\beta^2|-1=2d+|\beta|-3=r+2(s-1)=r_1+2s_1+r_2+2s_2$ it
  follows that $2d_j+|\beta^j|-1=r_j+2s_j$ for $j=1,2$.

  \ref {prop-shape-b} Now assume that $C_0$ contains a bounded edge. By lemma
  \ref{lem-general}, each connected component of $C$ minus the marked points
  contains exactly one non-fixed unmarked end. If $P_1$ is real, removing the
  marked end $x_1$ satisfying $h(x_1)=P_1$ from $\Gamma$ produces $2$ connected
  components; if it is complex it produces $3$ connected components. It follows
  that $C_0$ contains at most $2$ non-fixed ends of $C$ if $P_1$ is real, or
  $3$ if $P_1$ is complex. Ends of $C_0$ are of direction $(k,0)$ for some $k$
  (resulting from cutting horizontal bounded edges of $C$) or ends of $C$. If
  $C_0$ contains a bounded edge then $ C_0 $ cannot lie entirely in a
  horizontal line, since otherwise the length of such a bounded edge could not
  be fixed by our conditions. It follows by the balancing condition that $C_0$
  must have ends of direction $(0,-1)$ and $(1,1)$, and in fact an equal number
  of them. But since ends of these directions are non-fixed and we have at most
  $3$ non-fixed ends of $C$ in $C_0$, we conclude that there is exactly one end
  of direction $(0,-1)$ and $(1,1)$ each. Since all other ends of $C_0$ are
  horizontal, it follows from the balancing condition that the directions of
  the bounded edges of $C_0$ are $\pm(a,1)$ for some $a$. In particular, they
  are all odd.

  If $P_1$ is real, $C_0$ cannot have more non-fixed ends of $C$ than the two
  ends of direction $(0,-1)$ and $(1,1)$. So then all left and right ends of
  $C_0$ are fixed, and we are in case (D). If $P_1$ is complex, there can be
  one non-fixed left end of $C_0$, which then has to be adjacent to $P_1$ as in
  case (E). Otherwise, $P_1$ has to be adjacent to a horizontal edge connecting
  $C_0$ with $\tilde{C}$. This is true because by the directions of the ends of
  $C_0$ and the balancing condition we can conclude that every vertex of $C_0$
  is adjacent to an edge of direction $(k,0)$ for some (positive or negative)
  $k$. Thus we are then in case (F).

  Assume we have to cut $l$ edges to produce $C_0$ and $\tilde{C}$, then
  $\tilde{C}$ consists of $l$ connected components. Each connected component is
  a curve of some type $(\alpha^j,\beta^j)$ with $I\alpha^j+I\beta^j=d_j$.
  It follows from the balancing condition that $\sum_{j=1}^l d_j=d-1$. The
  equations $ 2d_j+|\beta^j|-1=r_j+2s_j $ follow as in part \ref
  {prop-shape-a}.
\end{proof}

\begin{example}\label{ex-leftright}
  The picture shows an example of a curve $C$ decomposing into a floor $C_0$ of
  type (D) on the left and a reducible curve $\tilde{C}$ on the right. $C$ is
  of type $ ((3,1),(3,1)) $ passing through $r=7$ real and $s=8$ complex
  points satisfying $2d+|\beta|-1=20+4-1=23=r+2s$. We have chosen to move a
  real point to the left of the others.

  \begin {center} \input {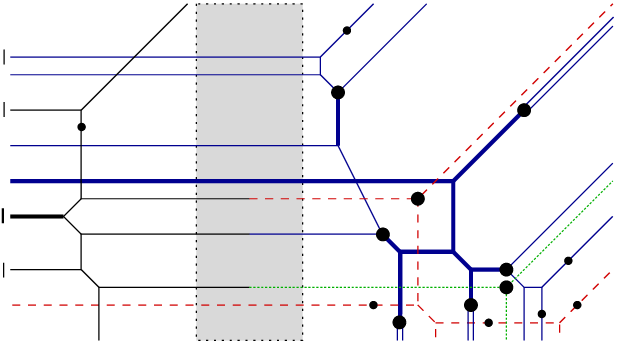} \end {center}

  The reducible curve $\tilde{C}$ consists of three connected components, $C_1$
  (green dotted), $C_2$ (red dashed) and $C_3$ (blue solid). $C_1$ is a curve
  of type $((0),(1))$ passing through $s_1=1$ complex points, satisfying
  $2d_1+|\beta^1|-1=2+1-1=2=r_1+2s_1$. $C_2$ is a curve of type $((0),(2))$
  passing through $r_2=3$ real and $s_2=1$ complex points satisfying
  $2d_2+|\beta^2|-1=4+2-1=5=r_2+2s_2$. $C_3$ is a curve of type $((1),(3,1))$
  passing through $r_3=3$ real and $s_3=6$ complex points satisfying
  $2d_3+|\beta^3|-1=12+4-1=15=r_3+2s_3$. We have $d_1+d_2+d_3=1+2+6=d-1 $.
  All three curves are connected to $C_0$ via a horizontal edge of weight $1$.
  We have $\beta=(3,1)=\beta^1+\beta^2+\beta^3-3e_1$ and $ \alpha^1+\alpha^2+
  \alpha^3=(1)<\alpha=(3,1)$.
\end{example}

Note that in the situation above there is always a unique possibility for $
C_0 $ once we are given the left and right ends of $ C_0 $ (together with their
position for fixed ends) as well as the position of $ P_1 $. Thus, to determine
$ N^d(\alpha,\beta,s) $, we just have to determine the different contributions
from all possibilities for $\tilde{C}$. This is the content of the following
theorem. 

\begin{theorem}[Caporaso-Harris formula for $N^d(\alpha,\beta,s)$]
    \label{thm-brocc-ch}
  The following two recursive formulas hold for the invariants $
  N^d(\alpha,\beta,s) $, where we use the notation $ r := 2d+|\beta|-2s-1 $
  (resp.\ $ r_j := 2d_j+|\beta_j|-2s_j-1 $ for all $j$) for the corresponding
  number of real markings in the invariant:
  \begin {enumerate}
  \item \label {thm-brocc-ch-a} (Moving a real point to the left) If $ r>0 $
    then
    \begin {align*} \qquad \quad
      N^d(\alpha,\beta,s) = \quad
         & \sum_{k \text { odd}}
	     N^d(\alpha+e_k,\beta-e_k,s) \tag A \\
        +& \sum \frac 1{l!} \,
             \binom {s}{s_1,\dots,s_l}
             \binom {r-1}{r_1,\dots,r_l}
             \binom {\alpha}{\alpha^1,\dots,\alpha^l} \,
             \prod_{m \text { even}} (-m)^{\alpha'_m} \,
             \prod_{\substack{j=1 \\ k_j \text { even}}}^l k_j \\
         & \qquad \cdot
             \prod_{j=1}^l \left( \beta_{k_j}^j \,
               N^{d_j}(\alpha^j,\beta^j,s_j) \right) \tag D
    \end{align*}
    where we set $ \alpha':= \alpha - \sum_{j=1}^l \alpha^j $, and where the
    sum in (D) runs over all $ l \ge 0 $ and all $ \alpha^j, \beta^j, k_j \ge
    1, d_j \ge 1, s_j \ge 0 \, $ for $ 1 \le j \le l $ satisfying $ \sum_j
    \alpha^j < \alpha $, $ \sum_j (\beta^j-e_{k_j}) = \beta $, $ \sum_j d_j =
    d-1 $, $ \sum_j s_j=s $.
  \item \label {thm-brocc-ch-b} (Moving a complex point to the left) If $ s>0 $
    then
    \begin{align*} \qquad
      N^d(\alpha,\beta,s) = \quad
         & \sum - \frac 12 \,
             N^d(\alpha+e_{k_1+k_2},\beta-e_{k_1}-e_{k_2},s-1) \tag B \\
        +& \sum \frac 12 \,
             \binom {s-1}{s_1,s_2}
             \binom {r}{r_1,r_2}
             \binom {\alpha}{\alpha^1,\alpha^2} \cdot
             \prod_{j=1}^2 N^{d_j}(\alpha^j+e_{k_j},\beta^j,s_j) \tag C \\
        +& \sum \frac 1{l!} \,
             \binom {s-1}{s_1,\dots,s_l}
             \binom {r}{r_1,\dots,r_l}
             \binom {\alpha}{\alpha^1,\dots,\alpha^l} \,
             M_k \,
             \prod_{m \text { even}} (-m)^{\alpha'_m} \,
             \prod_{\substack{j=1 \\ k_j \text { even}}}^l k_j \\
         & \qquad \cdot
             \prod_{j=1}^l \left( \beta_{k_j}^j \,
               N^{d_j}(\alpha^j,\beta^j,s_j) \right) \tag E \\
        +& \sum \frac 1{(l-1)!} \,
             \binom {s-1}{s_1,\dots,s_l}
             \binom {r}{r_1,\dots,r_l}
             \binom {\alpha}{\alpha^1,\dots,\alpha^l} \,
             \tilde M_{k_1} \,
             \prod_{m \text { even}} (-m)^{\alpha'_m} \,
             \prod_{\substack{j=2 \\ k_j \text{ even}}}^l k_j \\
         & \qquad \cdot
             N^{d_1}(\alpha^1+e_{k_1},\beta^1,s_1) \,
             \prod_{j=2}^l \left( \beta_{k_j}^j \,
               N^{d_j}(\alpha^j,\beta^j,s_j) \right) \tag F
    \end{align*}
    where as above $ \alpha':= \alpha - \sum_{j=1}^l \alpha^j $, and where the
    sums run over
    \begin {enumerate}
    \item [(B)] all $ k_1,k_2 \ge 1 $ such that at least one of them is odd;
    \item [(C)] all $ \alpha^j, \beta^j, k_j \ge 1, d_j \ge 1, s_j \ge 0 $ for
      $ j \in \{1,2\} $ such that at least one of $ k_1,k_2 $ is odd, $ \sum_j
      \alpha^j = \alpha $, $ \sum_j \beta^j = \beta-e_{k_1+k_2} $, $ \sum_j d_j
      = d $, $ \sum_j s_j=s-1 $.
    \item [(E)] all $ l \ge 0 $ and all $ \alpha^j, \beta^j, k \ge 1, k_j \ge
      1, d_j \ge 1, s_j \ge 0 $ for $ 1 \le j \le l $ such that $ \sum_j
      \alpha^j \leq \alpha $, $ \sum_j (\beta^j-e_{k_j}) = \beta-e_k $, $ \sum_j
      d_j=d-1 $, $ \sum_j s_j=s-1 $.
    \item [(F)] all $ l \ge 1 $ and all $ \alpha^j, \beta^j, k_j \ge
      1, d_j \ge 1, s_j \ge 0 $ for $ 1 \le j \le l $ such that $ \sum_j
      \alpha^j < \alpha $, $ \beta^1 + \sum_{j>1} (\beta^j-e_{k_j}) = \beta $,
      $ \sum_j d_j=d-1 $, $ \sum_j s_j=s-1 $.
    \end {enumerate}
    Here, the numbers $ M_k $ and $ \tilde M_k $ are defined by
      \[ \qquad \quad
         M_k = \begin {cases}
             k  & \text {if $k$ odd}, \\
             -1 & \text {if $k$ even}
         \end {cases}
         \qquad \text {and} \qquad
         \tilde M_k = \begin {cases}
             k  & \text {if $k$ odd}, \\
             1 & \text {if $k$ even}.
         \end {cases} \]
  \end {enumerate}
  Of course, for both equations it is assumed that the sums are taken only over
  choices of variables such that all occurring sequences have only non-negative
  entries and all relative broccoli invariants satisfy the dimension condition.
\end {theorem}

\begin{proof}
  As we have mentioned already we move one of the point conditions to the far
  left, so that each curve satisfying the conditions decomposes into a left
  part $ C_0 $ and a right part $ \tilde C $. Since we have studied the
  possibilities for $C_0$ and $\tilde{C}$ in proposition \ref{prop-shape}
  already it only remains to understand the different contributions to the
  relative broccoli invariant from each of these cases.

  \ref {thm-brocc-ch-a} The first formula arises from moving a real point to
  the left, so we have the cases (A) and (D).
  \begin {enumerate}
  \item [(A)] $C_0$ consists of one vertex of multiplicity $1$, and $\tilde{C}$
    has the same ends as $C$, with one odd non-fixed left end replaced by a
    fixed one. Thus we only have to sum over all possibilities of weights of
    this left end.
  \item [(D)] We have to sum over all possibilities for $\tilde{C}$ to split
    into $l$ connected components $C_1,\ldots,C_l$, where $C_j$ is of type
    $ (\alpha^j,\beta^j) $ with $I\alpha^j+I\beta^j=d_j$ and passes through
    $r_j$ real and $s_j$ complex points of $ P_2,\dots,P_{r+s} $. The right
    ends of $ C_0 $ are the gluing points for $ C_1,\dots,C_l $. They are fixed
    for $ C_0 $ and thus non-fixed for $ C_1,\dots,C_l $, i.e.\ they belong to
    $ \beta^1,\dots,\beta^l $. Let $k_j$ be the weight of the edge with
    which $C_0$ and $C_j$ are connected. Then we have $ \sum_{j=1}^l
    (\beta^j-e_{k_j}) = \beta $. Also, we have $ \sum_{j=1}^l\alpha^j<\alpha $,
    and $ \alpha'=\alpha-\sum_{j=1}^l\alpha^j $ is the sequence of fixed left
    ends adjacent to $C_0$. The multinomial coefficient $\binom
    {s}{s_1,\ldots,s_l}$ gives the number of possibilities how the $s$
    complex points of $ P_2,\dots,P_{r+s} $ can be distributed among the $C_j$.
    The second and third multinomial coefficient give the corresponding number
    for the real points and the fixed left ends, respectively.

    It remains to take care of different multiplicity factors. First of all
    note that every fixed left end adjacent to $C_0$ (described by $ \alpha' $)
    is not a fixed end of $\tilde{C}$ any more, so when counting the
    contribution from $\tilde{C}$ instead of $C$ we lose a factor of $i^{k-1}$
    for every such end of weight $k$ (remember that the weights of the ends of
    a curve $C$ enter into the multiplicity $m_C$, see definition \ref
    {def-vertex}). Also, each such fixed end is adjacent to a vertex of $C_0$
    whose multiplicity is $i^{k-1}\cdot k$ if $k$ is even and $i^{k-1}$
    if $k$ is odd. Thus, we lose a factor $ i^{2k-2} = (-1)^{k-1} = 1 $ if
    $k$ is odd, and $ k \cdot i^{2k-2} = k\cdot (-1)^{k-1} = -k $ if $k$ is
    even. Therefore we have to multiply by $ \prod_{m \text{ even}} (-m)^{
    \alpha_m'} $.

    Similarly, for $ j=1,\dots,l $ the end of weight $k_j$ with which $C_j$ is
    connected to $C_0$ yields a factor of $i^{k_j-1}$ in the multiplicity of
    $\tilde{C}$ that we do not need for $C$. The vertex of $C_0$ adjacent to
    such an edge has multiplicity $k_j\cdot i^{k_j-1}$ if $k_j$ is even, and
    $i^{k_j-1}$ if $k_j$ is odd. Thus we need to multiply by $\prod_{j: \,
    k_j \text{ even}}^l k_j $.

    The factors $\beta^j_{k_j}$ stand for the number of possibilities with
    which of the $\beta^j_{k_j}$ non-fixed ends of weight $k_j$ the component
    $C_j$ is connected to $C_0$. The factor $ \frac 1{l!} $ takes care of the
    overcounting due to the labeling of the components $ C_1,\dots,C_l $. As $
    C_0 $ has one end of direction $ (0,-1) $ and $ (1,1) $ each it is clear
    that we must have $ \sum_j d_j = d-1 $.
  \end {enumerate}

  \ref {thm-brocc-ch-b} In the second formula we move a complex point to the
  left, so we have four summands corresponding to the possibilities (B), (C),
  (E), and (F).
  \begin {enumerate}
  \item [(B)] We have to sum over all possibilities $k_1$ and $k_2$ for the
    weights of the two left ends which are adjacent to $P_1$. If we sum over
    all tuples $(k_1,k_2)$, we overcount by a factor of $2$ since these two
    weights are unordered. Therefore we multiply by $\frac{1}{2}$. For summands
    with $k_1=k_2$, the $\frac{1}{2}$ takes care of the factor of $\frac{1}{2}$
    in the multiplicity of the vertex of $C_0$ that we have to include when
    counting curves without labels at the unmarked ends (see remark \ref
    {rem-unlabeled}). We lose factors of $i^{k_1-1}$ and $i^{k_2-1}$ since
    these two ends are not ends of $\tilde{C}$, and we lose a factor of $
    i^{-1} $ for the vertex of $C_0$. Instead, we have a factor of
    $i^{k_1+k_2-1}$ for the end of $\tilde{C}$ with which it is glued to $C_0$.
    Thus, we have to multiply by $-1$.
  \item [(C)] In this case we have to sum over all choices of the connecting
    weights $ k_1 $ and $ k_2 $ (which are fixed ends for $C_1$ and $C_2$),
    degrees $ d_1 $ and $ d_2 $, and numbers $ s_1 $ and $ s_2 $ of complex
    markings on each component. The symmetry factor $ \frac 12 $ cancels the
    overcounting due to the labeling of the two components. The binomial
    factors count the possibilities how the complex and real points and the
    fixed ends can be distributed among $C_1$ and $C_2$. In $C_0$, we have the
    left end contributing $i^{k_1+k_2-1}$ and a vertex contributing $i^{-1}$,
    in $\tilde{C}$ we have instead the two ends contributing $i^{k_1-1}$ and
    $i^{k_2-1}$. So we do not need to multiply by a factor to take care of
    these multiplicities.
  \item [(E)] The terms are essentially as in (D) above, except that in
    addition we have to sum over all possibilities for the weight $k$ of the
    non-fixed left end adjacent to $P_1$. Also, this non-fixed end is not an
    end of any of the $C_j$, so the condition $ \sum_j (\beta^j-e_{k_j}) =
    \beta $ has to be changed to $ \sum_j (\beta^j-e_{k_j}) = \beta-e_k $. In
    addition to the factors of (D) we lose a factor of $i^{k-1}$ for the end,
    and of $i^{k-1}$ if $k$ is even and $k\cdot i^{k-1}$ if $k$ is odd for the
    vertex at $P_1$. So altogether we have to multiply by $ i^{2k-2} =
    (-1)^{k-1} = -1 $ if $k$ is even and by $k$ if $k$ is odd.
  \item [(F)] We get again a similar summand as in (E). However, here instead
    of summing over the possibilities for $k$ we now have to choose one of the
    $C_j$ --- call it $C_1$ --- which is adjacent to $P_1$. This component will
    then have an additional fixed end of weight $k_1$. So in the invariant for
    $ C_1 $ we have to replace $ \alpha^1 $ by $ \alpha^1+e_{k_1} $; at the
    same time however we do not have to multiply this invariant by $
    \beta^1_{k_1} $ as $C_1$ is connected to $C_0$ by a fixed end. The fixed
    end of weight $k_1$ of $C_1$ contributes a factor of $i^{k_1-1}$ to
    $\tilde{C}$. We lose the multiplicity of the vertex at $P_1$ which is $
    i^{k_1-1} $ if $k_1$ is even and $ k_1 \cdot i^{k_1-1}$ if $k_1$ is odd.
    Hence we have to multiply by $ \tilde M_{k_1} $. \qedhere
  \end {enumerate}
\end{proof}

Of course, theorem \ref{thm-brocc-ch} now gives recursive formulas for all
broccoli invariants $ N^d(\alpha,\beta,s) $, and thus in particular by remark
\ref{rem-relbroccoliwelschinger} also for the Welschinger numbers
$W_{\PP^2}(d,3d-2s-1,s)$.

\begin{example}[Relative broccoli invariants in degree 3]
  The following table shows all invariants $ N^d(\alpha,\beta,s) $ for $ d=3 $,
  as computed by theorem \ref{thm-brocc-ch}. The numbers in the last line are
  those that correspond to the degree-$3$ Welschinger invariants. The entries
  in the second last line are all $0$ in accordance with example
  \ref{ex-welschinger-even} \ref{ex-welschinger-even-b}.
    \[ \begin{array}{|c|ccccc|} \hline
         \alpha,\beta & s=0 & s=1 & s=2 & s=3 & s=4 \\ \hline
         (0,0,1),(0)  & 3   & 1   & -1  &     &     \\
         (0,1),(1)    & -12 & -8  & -4  & 0   &     \\
         (1,1),(0)    & -8  & -4  & 0   &     &     \\
         (1),(0,1)    & 0   & 0   & 0   & 0   &     \\
         (1),(2)      & 8   & 6   & 4   & 2   &     \\
         (2),(1)      & 8   & 6   & 4   & 2   &     \\
         (3),(0)      & 6   & 4   & 2   &     &     \\
         (0),(0,0,1)  & 3   & 1   & -1  & -3  &     \\
         (0),(1,1)    & 0   & 0   & 0   & 0   &     \\
         (0),(3)      & 8   & 6   & 4   & 2   & 0   \\ \hline
       \end{array} \]
  Inge Sandstad Skrondal implemented the formula of theorem \ref{thm-brocc-ch}
  in Java for his Master thesis \cite{Skr12} and got results up to degree $6$.
  They agree with the computations of absolute Welschinger numbers in
  \cite{ABM08}. He also found analogous formulas for $\PP^1 \times \PP^1$ and
  $\PP^2_k$ for $k \leq 2$.
\end{example}

  \bibliographystyle {amsalpha}
  \bibliography {references}

\end {document}